%% file: ms.tex
\documentclass{article}
\usepackage{fullpage,bbm}
\usepackage{amsmath,amssymb,amscd,amsxtra, amsfonts, amsthm, amsbsy,wasysym,graphics, graphicx}
\usepackage{epstopdf}
\usepackage{algpseudocode}
\usepackage{algorithm}
\usepackage{latexsym,enumerate}
\usepackage{float}
\usepackage{color}

\newcommand{\sd}{\Sigma\Delta}
\newcommand{\bitem}{\begin{itemize}}
\newcommand{\eitem}{\end{itemize}}
\newcommand{\benum}{\begin{enumerate}}
\newcommand{\eenum}{\end{enumerate}}
\newcommand{\beq}{\begin{equation}}
\newcommand{\eeq}{\end{equation}}
\newcommand{\bal}{\begin{align*}}
\newcommand{\eal}{\end{align*}}
\newcommand{\bm}{\left[\begin{matrix}}
\newcommand{\enm}{\end{matrix}\right]}
\newcommand{\wtl}{\widetilde}
\newtheorem{theorem}{Theorem}[section]
\newtheorem{pro}[theorem]{Proposition}

\newtheorem{lm}[theorem]{Lemma}

\newtheorem{corollary}[theorem]{Corollary}
\theoremstyle{definition}
\newtheorem{definition}[theorem]{Definition}

\newtheorem{ex}{Example}
\usepackage{xcolor}

\newcommand{\be}{\begin{equation}}
\newcommand{\nde}{\end{equation}}
\theoremstyle{remark}

\numberwithin{equation}{section}

\newcommand{\R}{\mathbb R}

\def\x{{\mathbf x}}
\def\y{{\mathbf y}}
\def\q{{\mathbf q}}
\def\v{{\mathbf v}}
\def\u{{\mathbf u}}

\begin{document}
\author{Theodore Faust\thanks{University of California Los Angeles, Department of Mathematics, \textsc{tfaust@math.ucla.edu}.}, 
Mark Iwen\thanks{Michigan State University, Department of Mathematics, and the Department of Computational Mathematics, Science and Engineering (CMSE), \textsc{markiwen@math.msu.edu}. Supported in part by NSF DMS 1912706.}, %
Rayan Saab\thanks{University of California San Diego, Department of Mathematics, and the Hal{\i}c{\i}o{g}lu Data Science Institute, \textsc{rsaab@ucsd.edu}. Supported in part by NSF DMS 2012546.},
Rongrong Wang\thanks{Michigan State University, Department of Computational Mathematics, Science and Engineering (CMSE), and the Department of Mathematics, \textsc{wangron6@msu.edu}. Supported in part by NSF CCF-1909523.}}
\title{On the $\ell^{\infty}$-norms of the Singular Vectors of Arbitrary Powers of a Difference Matrix with Applications to Sigma-Delta Quantization}
%\title{Singular vectors of finite difference matrices with arbitrary order}
\maketitle
\begin{abstract}
Let $\| A \|_{\max} := \max_{i,j} |A_{i,j}|$ denote the maximum magnitude of entries of a given matrix $A$.  In this paper we show that 
\[
\max \left\{ \|U_r \|_{\max},\|V_r\|_{\max} \right\} \le \frac{(Cr)^{6r}}{\sqrt{N}} ,  
\]   
where $U_r$ and $V_r$ are the matrices whose columns are, respectively, the left and right singular vectors of the $r$-th order finite difference matrix $D^{r}$ with $r \geq 2$, and where $D$ is the $N\times N$ finite difference matrix with $1$ on the diagonal, $-1$ on the sub-diagonal, and $0$ elsewhere.  Here $C$ is a universal constant that is independent of both $N$ and $r$.  Among other things, this establishes that both the right and left singular vectors of such finite difference matrices are {\it Bounded Orthonormal Systems (BOSs)} with known upper bounds on their BOS constants, objects of general interest in classical compressive sensing theory.  Such finite difference matrices are also fundamental to standard $r^{\rm th}$ order Sigma-Delta quantization schemes more specifically, and as a result the new bounds provided herein on the maximum $\ell^{\infty}$-norms of their $\ell^2$-normalized singular vectors allow for several previous Sigma-Delta quantization results to be generalized and improved.
\end{abstract}

%%%%%%%%%%%%%%%%%%%%%%%%%%%%%%%%%%%%%%%%%%%%%%%%%%%%%%%%%%%%%%
%%%%%%%%%%%%%%%%%%%%%%%%%%%%%%%%%%%%%%%%%%%%%%%%%%%%%%%%%%%%%%
\section{Introduction}
\label{sec:Intro}

\input{Intro.tex}
\section{Proof of the Main Result (Theorem~\ref{thm:main0})}
\label{sec:Proof}

Below we will denote the set $\{ 1, 2, \dots, n \} \subset \mathbbm{N}$ by $[n]$.  For any matrix $A \in \mathbbm{R}^{m \times N}$ we will denote the $j^{\rm th}$ column of $A$ by ${\bf a}_j \in \mathbbm{R}^{m}$.  The transpose of a matrix, $A \in \mathbbm{R}^{m \times N}$, will be denoted by $A^{\rm T} \in \mathbbm{R}^{N \times m}$, and the singular values of any matrix $A \in \mathbbm{R}^{m \times N}$ will always be ordered  as $\sigma_1(A) \geq \sigma_2(A) \geq \dots \geq \sigma_{\min(m,N)}(A) \geq 0.$  We will denote the standard indicator function by
$$\delta_{i,j} := \left\{ \begin{array}{ll} 1 & \textrm{if}~i=j\\ 0 & \textrm{if}~i \neq j\end{array} \right.,$$
for $i,j \in \mathbbm{N}$.  Given a matrix $A$ with a singular value decomposition $A = U \Sigma V^*$, we use $\mathbf{u}_j$ (resp. $\v_j$) to denote the columns of $U$ (resp. $V$).

To begin the proof, it is straightforward to verify that with reversed row and column orders, $D^r(D^r)^T$ coincides with $(D^r)^TD^r$.  That is, the $(i,j)$th element of $D^r(D^r)^T$ is equal to the $(N-i,N-j)$th element of $(D^r)^TD^r$ for all $r \geq 1$. Then, since the eigenvectors of $D^r(D^r)^T$ and $(D^r)^TD^r$ are the left and right singular vectors of $D^r$, respectively, this then implies that the left singular vectors of $D^r$ are just the right singular vectors with reversed entries.  As a result, we have the following lemma, which we prove in detail in Appendix~\ref{sec:Proofoflem:Rightinf=Leftinf}.

\begin{lm}
Suppose that $D^r$ has singular value decomposition $D^r = U \Sigma V^*$ for $r \geq 1$.  Then,
$\|\mathbf{u}_j\|_{\infty}= \|\mathbf{v}_j\|_{\infty}$ for all $j \in [N]$.
\label{lem:Rightinf=Leftinf} 
\end{lm}

In light of Lemma~\ref{lem:Rightinf=Leftinf} it suffices to prove the following result bounding just the $\ell^\infty$-norms of the right singular vectors of $D^r$ in order to obtain a proof of Theorem~\ref{thm:main0}.

\begin{theorem}\label{thm:main}
Suppose that $r \ge 2$, and let $\sigma_j := \sigma_j \left( D^r \right)$ have associated right singular vector $\mathbf{v}_j \in \mathbbm{R}^N$ for all $j \in [N]$. There exists absolute constants $C,C_3 > 0$ such that if $N \ge C_3^r$, we have $||\mathbf{v}_j||_\infty \le \frac{(Cr)^{6r}}{\sqrt{N}}$ for all $j \in [N]$.
\end{theorem}

The proof of Theorem \ref{thm:main} will be broken up into cases depending on the size of $\sigma_j$, the $j$th singular value of $D^r$. Thus, we begin by providing bounds for each $\sigma_j$.
\begin{lm}\label{lm:sigD}
Let $\sigma_j := \sigma_j \left( D^r \right)$ be the $j^{th}$ singular value of $D^r \in \mathbbm{R}^{N \times N}$ where $j \in [N]$, and $\sigma_1 \geq \sigma_2 \geq \cdots \geq \sigma_N$. Then, 
$$0<\sigma_j < \left(2\cos\left(\frac{\pi}{2N+1}\right)\right)^{r}$$
holds for all $j \in [N]$.
Moreover, there exist absolute constants $c,C \in \mathbbm{R}^+$ such that
$$r^{-r}c \left(\frac{j}{N}\right)^{r}\leq \sigma_{N-j+1} \leq r^{r} C\left(\frac{j}{N}\right)^{r}$$
also holds for all $j \in [N]$.
\end{lm}
\begin{proof}
Since $D$ is of full rank we have that $\sigma_j>0$ for all $j \in [N]$. In addition $\|D \mathbf{v}\|_2 < 2{\cos\left(\frac{\pi}{2N+1}\right)} \| \mathbf{v}\|_2$ holds for all $\mathbf{v} \in \mathbbm{C}^N$ (see, e.g., \cite{Iwen2013}), which implies that $\|D^r \mathbf{v} \|_2 < \left(2{\cos\left(\frac{\pi}{2N+1}\right)}\right)^r \| \mathbf{v}\|_2 $  for all $\mathbf{v} \in \mathbbm{C}^N$, and hence that $\sigma_1< \left(2{\cos\left(\frac{\pi}{2N+1}\right)}\right)^r$. The second item is a direct consequence of Proposition 2.2 in \cite{KSY2014}.\end{proof}

We see that this result implies that $\sigma_j \in (0,2^r)$. Going forward we will prove Theorem~\ref{thm:main} by bounding $||\mathbf{v}_j||_\infty$ in two separate cases:  the case where $\sigma_j$ is ``small'' (namely $0 < \sigma_j^{1/r}<\frac{C_2r^6}{N}$ for a constant $C_2>0$), and the case where $\sigma_j$ is ``large'' ($2 > \sigma_j^{1/r} \ge \frac{C_2r^6}{N}$). In Section \ref{sec:Proof} below we prove the result for the case of ``small'' $\sigma_j$. This proof is a fairly straightforward application of results about $D$ together with a simple lemma concerning discrete dynamical systems. We also state the result for ``large'' $\sigma_j$ and then formally prove Theorem \ref{thm:main} given that the stated result holds.  The remainder of paper is then dedicated to proving that the stated result for the case of  ``large'' $\sigma_j$ actually does indeed hold.

To begin the proof of the ``large'' $\sigma_j$ case (i.e., Lemma \ref{lm:large} below), we first find a formula for the right singular vectors of $D^r$ in Section \ref{sec:EqnsForEigenvector}. To achieve this goal, we extend each singular vector $\mathbf{v}_j$ to an infinite sequence $\tilde{\mathbf{v}}_j$, and then use techniques from the solution of difference equations to find a formula for each entry of $\mathbf{v}_j$. In particular, we are able to write each $\mathbf{v}_j$ in terms of the roots $\rho_{k, \ell}$ of a characteristic polynomial $p(x)$ (which differs for each $j$) in the following way: 
\begin{equation}\label{eq:vjintro}
\left(\v_{j}\right)_i = \sum^{1}_{\ell=0} \sum^{r-1}_{k = 0} c_{k,\ell} \cdot \rho^i_{k,\ell}.
\end{equation}
The rest of the section is then devoted to proving results about the roots $\rho_{k, \ell}$ which then ultimately allow us to bound $||\mathbf{v}_j||_\infty$.

Finally, in Section \ref{sec:proof_BIG_lemma}, we prove the main result in the case that $\sigma_j$ is ``large'', and so complete the proof of Theorem~\ref{thm:main} as a result. To do this, we first seek to find a bound on the constants $c_{k,\ell}$ from the above expression \eqref{eq:vjintro} corresponding to roots $\rho_{k,\ell}$ with $|\rho_{k,\ell}| \ge 1$. This proof is rather involved, and so is contained in Appendix \ref{sec:Proofoflem:coeff}. After this, we use this bound and the properties of the infinite sequence $\tilde{\mathbf{v}}_j$ to bound $c_{k,\ell}$ in the case that $|\rho_{k,\ell}| < 1$, which gives us a bound on $c_{k,\ell}$ for all $k,\ell$. We then use those bounds together with the properties of each infinite sequence $\tilde{\mathbf{v}}_j$ to bound $c_{k,\ell}\rho_{k,\ell}^{N+1-r}$ for all $k,\ell$. Combining these bounds, we are then able to prove the main result in the case that $\sigma_j$ is ``large'', thereby completing the proof of Lemma \ref{lm:large} below (and, therefore, proving Theorem \ref{thm:main} as well).

We next begin by proving the result in the ``small'' $\sigma_j$ case. To do this, we will utilize the following general result concerning the $\ell^{\infty}$-norms of the $\ell^2$-normalized right singular vectors of an arbitrary matrix power $A^r \in \mathbbm{C}^{N \times N}$.  More specifically, the following lemma can be used to show that the right singular vectors of $A^r$ associated with its smallest singular values will always be ``flat'' (i.e., have $\ell^{\infty}$-norms on the order of $\sim \frac{1}{\sqrt{N}}$) when the rows of $A^{-1}$ are all sufficiently small in $\ell^2$-norm.

\begin{lm}
Let $r \in \mathbbm{Z}^+$, $A \in \mathbbm{C}^{N \times N}$, and $A^r$ have the singular value decomposition $A^r = U \Sigma V$.  If $\left( \sigma_j\left( A^r \right) \right)^{1/r} \leq \alpha \cdot \sigma_N \left( A \right)$ holds for some $j \in [N]$ and $\alpha \in \mathbbm{R}^+$, then $$\| \mathbf{v}_j \|_{\infty} \leq \alpha^r \sigma_N \left( A \right) \cdot \max_{k \in [N]} \| (A^{-1})^* \mathbf{e}_k \|_{2}.$$
\label{lem:DynamicalSys}
\end{lm}

\begin{proof}
Consider the discrete dynamical system defined by $\Phi_{\mathbf{v}_j}(k + 1) := \left( \sigma_j \left( A^r \right) \right)^{-1/r} A \Phi_{\mathbf{v}_j}(k)$ for all $k \in \mathbbm{Z}^+$ with $\Phi_{\mathbf{v}_j}(0) := \mathbf{v}_j$.  It is not difficult to see that both
\begin{itemize}
\item $\| \Phi_{\mathbf{v}_j}(r) \|_2 = \| \mathbf{u}_j \|_2 = 1$, and  
\item $\| \Phi_{\mathbf{v}_j}(k) \|_2 = \| A^{-1} \left( \sigma_j \left( A^r \right) \right)^{1/r} \Phi_{\mathbf{v}_j}(k+1) \|_2 \leq \frac{ \left( \sigma_j \left( A^r \right) \right)^{1/r}}{\sigma_N \left( A \right)} \| \Phi_{\mathbf{v}_j}(k+1) \|_2 \leq \alpha \| \Phi_{\mathbf{v}_j}(k+1) \|_2$
\end{itemize}
hold for all $k \in \mathbbm{Z}^+$ since $\left( \sigma_j\left( A^r \right) \right)^{1/r} \leq \alpha \cdot \sigma_N \left( A \right)$.  As a consequence, $\| \Phi_{\mathbf{v}_j}(1) \|_2 \leq \alpha^{r-1}$ must also hold.

Continuing, we can now see that
\begin{align*}
\| \mathbf{v}_j \|_{\infty} &= \| A^{-1} \left( \sigma_j \left( A^r \right) \right)^{1/r}  \Phi_{\mathbf{v}_j}(1) \|_{\infty} = \left( \sigma_j \left( A^r \right) \right)^{1/r} \max_{k \in [N]} \left| \langle A^{-1} \Phi_{\mathbf{v}_j}(1), \mathbf{e}_k \rangle \right|\\
& = \left( \sigma_j \left( A^r \right) \right)^{1/r} \max_{k \in [N]} \left| \langle \Phi_{\mathbf{v}_j}(1), (A^{-1})^* \mathbf{e}_k \rangle \right|\\
& \leq \left( \sigma_j \left( A^r \right) \right)^{1/r} \max_{k \in [N]} \| \Phi_{\mathbf{v}_j}(1) \|_2 \| (A^{-1})^* \mathbf{e}_k \|_{2}\\
& \leq \alpha^r \sigma_N \left( A \right) \cdot \max_{k \in [N]} \| (A^{-1})^* \mathbf{e}_k \|_{2},
\end{align*}
where the last inequality uses both our assumed upper bound on $\left( \sigma_j\left( A^r \right) \right)^{1/r}$, and the fact that $\| \Phi_{\mathbf{v}_j}(1) \|_2 \leq \alpha^{r-1}$.
\end{proof}

With Lemma~\ref{lem:DynamicalSys} in hand we may immediately obtain the following result, which we will use to bound $\mathbf{v}_j$ in the case that $\sigma_j$ is ``small''. 

\begin{corollary}\label{cor:smallvjbound}
Let $\sigma_j := \sigma_j \left( D^r \right)$ have associated right singular vector $\mathbf{v}_j \in \mathbbm{R}^N$ for all $j \in [N]$.  If $\sigma^{1/r}_j \leq \frac{\alpha}{N}$ for some $\alpha \in \mathbbm{R}^+$, then there exists an absolute constant $c_2 \in \mathbbm{R}^+$ such that $\| \mathbf{v}_j \|_{\infty} \leq \frac{\left( c_2 \alpha \right)^r}{\sqrt{N}}$.
\end{corollary}

\begin{proof}
We apply Lemma~\ref{lem:DynamicalSys} with $A = D$.  Note that
\begin{equation*}
D^{-1}_{i,j} := \left\{ \begin{array}{ll} 1 & \textrm{if}~i \leq j\\ 0 & \textrm{otherwise}\end{array} \right..
\end{equation*}
Thus, $\max_{k \in [N]} \| (D^{-1})^* \mathbf{e}_k \|_{2} = \sqrt{N}$.  Furthermore, Lemma~\ref{lm:sigD} (with $r=1$) tells us that $\frac{c}{N}\leq \sigma_N \left( D \right) \leq \frac{C}{N}$ so that $\sigma^{1/r}_j \leq \frac{\alpha}{c} \cdot \frac{c}{N} \leq \frac{\alpha}{c} \cdot \sigma_N \left( D \right)$.  Thus, Lemma~\ref{lem:DynamicalSys} allows us to conclude that
$$\| \mathbf{v}_j \|_{\infty} \leq \left( \frac{\alpha}{c} \right)^r \sigma_N \left( D \right) \sqrt{N} \leq \left( \frac{\alpha}{c} \right)^r \frac{C}{\sqrt{N}} \leq \frac{\left( c_2 \alpha \right)^r}{\sqrt{N}}$$
where $c_2 := \frac{\max\{1,C\}}{c} \geq \frac{C^{1/r}}{c} $ for  $c,C \in \mathbbm{R}^+$ as in Lemma~\ref{lm:sigD}.
\end{proof}

We see from Corollary \ref{cor:smallvjbound} that we can bound the $\ell^{\infty}$-norm of $\mathbf{v}_j$ in the case that $\sigma_j$ is ``small''. As discussed previously, the remaining sections of the paper will be devoted to proving the following main result, which bounds $\mathbf{v}_j$ in the case that $\sigma_j$ is ``large'':
\begin{lm}\label{lm:large} 
Let $\sigma_j := \sigma_j \left( D^r \right)$ have associated right singular vector $\mathbf{v}_j \in \mathbbm{R}^N$ for all $j \in [N]$. There exist absolute universal constants $C_0, C_2, C_3 \in \mathbbm{R}^+$ such that for all $r \geq 2$ and $N \geq C_3^r$,
$\|\v_j\|_{\infty}\leq \frac{(C_0r)^{4r-5}}{\sqrt{N}}$ holds for all $j \in [N]$ with $\sigma_j^{1/r}  \geq \frac{C_2 r^6}{N}$.
% There exist constants absolute $C_0$, $C_1$ and $C_2$ independent of $N$, such that provided $s_i^{1/r} >\frac{C_0 (C_1^{r}+r^4) r^2\log r }{N}$, we have
%$\|V_i\|_{\infty}\leq C_2(C_1 r)^{4r} /\sqrt N$. 
%$|\alpha_j| \leq c \sqrt{\left|\frac{1-\lambda_{j}^2}{1-\lambda_j^{2N} }\right|}$
\end{lm}

\begin{proof}
This proof is quite involved.  See Sections~\ref{sec:EqnsForEigenvector}~and~\ref{sec:proof_BIG_lemma} below.
\end{proof}

Using Corollary \ref{cor:smallvjbound} and Lemma \ref{lm:large} we can now prove Theorem \ref{thm:main}, thereby establishing our main result. 

\begin{proof}[\textbf{Proof of Theorem \ref{thm:main}}]
Suppose that $N \ge C_3^r$ for $C_3$ as in Lemma \ref{lm:large}. Then, if $\sigma_j^{1/r} \ge \frac{C_2 r^6}{N}$ (for $C_2$ as in Lemma \ref{lm:large}), Lemma \ref{lm:large} implies $\|\v_j\|_{\infty}\leq \frac{(C_0r)^{4r-5}}{\sqrt{N}}$.  If $\sigma_j^{1/r} <  \frac{C_2 r^6}{N}$, then by setting $\alpha =  C_2 r^6$ in Corollary \ref{cor:smallvjbound}, we have 
\[
\|\v_j\|_{\infty} \leq \frac{(C'r)^{6r}}{\sqrt{N}}.
\]
where $C'>0$ is an absolute constant chosen such that $(c_2C_2r^6)^r \le (C'r)^{6r}$ for all $r \ge 2$. 
Thus, choosing $C>0$ an absolute constant such that $(Cr)^{6r} \ge \max\{(C_0r)^{4r-5}, (C'r)^{6r}\}$ for all $r \ge 2$, we have 
\[
\|\v_j\|_{\infty} \leq \frac{(Cr)^{6r}}{\sqrt{N}}.
\]
for all $j \in [N]$.
\end{proof}

The remaining sections of the paper are dedicated to proving the result Lemma \ref{lm:large} which will allow us to complete the proof of Theorem \ref{thm:main} (and therefore, to complete the proof of our main result).

%%%%%%%%%%%%%%%%%%%%%%%%%%%%%%%%%%%%%%%%%%%%%%%%%%%%%%%%%%%%%%
%%%%%%%%%%%%%%%%%%%%%%%%%%%%%%%%%%%%%%%%%%%%%%%%%%%%%%%%%%%%%%
\section{Toward the Proof of Lemma~\ref{lm:large}:  A Formula for the Right Singular Vectors of $D^r$ }
\label{sec:EqnsForEigenvector}

Before we can prove Lemma~\ref{lm:large} we will need some basic facts about the structure of the right singular vectors of $D^r \in \mathbbm{Z}^{N \times N}$ for any $r \in \mathbbm{Z}^+$.  Note that these singular vectors will be identical to the eigenvectors of the related symmetric matrix $(D^r)^TD^r$.  As a result, the remainder of this section will be devoted to studying the structure of $(D^r)^TD^r$.  %As one might suspect, we will begin by first established a general formula for the entries of $(D^r)^TD^r$. 
The next lemma begins our study of $(D^{r})^T D^r$ by establishing a general formula for its entries, which turn out to be closely related to the alternating binomial coefficients.   
   
\begin{lm}
Let $r, N \in \mathbbm{Z}^+$ be such that $r < N / 2$.  All the entries of $(D^r)^TD^r \in \mathbbm{Z}^{N \times N}$ are given by $$((D^r)^TD^r)_{j-m,j} = \left\{
\begin{array}{ll}
        \displaystyle    (-1)^m  {2r \choose r-m} & \text{ if } j-m \leq N -r , ~m \in \{ 0, 1, \cdots, r \} \\
& \\
      
      \displaystyle (-1)^m \sum_{l=0}^{N-j} {r \choose l+m}{r \choose l} &  \text{ if } j-m > N-r,~m \in \{ 0, 1, \cdots, r\} \\
      \hspace{.5in}0 & \text{ if } m > r
\end{array} 
\right.,$$
combined with the fact that $(D^r)^TD^r$ is symmetric.
\label{lem:EntriesofDrTDr}
\end{lm}

\begin{proof} See Appendix~\ref{sec:ProofofEntriesofDrTDr}. \end{proof}

With Lemma~\ref{lem:EntriesofDrTDr} in hand we are now ready to study the eigenvectors of $\left(D^r \right)^{\rm T} D^r$.  Let $\lambda \in \mathbbm{R}^+$ be an eigenvalue with associated eigenvector ${\bf v}$.  Note that we know $\lambda \in (0,4^r)$ from Lemma~\ref{lm:sigD}. Considering the equation $\left(D^r \right)^{\rm T} D^r {\bf v} = \lambda {\bf v}$ in light of Lemma~\ref{lem:EntriesofDrTDr}, we can see that 
\begin{equation}
\sum^{2r}_{k=0} (-1)^{k+r}{2r \choose k} v_{i-r+k} = \lambda v_{i}
\label{eq:RecRel1}
\end{equation}
holds for all $N-r \geq i > r$.  Our strategy going forward will be to extend ${\bf v}$ to an infinite sequence $\tilde{\bf v}$ which satisfies the simple recurrence \eqref{eq:RecRel1} for all $i \in \mathbbm{Z}$, instead of just for $i \in (r,N-r]$.  That is, we want to construct an infinite sequence $\tilde{\bf v}$ such that both 
$$\tilde{v}_{i} = v_i~\textrm{for all}~i \in [1,N],$$
and
\begin{equation}
\sum^{2r}_{k=0} (-1)^{k+r}{2r \choose k} \tilde{v}_{i-r+k} = \lambda \tilde{v}_{i} ~\textrm{for all}~i \in \mathbbm{Z},
\label{eq:RecRel}
\end{equation}
hold.  Once we have managed to complete this task we will then be able to use standard techniques for the solution of difference equations (see, e.g., \cite{DiffEqu_Goldberg}) in order to construct a simple formula for every entry of $\tilde{\bf v}$.  This same formula will then also generate every entry of ${\bf v}$.  Finding such a formula is the ultimate goal of this section.

\subsection{Extending ${\bf v}$ to a Sequence $\tilde{\bf v}$ that Satisfies \eqref{eq:RecRel}}

We will extend ${\bf v}$ to an infinite sequence $\tilde{\bf v}$ as follows.  Let $\tilde{\bf v}$ be the sequence of real numbers whose entries $i \in [1-r,N+r]$ are given by
\begin{equation}
\tilde{v}_{i} = \left\{ \begin{array}{ll} 0 & \textrm{if}~ 1-r \leq i \leq 0\\ 
v_i & \textrm{if}~ 1 \leq i \leq N \\
v'_i & \textrm{if}~ N+1 \leq i \leq N+r \end{array} \right. ,
\label{eq:RRecRel1}
\end{equation}
where the $v'_{N+1}, \dots, v'_{N+r} \in \mathbbm{R}$ above are chosen so that 
\begin{equation}
\sum^r_{k=0} (-1)^{k}{r \choose k} \tilde{v}_{i-k} = 0
\label{eq:RRecRel1BC}
\end{equation}
holds for all $i \in [N+1,N+r]$.  Note that these $\tilde{v}_i = v'_i$ are uniquely defined by ${\bf v}$ together with \eqref{eq:RRecRel1BC} for all $i \in [N+1,N+r]$.  Having extended ${\bf v}$ to the larger index set $[1-r,N+r]$ in this fashion, we may now finish extending ${\bf v}$ to all of $\mathbbm{Z}$ by inductively setting
\begin{equation}
\tilde{v}_{i} = \left\{ \begin{array}{ll} (-1)^r \lambda \tilde{v}_{i+r} - \sum^{2r}_{k=1} (-1)^k {2r \choose k} \tilde{v}_{i+k} & \textrm{if}~ i \leq -r \\ (-1)^r \lambda \tilde{v}_{i-r} - \sum^{2r-1}_{k=0} (-1)^k {2r \choose k} \tilde{v}_{i-2r+k} & \textrm{if}~ i > N+r \end{array} \right..
\label{eq:RRecRel2}
\end{equation}

\begin{lm}
Equations $\eqref{eq:RRecRel1}-\eqref{eq:RRecRel2}$ imply that \eqref{eq:RecRel} holds.
\label{lem:proofofeq:RecRel1}
\end{lm}

\begin{proof}
See Appendix~\ref{sec:Proofoflem:proofofeq:RecRel1}.
\end{proof}

In the next subsection we will solve \eqref{eq:RecRel} via its characteristic polynomial.
Note that real solutions are guaranteed to exist for \eqref{eq:RecRel} whenever $\lambda$ is an eigenvalue of $\left(D^r \right)^{\rm T} D^r$, and we can always find them via the approach below (see, e.g., \cite{DiffEqu_Goldberg}).

\subsection{Solving the Related Recurrence Relation for $\tilde{\bf v}$}

Before we can write a formula for $\tilde{\bf v}$ we must first find the roots of the characteristic polynomial of \eqref{eq:RecRel} (for simplicity, we multiply each side of \eqref{eq:RecRel} by $(-1)^r$)
\begin{equation}
p(x) ~=~ \sum^{2r}_{k=0} {2r \choose k} (-x)^k - (-1)^{r} \lambda x^r~=~ (1-x)^{2r} - (-1)^{r} \lambda x^r.
\label{eq:charPoly}
\end{equation}
By considering \eqref{eq:charPoly} when $p(x) = 0$ it is not difficult to see that
\begin{equation}
p(x) ~=~ \prod^{r-1}_{k=0} \left( -(1-x)^2 - \lambda^{1/r} \mathbbm{e}^{2 k \pi \mathbbm{i}/r} x \right).
\label{eq:charPolyFact}
\end{equation}
 Examining \eqref{eq:charPolyFact}, one can now easily deduce the following lemma concerning the roots of the characteristic polynomial $p(x)$.
\begin{lm}\label{lem:charPolyProps1}
The roots of $p(x)$ are given by 
\begin{equation}
\rho_{k,0} := \frac{2 - \lambda^{1/r} \mathbbm{e}^{2 k \pi \mathbbm{i}/r} + \sqrt{\lambda^{2/r} \mathbbm{e}^{4 k \pi \mathbbm{i}/r}-4\lambda^{1/r} \mathbbm{e}^{2 k \pi \mathbbm{i}/r}}}{2}
\label{eq:root1}
\end{equation}
and 
\begin{equation}
\rho_{k,1} := \frac{2 - \lambda^{1/r} \mathbbm{e}^{2 k \pi \mathbbm{i}/r} - \sqrt{\lambda^{2/r} \mathbbm{e}^{4 k \pi \mathbbm{i}/r}-4\lambda^{1/r} \mathbbm{e}^{2 k \pi \mathbbm{i}/r}}}{2}
\label{eq:root2}
\end{equation}
for $k \in \{ 0, 1, \dots, r-1 \}$. Moreover, it is also not difficult to see that both 
$$\rho_{k,0} =\rho^{-1}_{k,1} \quad \text{ for all } \quad k \in \{ 0, 1, \dots, r-1 \}$$
and
$$\rho_{k,j} = \overline{\rho_{r-k,j}} \quad\text{ for all } \quad k \in \{ 1, \dots, r-1 \}$$
are true.  That is, both the multiplicative inverse and complex conjugate of every root are also a root. 
\end{lm}
\begin{proof}
The lemma can be directly verified by substitutions using \eqref{eq:charPolyFact} -- \eqref{eq:root2}, coupled with the fact that $\lambda \in (0,4^r)$ by Lemma~\ref{lm:sigD}.  \end{proof}

In the following lemmas we will establish several other important properties of the roots of $p(x)$, including their uniqueness for all $r \geq 2$. These properties will be useful later.  In particular, the fact that each root of $p(x)$ is unique (i.e., has multiplicity one) will be crucial to our ability to write down a simple formula for each entry of $\tilde{\bf v}$, and therefore, will also be crucial to our discovery of a compact formula for each eigenvector of $\left(D^r \right)^{\rm T} D^r$.
%%%
\begin{lm}\label{lem:BasicRootsEqns}
Let $\rho$ be any root of $p(x)$ as given in \eqref{eq:root1} and \eqref{eq:root2}. Then, for $r\geq 2$ 
\begin{equation} \label{eq:rhoBound2}
(1+\sqrt{2})^{-2} \leq |\rho| \leq (1+\sqrt{2})^2
\end{equation}
and
\begin{equation} \label{eq:rhoBound1}
(1+\sqrt{2})^{-1} \lambda^{\frac{1}{2r}} \leq | \rho-1| \leq (1+\sqrt{2}) \lambda^{\frac{1}{2r}}.
\end{equation}
Furthermore, using \eqref{eq:root1} and \eqref{eq:root2} we see that $\rho_{0,0} \neq \rho_{0,1}$,
$$|\rho_{0,0}| =|\rho_{0,1}|=1,$$ and 
$$|\rho_{k,j}| \neq 1~\text{holds unless}~ k = 0.$$
\end{lm}

\begin{proof}
We will begin with \eqref{eq:rhoBound1} and \eqref{eq:rhoBound2}. 
Examining \eqref{eq:charPolyFact}, we have that  
\begin{equation}
c_k=c_k(\rho):= \sqrt{\rho}-\frac{1}{\sqrt{\rho}} = \pm i \lambda^{\frac{1}{2r}}\mathbbm{e}^{k\pi i/r}.\label{eq:ck}
\end{equation}
Recalling again that $0<\lambda<4^{r}$ by Lemma~\ref{lm:sigD}, we have $|c_k|<2$.  For each $k\in\{0,...,r-1\}$ note that 
$\sqrt {\rho_{k,0}}=\frac{c_k+\sqrt{c_k^2+4}}{2}$, and $\sqrt{\rho_{k,1}} = \frac{c_k-\sqrt{c_k^2+4}}{2}$ are the two solutions of 
\[
z^2-c_kz-1=0,
\]
where $c_k$ is defined in \eqref{eq:ck}.  Using that $|c_k|<2$ we can now see that $|\sqrt\rho | \leq 1+\sqrt 2$ holds. As $\frac{1}{\rho}$ is also a root of \eqref{eq:charPolyFact} by Lemma \ref{lem:charPolyProps1}, we also have $|\sqrt \rho | \geq 1/(1+\sqrt 2)$. This establishes \eqref{eq:rhoBound2}. 
To obtain \eqref{eq:rhoBound1} note that since $p(\rho)=0$, it follows by \eqref{eq:charPolyFact} that $|\rho-1| =\lambda^{\frac{1}{2r}} |\sqrt \rho |$, and hence \eqref{eq:rhoBound2} implies the desired result.

The fact that $\rho_{0,0} \neq \rho_{0,1}$ and that $|\rho_{0,0}| =|\rho_{0,1}|=1$ can be readily obtained by direct calculation using \eqref{eq:root1} and \eqref{eq:root2} together with Lemma~\ref{lm:sigD} (to prove $|\rho_{0,0}|=|\rho_{0,1}|=1$, we note that $\lambda^{1/r} \in (0,4)$, and hence $\sqrt{\lambda^{2/r}-4\lambda^{1/r}}=\mathbbm{i}\sqrt{4\lambda^{1/r}-\lambda^{2/r}}$).  To finish we may now use the calculations in the paragraph above to see that
\begin{align}
 |  \rho_{k,0}| + | \rho_{k,1}| \notag  & = \frac{c_k+\sqrt{c_k^2+4}}{2}\cdot \frac{\overline{c_k}+\sqrt{\overline{c_k}^2+4}}{2}  + \frac{c_k - \sqrt{c_k^2+4}}{2}\cdot \frac{\overline{c_k}-\sqrt{\overline{c_k}^2+4}}{2} \notag \\ 
& = \frac{|c_k|^2 + \sqrt {|c_k|^4+ 4(c_k^2+\overline{c_k}^2)+16}}{2}. \label{eq:rhok} 
\end{align}
Then, we first observe that $|c_k|^4 = (|c_k^2|)^2 \ge (\textrm{Re}(c_k^2))^2$ (since $|c_k^2| \ge |\textrm{Re}(c_k^2)|$), and therefore we have
\begin{align}
|c_k|^4 +4 (c_k^2+ \overline{ c_k}^2 )+16 & \geq (\textrm{Re}(c_k^2))^2 +8 \textrm{Re}(c_k^2)+16 =  (\textrm{Re}(c_k^2)+4)^2 = (\textrm{Re}(c_k^2+4))^2 \notag \\ 
& =  (4 - \lambda^{1/r} \cos (2k \pi/ r))^2  \label{eq:est}
.\end{align}
Combining \eqref{eq:rhok} and \eqref{eq:est} and Lemma~\ref{lem:charPolyProps1} we can now see that 
\begin{align*}
|\rho_{k,0} |+  \left|\frac{1}{ \rho_{k,0}} \right| &= \frac{|c_k|^2 + \sqrt {|c_k|^4+ 4(c_k^2+\overline{c_k}^2)+16}}{2} = \frac{ \lambda^{1/r} + \sqrt {|c_k|^4+ 4(c_k^2+\overline{c_k}^2)+16}}{2} \\
 & \geq     \frac{ \lambda^{1/r} + 4 -\lambda^{1/r} \cos (2k \pi/ r)}{2}= 2+ \frac{\lambda^{1/r}  (1-  \cos (2k \pi/ r))}{2} > 2, \text{ for all } k \in \{1,2, \cdots, r-1\}.
\end{align*}
Thus, $|\rho_{k,0}| \neq 1, \text{ for all } k \in \{1,2, \cdots, r-1\} $.  The desired result follows.
\end{proof}

Lemma~\ref{lem:BasicRootsEqns} above tells us that all of the roots of the characteristic polynomial $p$ in \eqref{eq:charPoly} are contained in a disk of radius $(1+\sqrt{2}) \lambda^{\frac{1}{2r}}$ centered at $1$.  This information alone is enough for us to easily upper bound the distance between any two roots of $p$ by $2(1+\sqrt{2}) \lambda^{\frac{1}{2r}}$.  Obtaining {\it lower bounds} between the distances of the roots of $p$ from one another is a much more difficult task, however.  We will now begin the process of computing such lower bounds with the following lemma.  It establishes that all of the roots of the characteristic polynomial $p$ are unique so that their pairwise distances are nonzero.

\begin{lm}
The characteristic polynomial \eqref{eq:charPolyFact} always has $2r$ unique roots (with multiplicity one).
\label{lem:Uniquerootsforp}
\end{lm}

\begin{proof}
Suppose that $\rho$ is root of $p$ with multiplicity $> 1$.  We will consider two cases based on \eqref{eq:charPolyFact}.  First, suppose that 
$$-(1-\rho)^2 - \lambda^{1/r} \mathbbm{e}^{2 k \pi \mathbbm{i}/r} \rho = -(1-\rho)^2 - \lambda^{1/r} \mathbbm{e}^{2 l \pi \mathbbm{i}/r} \rho = 0$$
for $k \neq l$.  This can only occur if $\rho = 0$ since $\lambda > 0$ by Lemma~\ref{lm:sigD}, which then implies that $-(1-0)^2 = 0$ (a contradiction).  

Thus, it must instead be the case that 
$$-x^2  + (2 - \lambda^{1/r} \mathbbm{e}^{2 k \pi \mathbbm{i}/r})x - 1 = -(1-x)^2 - \lambda^{1/r} \mathbbm{e}^{2 k \pi \mathbbm{i}/r} x = c(x-\rho)^2 = c x^2 -2c \rho x + c \rho^2$$
for some $c \in \mathbbm{C}$ and $k \in \{ 0, \dots, r-1 \}$.  This in turn implies that $c = -1$ and $\rho^2 = 1$ must be true.  However, this also can't be the case since then we'd have
$$2 - \lambda^{1/r} \mathbbm{e}^{2 k \pi \mathbbm{i}/r} = \pm 2 \implies \textrm{either }\lambda^{1/r} \mathbbm{e}^{2 k \pi \mathbbm{i}/r} = 0 ~{\rm or }~ \lambda^{1/r} \mathbbm{e}^{2 k \pi \mathbbm{i}/r} = 4,$$
both of which are impossible since $\lambda \in (0,4^r)$ by Lemma~\ref{lm:sigD}.
\end{proof}

As a consequence of Lemma~\ref{lem:Uniquerootsforp} together with the discussion above, we can see that all $2r$ roots provided by \eqref{eq:root1} and \eqref{eq:root2} above are unique (i.e., with multiplicity one).  Therefore, the general solution to the recurrence relation \eqref{eq:RecRel} is 
\begin{equation}
\tilde{v}_{i} = \sum^{1}_{j=0} \sum^{r-1}_{k = 0} c_{k,j} \cdot \rho^i_{k,j}
\label{equ:RecSoln}
\end{equation}
for all $i \in \mathbbm{Z}$, where the $c_{k,j} \in \mathbbm{C}$ are chosen so that that the first line of \eqref{eq:RRecRel1} together with \eqref{eq:RRecRel1BC} both hold.

\subsection{Additional Properties of the Roots of the Characteristic Polynomial \eqref{eq:charPoly}}
\label{sec:RootsResCharpoly}

Unfortunately, the uniqueness of the roots of $p$ alone will ultimately not be enough for our purposes below.  We will also require lower bounds on their distances from one another.  The following lemmas provide such estimates.

\begin{lm}\label{lm:lambda}
For any two roots of $p(x)$, $\rho \neq \tilde{\rho}$, either $\rho = \tilde{\rho}$, $\bar\rho = \tilde\rho$, $\rho^{-1} = \tilde\rho$, $\overline{\rho}^{-1} = \tilde \rho$, or 
\begin{equation}\label{eq:modulus_difference_bound}
c r^{-2}\lambda^{1/2r}\leq \left||\tilde{\rho}|-|{\rho}|\right|\leq C \lambda^{1/2r},
\end{equation}
where $c, C \in \mathbbm{R}^+$ are both absolute constants (i.e., universal constants independent of $N, r, \lambda,$ etc.).
\end{lm}
%%%
\begin{proof}
See Appendix~\ref{sec:ProofofRootBounds1}.
\end{proof}

\begin{lm} \label{lm:lambda2}
Let $\rho, \rho' \in \mathbbm{C}$, $\rho \neq \rho'$, be two roots of \eqref{eq:charPoly} with $|\rho| \neq 1$. Then, there exist absolute constants $C,c,c_1, c_2 \in \mathbbm{R}^+$, $c_1 > 1$, such that
\begin{equation}
t_r(\rho,\rho')\lambda^{\frac{1}{2r}} \leq |\rho-\rho'| \leq C\lambda^{\frac{1}{2r}},
\label{eq:upper_and_lower_bound_diff}
\end{equation}
where $t_r(\rho,\rho') \geq c_2 c_1^{-r}$ if either $\bar\rho= \rho'$ or $\rho' = \rho^{-1}$ holds, and $t_r(\rho,\rho') \geq c r^{-2}$ otherwise.
\end{lm}

\begin{proof}  
See Appendix~\ref{sec:ProofofRootBounds2}.
\end{proof}

Lemmas~\ref{lm:lambda} and~\ref{lm:lambda2} collectively bound the distances between all roots of $p$ from below \textit{except for $|\rho_{0,0} - \rho_{0,1}|$}, the distance between the two unimodular roots of $p$.  Thankfully, however, simply knowing that this single distance is nonzero will suffice below.
Finally, we conclude this section with a corollary of Lemma~\ref{lem:charPolyProps1}.  It characterizes when the roots of the characteristic polynomial $p$ will be real.

\begin{corollary}\label{cor:onlyRealRoots} 
The roots $\rho_{k,j} \in \mathbb{R}$ if and only if $r$ is even and $k=r/2$.
\end{corollary}
\begin{proof}
First, we see by Lemma~\ref{lem:charPolyProps1} and (\ref{eq:charPolyFact}) that $\rho_{k,0},\rho_{k,1}$ are the roots of 
\[q_k(x):=-(1-x)^2 - \lambda^{1/r} \mathbbm{e}^{2 k \pi \mathbbm{i}/r} x.\]
It is clear that $q_k(0)=-1$, so $q_k$ does not have a root at $0$. We also note that if $x \in \mathbb{R} \setminus \{0\}$, then $q_k(x) \not \in \mathbb{R}$ if $\mathbbm{e}^{2 k \pi \mathbbm{i}/r} \not \in \mathbb{R}$, so no such values of $k$ will lead to real roots of $q_k$. Therefore, it suffices to only consider values of $k$ for which $\mathbbm{e}^{2 k \pi \mathbbm{i}/r} \in \mathbb{R}$, namely $k=0$ and $k=r/2$ for $r$ even. 

By (\ref{eq:root1}) and (\ref{eq:root2}), we can see that the roots of $q_k$ are
\[ \frac{2 - \lambda^{1/r} \mathbbm{e}^{2 k \pi \mathbbm{i}/r} \pm \sqrt{\lambda^{2/r} \mathbbm{e}^{4 k \pi \mathbbm{i}/r}-4\lambda^{1/r} \mathbbm{e}^{2 k \pi \mathbbm{i}/r}}}{2}\]
and so it suffices to check the sign of $\lambda^{2/r} \mathbbm{e}^{4 k \pi \mathbbm{i}/r}-4\lambda^{1/r} \mathbbm{e}^{2 k \pi \mathbbm{i}/r}$ to determine whether or not these roots are real in this case.  We note that $x^2-4x = x(x-4) < 0$ if and only if $x \in (0,4)$ and $x^2-4x > 0$ if and only if $x \notin (0,4)$. Let $x:= \lambda^{1/r} \mathbbm{e}^{2 k \pi \mathbbm{i}/r}$ and note that $\lambda^{1/r} \in (0,4)$ by Lemma \ref{lm:sigD}. When $k=0$ we have $x = \lambda^{1/r}$ and $x(x-4) = \lambda^{2/r}-4\lambda^{1/r}$ which is negative since $x = \lambda^{1/r} \in (0,4)$ when $k = 0$. On the other hand, when $k=r/2$ we have $x = -\lambda^{1/r}$ and $x(x-4) =  \lambda^{2/r}+4\lambda^{1/r} > 0$ since $x = -\lambda^{1/r} \notin (0,4)$ in this case. Thus, $\rho_{k,j} \in \mathbb{R}$ if and only if $r$ is even and $k=r/2$.
\end{proof}

We are now prepared to begin proving Lemma~\ref{lm:large}.

%%%%%%%%%%%%%%%%%%%%%%%%%%%%%%%%%%%%%%%%%%%%%%%%%%%%%%%%%%%%%%
%%%%%%%%%%%%%%%%%%%%%%%%%%%%%%%%%%%%%%%%%%%%%%%%%%%%%%%%%%%%%%
\section{The Proof of Lemma~\ref{lm:large}}
\label{sec:proof_BIG_lemma}

\input{ProofofBIGlemma.tex}

\section*{Acknowledgements}
The authors would like to thank Wei-Hsuan Yu for reading and commenting on an initial draft of the proof of Lemma~\ref{lm:large} while he was a postdoc at Michigan State University.

\bibliographystyle{abbrv}
\bibliography{ms}

\begin{appendix}

%%%%%%%%%%%%%%%%%%%%%%%%%%%%%%%%%%%%%%%%%%%%%%%%%%%%%%%%%%%%%%
%%%%%%%%%%%%%%%%%%%%%%%%%%%%%%%%%%%%%%%%%%%%%%%%%%%%%%%%%%%%%%
\section{Proving the Basic Results:  Lemmas~\ref{lem:Rightinf=Leftinf} and~\ref{lem:EntriesofDrTDr}}

We will begin with the proof of Lemma \ref{lem:Rightinf=Leftinf}.

\subsection{Proof of Lemma~\ref{lem:Rightinf=Leftinf}}
\label{sec:Proofoflem:Rightinf=Leftinf}

To begin the proof that $\|\mathbf{u}_j\|_{\infty}= \|\mathbf{v}_j\|_{\infty}$ for all $j \in [N]$, we first observe that by Definition \ref{Def:diffMat}, $D_{i,j} = D_{N-j,N-i}$. We now prove by induction that 
\begin{equation}\label{eq:Dproperty1}
(D^r)_{i,j} = (D^r)_{N-j,N-i}
\end{equation}
 for any $r \in \mathbb{N}$. Suppose that $(D^r)_{i,j} = (D^r)_{N-j,N-i}$. Then we have 
\begin{align}\label{eq:Dproperty1res1}
(D^{r+1})_{i,j} &= \sum_{k=1}^N (D^r)_{i,k} (D)_{k,j} \nonumber \\
&= \sum_{k=1}^N (D^r)_{N-k,N-i} (D)_{N-j,N-k},
\end{align}
where in the last equality we used the inductive hypothesis $(D^r)_{i,j} = (D^r)_{N-j,N-i}$ and the fact that $D_{i,j} = D_{N-j,N-i}$. Then making the change of variables $k' = N-k$, we have
\begin{align}\label{eq:Dproperty1res2}
\sum_{k=1}^N (D^r)_{N-k,N-i} (D)_{N-j,N-k} &=  \sum_{k'=1}^N (D^r)_{k',N-i} (D)_{N-j,k'} \nonumber\\
&=   \sum_{k'=1}^N (D)_{N-j,k'}(D^r)_{k',N-i} \nonumber\\
&= (D^{r+1})_{N-j,N-i},
\end{align}
and hence combining (\ref{eq:Dproperty1res1}) and (\ref{eq:Dproperty1res2}) we have $(D^{r+1})_{i,j} =  (D^{r+1})_{N-j,N-i}$, completing the proof by induction.

Next, we claim that 
\begin{equation}\label{eq:Dproperty2}
(D^r(D^r)^T)_{i,j} = ((D^r)^TD^r)_{N-i,N-j}.
\end{equation}
We have
\begin{align}\label{eq:Dproperty2res1}
(D^r(D^r)^T)_{i,j} &= \sum_{k=1}^N (D^r)_{i,k} (D^r)_{j,k} \nonumber\\
&= \sum_{k=1}^N (D^r)_{N-k,N-i} (D^r)_{N-k,N-j}
\end{align}
where the last equality holds by (\ref{eq:Dproperty1}).
Now, making the change of variables $k'=N-k$, we have
\begin{align}\label{eq:Dproperty2res2}
\sum_{k=1}^N (D^r)_{N-k,N-i} (D^r)_{N-k,N-j} &=\sum_{k'=1}^N (D^r)_{k',N-i} (D^r)_{k',N-j}\nonumber\\
&= ((D^r)^TD^r)_{N-i,N-j}, \nonumber\\
\end{align}
so combining \eqref{eq:Dproperty2res1} and \eqref{eq:Dproperty2res2}, we see that $(D^r(D^r)^T)_{i,j} = ((D^r)^TD^r)_{N-i,N-j}$, verifying \eqref{eq:Dproperty2}.

Finally, for ease of notation, we let $\tilde{\mathbf{v}}$ be the vector $\mathbf{v}$ written in reverse order, i.e.
\[(\tilde{\mathbf{v}})_i := \mathbf{v}_{N-i}.\] 
Observe that \eqref{eq:Dproperty2} implies that $D^r(D^r)^T$ and $(D^r)^TD^r$ have reversed row and column orders. Hence, if $\mathbf{v}$ is an eigenvector of $D^r(D^r)^T$ with eigenvalue $\lambda$, then $\tilde{\mathbf{v}}$ is an eigenvector of $(D^r)^TD^r$ with eigenvalue $\lambda$. In other words, the eigenvectors of $(D^r)^TD^r$ can be obtained by reversing the order of each eigenvector of $D^r(D^r)^T$, and vice versa. Since the eigenvectors of $D^r(D^r)^T$ and $(D^r)^TD^r$ correspond to the left and right singular vectors of $D^r$, respectively, the same relationship holds between the left and right singular vectors of $D^r$. In particular, the left singular vectors can be obtained by reversing the order of each right singular vector, and vice versa. Thus, since reversing the order of a vector does not change its $\ell^\infty$-norm, we have
\[\|\mathbf{v}_j\|_{\infty}= \|\mathbf{u}_j\|_{\infty}\]
for all $j \in [N]$, the desired result.

\subsection{Proof of Lemma~\ref{lem:EntriesofDrTDr}}
\label{sec:ProofofEntriesofDrTDr}

Recall that we are seeking to prove the following lemma: 
\newtheorem*{lem:EntriesofDrTDr}{\textbf{\emph{Lemma \ref{lem:EntriesofDrTDr}}}}
\begin{lem:EntriesofDrTDr}
Let $r, N \in \mathbbm{Z}^+$ be such that $r < N / 2$.  All the entries of $(D^r)^TD^r \in \mathbbm{Z}^{N \times N}$ are given by $$((D^r)^TD^r)_{j-m,j} = \left\{
\begin{array}{ll}
        \displaystyle    (-1)^m  {2r \choose r-m} & \text{ if } j-m \leq N -r , ~m \in \{ 0, 1, \cdots, r \} \\
& \\
      
      \displaystyle (-1)^m \sum_{l=0}^{N-j} {r \choose l+m}{r \choose l} &  \text{ if } j-m > N-r,~m \in \{ 0, 1, \cdots, r\} \\
      \hspace{.5in}0 & \text{ if } m > r
\end{array} 
\right.,$$
combined with the fact that $(D^r)^TD^r$ is symmetric.
\end{lem:EntriesofDrTDr}
Define the upper and lower triangular nilpotent shift matrices $U \in \mathbb{R}^{N \times N}$ and $L \in \mathbb{R}^{N \times N}$ as
\begin{equation*}
U_{i,j} := \left\{ \begin{array}{ll} 1 & \textrm{if}~i = j-1\\ 0 & \textrm{otherwise}\end{array} \right.
\label{Def:upShift}
\end{equation*}
and
\begin{equation*}
L_{i,j} := \left\{ \begin{array}{ll} 1 & \textrm{if}~i = j+1\\ 0 & \textrm{otherwise}\end{array} \right..
\label{Def:downShift}
\end{equation*}
Note that $D = I - L$ so that 
\begin{equation*}
D^r = \sum^{r}_{k=0} {r \choose k} (-1)^k L^k.
\end{equation*}
Similarly, $D^{\rm T} = I - U$ so that 
\begin{equation*}
\left(D^r \right)^{\rm T} = \sum^{r}_{k=0} {r \choose k} (-1)^k U^k.
\end{equation*}
Since we are interested in the right singular vectors of $D^r$ we will consider the symmetric matrix
\begin{equation}
\left(D^r \right)^{\rm T} D^r ~=~ \left( \sum^{r}_{k=0} {r \choose k} (-1)^k U^k \right) \left( \sum^{r}_{k=0} {r \choose k} (-1)^k L^k \right) ~=~ \sum^{r}_{k,l=0} {r \choose k}{r \choose l} (-1)^{k+l} U^kL^l.
\label{eq:DrSym}
\end{equation}

Note that $U^kL^l = X_{k,l}$, where
\begin{equation}
\left(X_{k,l} \right)_{i,j} = \left\{ \begin{array}{ll} 1 & \textrm{if}~j \leq N - l~{\rm and}~i = j-k+l > 0\\ 0 & \textrm{otherwise}\end{array} \right..
\label{eq:DefX}
\end{equation}
Thus, if $j \leq N-r$ and $i = j - m$, $m \in \{0,1,\dots,r \}$, we will have
\begin{align*} 
\left(\left(D^r \right)^{\rm T} D^r \right)_{i,j} &=~ \left( \sum^{r}_{k,l=0} {r \choose k}{r \choose l} (-1)^{k+l} X_{k,l} \right)_{i,j} ~=~ (-1)^{m}\sum^{r-m}_{l=0} {r \choose l + m}{r \choose l}\\ &=~ (-1)^{m}\sum^{r-m}_{l=0} {r \choose (r-m)-l}{r \choose l}.
\end{align*}
Simplifying the expression above using Vandermonde's identity we can now see that
\begin{equation*}
\left(\left(D^r \right)^{\rm T} D^r  \right)_{j-m,j} ~=~ (-1)^{m}\sum^{r-m}_{l=0} {r \choose (r-m)-l}{r \choose l} ~=~ (-1)^{m}{2r \choose r-m}
\end{equation*}
for all $j  \leq N-r$, $m \in \{0,1,\dots,r \}$.  By inspecting \eqref{eq:DrSym} and \eqref{eq:DefX} it is not difficult to see that, more generally, we will have
\begin{equation}
\left(\left(D^r \right)^{\rm T} D^r \right)_{j-m,j} = \left\{ \begin{array}{ll} (-1)^{m}{2r \choose r-m} & \textrm{if}~ 0 \leq m \leq r\\ 0 & \textrm{if}~m > r\end{array} \right.
\label{Def:DrSym}
\end{equation}
for all $j \leq N-r$.
In fact, \eqref{Def:DrSym} gives $\left( \left(D^r \right)^{\rm T} D^r \right)_{i,j}$ for all $i,j \in [N]$ with $\max\{i,j\} \leq N-r$ by symmetry.  
If $j > N-r$ and $ i = j-m$ for $m \in \{0,1, \cdots, r\}$, then
\begin{align*}
\left(\left(D^r \right)^{\rm T} D^r \right)_{i,j} &~=~ \left( \sum^{r}_{k,l=0} {r \choose k}{r \choose l} (-1)^{k+l} X_{k,l} \right)_{i,j} ~=~   \sum^{r}_{k,l=0} {r \choose k}{r \choose l} (-1)^{k+l} (X_{k,l})_{i,j} \\
&=(-1)^{m}\sum^{\min \{ N-j, r-m \}}_{l=0} {r \choose l + m}{r \choose l} ~=~ (-1)^{m}{2r \choose r-m}
\end{align*}
whenever $r-m \leq N-j$, or equivalently, whenever $i \leq N -r$.  Otherwise, when $j > N -r + m$, or equivalently, when $i > N - r$, this last equation becomes
\begin{equation*}
\left(\left(D^r \right)^{\rm T} D^r \right)_{j-m,j} = (-1)^{m}\sum^{N-j}_{l=0} {r \choose l + m}{r \choose l}.
\end{equation*}
Also note that, by the same argument that was used in \eqref{Def:DrSym}, we will have $\left(\left(D^r \right)^{\rm T} D^r \right)_{j-m,j} = 0$ for $m>r$ in the case that $j>N-r$.
Utilizing symmetry now allows us to determine all the entries of $\left(D^r \right)^{\rm T} D^r$, and completes the proof of Lemma \ref{lem:EntriesofDrTDr}.

%%%%%%%%%%%%%%%%%%%%%%%%%%%%%%%%%%%%%%%%%%%%%%%%%%%%%%%%%%%%%%
%%%%%%%%%%%%%%%%%%%%%%%%%%%%%%%%%%%%%%%%%%%%%%%%%%%%%%%%%%%%%%
\section{Proof of Lemma~\ref{lem:proofofeq:RecRel1}}
\label{sec:Proofoflem:proofofeq:RecRel1}

We need to show that under the condition $(D^r)^TD^r\mathbf{v} = \lambda \mathbf{v}$, the infinite sequence $\tilde{\mathbf{v}}$ defined by  $\eqref{eq:RRecRel1}-\eqref{eq:RRecRel2}$ satisfies  \eqref{eq:RecRel} for each $i \in \mathbb{Z}$. We will divide $i \in \mathbb{Z}$ into five regimes: $ i\leq 0, 0< i \leq r , r < i \leq N-r, N-r+1\leq i \leq N, i > N $. In the first and last regimes, \eqref{eq:RecRel} trivially holds because of the way the sequence $\tilde{\bf v}$ is extended in \eqref{eq:RRecRel2}.  For the second and third regimes, it is easy to verify using Lemma~\ref{lem:EntriesofDrTDr} and \eqref{eq:RRecRel1} that  \eqref{eq:RecRel}  is exactly the $i^{\rm th}$ row of $(D^r)^TD^r \mathbf{v} = \lambda \mathbf{v}$, and hence holds true. It only remains to prove  \eqref{eq:RecRel} for $i \in \{N-r+1,\dots,N\}$.  In this case, we will show that  \eqref{eq:RecRel} is implied by the $i^{\rm th}$ equation in the system $(D^r)^TD^r \mathbf{v} = \lambda \mathbf{v}$ and  \eqref{eq:RRecRel1BC}. 

Let $ x= i-N+r$.  Lemma~\ref{lem:EntriesofDrTDr} then tells us that the $i^{\rm th}$ equation in $(D^r)^TD^r \mathbf{v} = \lambda \mathbf{v}$ is
\begin{align} \label{eq:ith_equ}
\lambda \widetilde{v}_i = \sum\limits_{k=0}^{r-x} (-1)^{r-k} \binom{2r}{k} \tilde{v}_{k+i-r} &+ \sum\limits_{k=r-x+1}^r (-1)^{r-k} \sum\limits_{l=0}^{r-x} \binom{r}{l+r-k}\binom{r}{l} \tilde{v}_{k+i-r} \\
&+ \sum\limits_{k=r+1}^{2r-x} (-1)^{k-r} \sum\limits_{l=0}^{2r-x-k} \binom{r}{l+k-r}\binom{r}{l} \tilde{v}_{k+i-r} 
 \notag 
\end{align}
for all $i \in \{N-r,...,N\}$ (so that $x = i-N+r \in \{ 0, \dots, r \} $).

Note that the righthand side of equation \eqref{eq:ith_equ} has three terms in accordance with Lemma~\ref{lem:EntriesofDrTDr}.  The first term involves entries $\tilde{v}_{k+i-r}$ with $k + i - r \leq N-r$ (i.e., before the entires associated with the irregular lower-right $r \times r$ submatrix of $\left(D^r \right)^{\rm T} D^r$), the second term involves entries $\tilde{v}_{k+i-r}$ with $i \geq k + i - r > N-r$ (left up to the diagonal), and the third involves entries $\tilde{v}_{k+i-r}$ with $k + i - r > i$ (right of the diagonal). Furthermore, the lefthand side of \eqref{eq:ith_equ} matches the righthand side of \eqref{eq:RecRel}.  Hence, if we can show that the righthand side of \eqref{eq:ith_equ} matches the lefthand side of \eqref{eq:RecRel} we will be finished with our proof.  We will accomplish this task below by showing that the difference between the lefthand side of \eqref{eq:RecRel} and the righthand side of \eqref{eq:ith_equ} is always zero.

Let the function $f: \{0,1,\cdots,r\} \times \{N-r+1, \cdots, N\} \rightarrow \mathbbm{C}$ be defined to be the related difference
\begin{align}\label{eq:deff}
f(x,i) :=& \sum\limits_{k=r-x+1}^r (-1)^{r-k} \sum\limits_{l=0}^{r-x} \binom{r}{l+r-k}\binom{r}{l} \tilde{v}_{k+i-r}\\ &+ \sum\limits_{k=r+1}^{2r-x} (-1)^{k-r} \sum\limits_{l=0}^{2r-x-k} \binom{r}{l+k-r}\binom{r}{l} \tilde{v}_{k+i-r} - \sum\limits_{k=r-x+1}^{2r} (-1)^{r-k} \binom{2r}{k} \tilde{v}_{k+i-r}  \notag 
\end{align}
where the first vacuous term is ignored when $x=0$. 
As per the preceding discussion, the lemma will be proven if we can show that $f(i-N+r,i) = 0$ for all $i \in \{N-r+1,\cdots, N\}$.  To show this, we will now prove that both
\begin{enumerate}[(a)]
\item \label{en:a}$f(0,i)=0$, and
\item\label{en:b} $f(x,i)-f(x-1,i) =0 $ for all $x\in \{1,\dots,i-N+r\}$
\end{enumerate}
hold for all $i \in \{N-r+1, N\}$.  As long as \eqref{en:a} and \eqref{en:b} above both hold, we can then deduce for any given $i \in \{N-r+1,\dots,N\}$ that 
$$f(i-N+r,i) = f(0,i) + \sum\limits_{x=1}^{i-N+r} f(x,i)-f(x-1,i)=0$$
as desired.  Thus, it suffices to prove both \eqref{en:a} and \eqref{en:b} in order to finish our proof of Lemma~\ref{lem:proofofeq:RecRel1}.

Both \eqref{en:a} and \eqref{en:b} can be verified by direct calculation. For \eqref{en:a}, we can see from \eqref{eq:deff} that
\begin{align*}
 f(0,i)& = \sum\limits_{k=r+1}^{2r} (-1)^{k-r} \sum\limits_{l=0}^{2r-k} \binom{r}{l+k-r}\binom{r}{l}\tilde{v}_{k+i-r} - \sum\limits_{k=r+1}^{2r} (-1)^{r-k}\binom{2r}{k} \tilde{v}_{k+i-r} \\
 & = \sum\limits_{k=r+1}^{2r} (-1)^{k-r} \left(\sum\limits_{l=0}^{2r-k} \binom{r}{l+k-r}\binom{r}{l}- \binom{2r}{k}\right)\tilde{v}_{k+i-r} = 0.
\end{align*}
Here the last equality is obtained by observing that all the coefficients of $\tilde{v}_{k+i-r}$ are 0 via Vandermonde's identity.
Thus, it suffices to prove \eqref{en:b} in order to finish our proof of Lemma~\ref{lem:proofofeq:RecRel1}.  

To verify \eqref{en:b} we can now use \eqref{eq:deff} to see that
\begin{align*}
f(x,i)\hspace{-.02in}-\hspace{-.04in}f(x \hspace{-.04in} - \hspace{-.04in}1,i) = &(-1)^{x-1} \sum\limits_{l=0}^{r-x} \binom{r}{l+x-1}\binom{r}{l}\tilde{v}_{i-x+1} - \hspace{-.15in} \sum\limits_{k=r-x+2}^r \hspace{-.05in} (-1)^{r-k} \binom{r}{2r-x-k+1}\binom{r}{x-1} \tilde{v}_{k+i-r} \\
&- (-1)^{r-x+1}\binom{r}{x-1}\tilde{v}_{r-x+1+i} -\sum\limits_{k=r+1}^{2r-x} (-1)^{k-r} \binom{r}{r-x+1}\binom{r}{2r-x-k+1} \tilde{v}_{k+i-r}\\
&-(-1)^{x-1} \binom{2r}{r-x+1} \tilde{v}_{i-x+1}
\end{align*}
where the first, second, and third lines of the righthand side above result form the differences between the first, second, and third terms in \eqref{eq:deff} for $f(x,i)$ and $f(x-1,i)$, respectively.  

Simplifying the equation directly above we get that
\begin{align*}
f(x,i)\hspace{-.02in}-\hspace{-.04in}f(x \hspace{-.04in} - \hspace{-.04in}1,i) = &(-1)^{x-1} \sum\limits_{l=0}^{r-x} \binom{r}{l+x-1}\binom{r}{l}\tilde{v}_{i-x+1} - \hspace{-.15in} \sum\limits_{k=r-x+2}^{2r-x} \hspace{-.05in} (-1)^{r-k} \binom{r}{2r-x-k+1}\binom{r}{x-1} \tilde{v}_{k+i-r} \\
&~~~~~- (-1)^{r-x+1}\binom{r}{x-1}\tilde{v}_{r-x+1+i} -(-1)^{x-1} \binom{2r}{r-x+1} \tilde{v}_{i-x+1} \\
= &(-1)^{x-1} \sum\limits_{l=0}^{r-x} \binom{r}{l+x-1}\binom{r}{l}\tilde{v}_{i-x+1} - \hspace{-.15in} \sum\limits_{k=r-x+2}^{2r-x+1} \hspace{-.05in} (-1)^{r-k} \binom{r}{2r-x-k+1}\binom{r}{x-1} \tilde{v}_{k+i-r} \\
&~~~~~-(-1)^{x-1} \binom{2r}{r-x+1} \tilde{v}_{i-x+1} \\
= &(-1)^{x-1} \sum\limits_{l=0}^{r-x} \binom{r}{l+x-1}\binom{r}{l}\tilde{v}_{i-x+1} - \hspace{-.15in} \sum\limits_{k=r-x+2}^{2r-x+1} \hspace{-.05in} (-1)^{r-k} \binom{r}{2r-x-k+1}\binom{r}{x-1} \tilde{v}_{k+i-r} \\
&~~~~~-(-1)^{x-1} \binom{2r}{r-x+1} \tilde{v}_{i-x+1} +(-1)^{x-1}\binom{r}{x-1}\tilde{v}_{i-x+1} - (-1)^{x-1}\binom{r}{x-1}\tilde{v}_{i-x+1}\\
= &(-1)^{x-1} \hspace{-.1in} \sum\limits_{l=0}^{r-x+1} \binom{r}{l+x-1}\binom{r}{l}\tilde{v}_{i-x+1} - \hspace{-.15in} \sum\limits_{k=r-x+1}^{2r-x+1} \hspace{-.05in} (-1)^{r-k} \binom{r}{2r-x-k+1}\binom{r}{x-1} \tilde{v}_{k+i-r} \\
&~~~~~-(-1)^{x-1} \binom{2r}{r-x+1} \tilde{v}_{i-x+1}.
\end{align*}
Using Vandermonde's identity we can now see that the first and third terms cancel.  Thus we have that
\begin{align*}
f(x,i)-f(x-1,i) &= - \sum\limits_{k=r-x+1}^{2r-x+1} (-1)^{r-k} \binom{r}{2r-x-k+1}\binom{r}{x-1} \tilde{v}_{k+i-r}\\
&= \binom{r}{x-1} (-1)^{x+r} \sum \limits_{j=0}^{r} (-1)^{j} \binom{r}{j} \tilde{v}_{i-x+r-j+1}.
\end{align*}
Note that $i-x+r+1 \in \{ N+1, \dots, N+r \}$ whenever $i \in \{N-r+1,\dots,N\}$ and $x\in \{1,\dots,i-N+r\}$.  Thus, \eqref{en:b} will hold by \eqref{eq:RRecRel1BC}.  This finishes the proof.

%%%%%%%%%%%%%%%%%%%%%%%%%%%%%%%%%%%%%%%%%%%%%%%%%%%%%%%%%%%%%%
%%%%%%%%%%%%%%%%%%%%%%%%%%%%%%%%%%%%%%%%%%%%%%%%%%%%%%%%%%%%%%
\section{Proof of Lemma~\ref{lm:lambda}}
\label{sec:ProofofRootBounds1}

We begin with some facts from the proof of Lemma~\ref{lem:BasicRootsEqns}.
For each $k\in\{0,...,r-1\}$ recall that $\sqrt {\rho_{k,0}}=\frac{c_k+\sqrt{c_k^2+4}}{2}$, and $\sqrt{\rho_{k,1}} = \frac{c_k-\sqrt{c_k^2+4}}{2}$ are the two solutions to 
\[
z^2-c_kz-1=0,
\]
where $c_k$ is defined in \eqref{eq:ck}.
Hence, by Vieta's formulas, we observe that $\sqrt{\rho_{k,0}}\sqrt{\rho_{k,1}}=-1$, i.e. $\sqrt{\rho_{k,0}}=-\frac{1}{\sqrt{ \rho_{k,1}}}$. Let $\tilde{\rho}_k$ be either one of the roots $\rho_{k,0}$ or $\rho_{k,1}$. 
To bound $|\tilde\rho_{k}|-|\tilde{\rho}_l|$ we will make use of the equality 
\begin{equation}  \left| | \tilde{\rho}_k| + |1/\tilde{\rho}_k| - (| \tilde{\rho}_l| + |1/ \tilde{\rho}_l| ) \right| = \left|(|\tilde\rho_{k}|-|\tilde{\rho}_l|)\frac{|\tilde{\rho}_k|| \tilde{\rho}_l|-1}{|\tilde{\rho}_k | |\tilde{\rho}_l|}\right|.\label{eq:helperEq}
\end{equation}
In particular, our strategy is to bound both the left hand side, and the term on the right hand side of \eqref{eq:helperEq} that multiplies $\left( |\tilde\rho_{k}|-|\tilde{\rho}_l| \right)$, from both above and below. 

We begin with the left hand side, and recall using \eqref{eq:rhok} that
\begin{align*}
 | \tilde\rho_{k}| + |1/\tilde\rho_{k}| &=|  \rho_{k,0}| + | \rho_{k,1}| \notag  = \frac{|c_k|^2 + \sqrt {|c_k|^4+ 4(c_k^2+\overline{c_k}^2)+16}}{2}.
\end{align*}
Consequently,
\begin{align}\label{eq:VenusFlyTrap}
| \tilde{\rho}_k| + |1/\tilde{\rho}_k| - (| \tilde{\rho}_l| + |1/ \tilde{\rho}_l| )& = \frac{|c_k|^2 + \sqrt {|c_k|^4+ 4(c_k^2+\overline{c_k}^2)+16}}{2}- \frac{|c_{l} |^2 + \sqrt {|c_{l}|^4+ 4(c_{l}^2+\overline{c_{l}}^2)+16}}{2} \notag \\
&=\frac{4(\textrm{Re}(c_k^2-c_{l}^2))}{ \sqrt {|c_k|^4+ 4(c_k^2+\overline{c_k}^2)+16} + \sqrt {|c_{l}|^4+ 4(c_{l}^2+\overline{c_{l}}^2)+16}}, 
\end{align}
where the last equality holds since $|c_k|=|c_l|$ by \eqref{eq:ck}.

We now observe that a lower bound for the left hand side of \eqref{eq:helperEq} holds because the denominator in \eqref{eq:VenusFlyTrap} satisfies
\[
0 < \sqrt {|c_k|^4+ 4(c_k^2+\overline{c_k}^2)+16} + \sqrt {|c_{l}|^4+ 4(c_{l}^2+\overline{c_{l}}^2)+16}\leq 16
\]
while the numerator, provided that $c_k \neq \overline{c}_{l}$ and $c_k \neq c_{l}$ s.t. $k+l \neq r$ and $k \neq l$,  satisfies 
\[
| \textrm{Re}(c_k^2-c_l^2) |  = \lambda^{1/r}| \cos( 2\pi k/r )-\cos(2\pi l /r)| = \lambda^{1/r} |2\sin(\pi(k-l)/r )\sin(\pi(k+l)/r)| \geq 2\sin^2 (\pi/2r) \lambda^{1/r} .
\]
This implies that
\begin{equation}\label{eq:lambda_kl_lower}
 d_1(r) \lambda^{1/r} \leq  \left| | \tilde{\rho}_k| + |1/\tilde{\rho}_k| - (| \tilde{\rho}_l| + |1/ \tilde{\rho}_l| ) \right| ,
\end{equation}
where  $d_1(r): = c'/r^2$ for $c'>0$ such that $c'/r^2 \le \frac{\sin^2 (\pi/2r)}{8}$ for all $r \ge 2$.
Now that we have established \eqref{eq:lambda_kl_lower}  we are almost done with proving \eqref{eq:modulus_difference_bound}. 

To finish the argument establishing \eqref{eq:modulus_difference_bound} note that at most one of $\rho =: \tilde{\rho}_l$ and $\tilde{\rho} =: \tilde{\rho}_k$ can be unimodular by Lemma~\ref{lem:BasicRootsEqns}. Thus, we only need consider two cases: 
First, assume without loss of generality that $| \tilde{\rho}_l |=1$ and $|\tilde{\rho}_k| \neq 1$. Then $k \neq l$, $k + l \neq r$, and
\[
 \left| | \tilde{\rho}_k| + |1/\tilde{\rho}_k| - (| \tilde{\rho}_l| + |1/ \tilde{\rho}_l| ) \right| = \left|\frac{(|\tilde{\rho}_k|-1)^2}{|\tilde{\rho}_k|}\right|.
\] 
Plugging the above into \eqref{eq:lambda_kl_lower} and noting that \eqref{eq:rhoBound2} guarantees $|\tilde\rho_{k}| \geq (1+\sqrt{2})^{-2}$, we obtain 
\begin{equation}\label{eq:lambda1}
(1+\sqrt{2}) \lambda^{1/2r} \geq |\tilde\rho_{k} -1| \geq | |\tilde{\rho}_k|-1 | \geq (1+\sqrt{2})^{-1} \sqrt{  d_1(r)} \lambda^{1/2r},
\end{equation}
which uses \eqref{eq:rhoBound1} for the upper bound.  

Second, assume that both $|\tilde\rho_{k}| \neq 1$ and $|\tilde\rho_{l}|\neq 1$. Then
$$  \left| | \tilde{\rho}_k| + |1/\tilde{\rho}_k| - (| \tilde{\rho}_l| + |1/ \tilde{\rho}_l| ) \right| = \left| (|\tilde{\rho}_k|-|\tilde{\rho}_l|)\frac{(|\tilde{\rho}_k|-1)| \tilde{\rho}_l|+(|\tilde{\rho}_l|-1)}{|\tilde{\rho}_k | |\tilde{\rho}_l|} \right|
.$$
If $\rho \neq \tilde{\rho}$, $\bar\rho \neq \tilde\rho$, $\rho^{-1} \ne \tilde\rho$, and $\overline{\rho}^{-1} \ne \tilde \rho$ all hold then we have $k \neq l$, $k + l \neq r$ by Lemmas~\ref{lem:charPolyProps1} and~\ref{lem:BasicRootsEqns}. Thus,
we may use \eqref{eq:lambda_kl_lower}, along with the bounds on $\tilde{\rho}_{k}$ and $\tilde\rho_{l}$ from  \eqref{eq:rhoBound2} and \eqref{eq:rhoBound1}, to see that
\begin{align}
||\tilde{\rho}_k|-|\tilde{\rho}_l|| &\geq \left|\frac{(|\tilde{\rho}_k|-1)| \tilde{\rho}_l|+(|\tilde{\rho}_l|-1)}{| \tilde{\rho}_k | |\tilde{\rho}_l|}\right|^{-1} d_1(r)\lambda^{1/r} \notag
\\& \geq \frac{| \tilde{\rho}_k | |\tilde{\rho}_l|} {(1+\sqrt{2}) \lambda^{1/2r}| \tilde{\rho}_l|+(1+\sqrt{2}) \lambda^{1/2r}}d_1(r)\lambda^{1/r} \notag
\\& \geq d_1(r) (1+\sqrt{2})^{-1} \frac{|\tilde{\rho}_k\tilde{\rho}_l| }{ |\tilde\rho_{l}|+1}  \lambda^{1/2r} \notag 
\\&\geq e(r) \lambda^{1/2r},
\end{align}
where $e(r):= c/r^2$ is such that $e(r) \le d_1(r)(1+\sqrt{2})^{-1} \frac{(1+\sqrt{2})^{-4}}{1+(1+\sqrt{2})^2}$. 
On the other hand, from \eqref{eq:lambda1}, we have
$$
| | \tilde{\rho}_k| - |\tilde{\rho}_l| | \leq \left| \tilde{\rho}_k-1\right| +  \left| \tilde{\rho}_l-1\right| \leq 2 (1+\sqrt{2}) \lambda^{1/2r}.
$$
This finishes the proof.

%%%%%%%%%%%%%%%%%%%%%%%%%%%%%%%%%%%%%%%%%%%%%%%%%%%%%%%%%%%%%%
%%%%%%%%%%%%%%%%%%%%%%%%%%%%%%%%%%%%%%%%%%%%%%%%%%%%%%%%%%%%%%
\section{Proof of Lemma~\ref{lm:lambda2}}
\label{sec:ProofofRootBounds2}

The upper bound of \eqref{eq:upper_and_lower_bound_diff} follows directly from Lemma \ref{lem:BasicRootsEqns} as $|\rho - \rho'| \leq |\rho-1|+|\rho'-1| \leq 2(1+\sqrt{2})\lambda^{1/2r}$.  As for the lower bound of \eqref{eq:upper_and_lower_bound_diff}, if $\bar\rho \neq \rho'$, $\rho^{-1} \neq \rho'$, $\bar{\rho}^{-1} \ne \rho'$, then Lemma \ref{lm:lambda} gives that $| \rho-\rho'| \geq ||\rho|-|\rho'||\geq c r^{-2}\lambda^{1/2r}$.  Thus, it only remains to lower bound  $|\rho- \rho'|$ in the case where one of $\bar\rho = \rho'$, $\rho^{-1} = \rho'$, or $\bar{\rho}^{-1} = \rho'$ holds.  We will do this by producing a new (weaker) lower bound for $|\rho-\rho'|$ which does not assume anything other than that $|\rho| \neq 1$.

Let $\Lambda := \{\rho_1,..., \rho_{2r}\}$ be the set of roots of $p(x)$ in \eqref{eq:charPoly}, and note that
\begin{align}\label{eq:rho_minus_prime}
|\rho - \rho'| = \frac{\left|\prod\limits_{\rho_t \in \Lambda \setminus \{ \rho \}} (\rho-\rho_t) \right|}{ \left| \prod\limits_{\rho_t \in \Lambda \setminus \{ \rho, \rho' \}} (\rho-\rho_t) \right|},
\end{align}
where both the denominator and numerator are strictly positive by Lemma~\ref{lem:Uniquerootsforp}.
Consequently, we can bound $|\rho-\rho'|$ by appropriately bounding the numerator and denominator in \eqref{eq:rho_minus_prime}. 

To deal with the numerator in \eqref{eq:rho_minus_prime}, note that \eqref{eq:charPoly} implies that the leading coefficient of $p(x)$ is $1$, and hence
\begin{align}
\frac{1}{2r}\prod\limits_{\rho_t \in \Lambda \setminus \{ \rho \}} (\rho-\rho_t)
& = \left.\frac{1}{2r} \frac{dp(z)}{dz} \right|_{z=\rho} \nonumber\\
&= (\rho-1)^{2r-1}-(-1)^r \frac{1}{2} \lambda \rho^{r-1} \nonumber \\
& = \frac{(-1)^r  \lambda \rho^{r} }{\rho-1}-(-1)^r \frac{1}{2} \lambda \rho^{r-1}  \nonumber\\
&=\frac{(-1)^r}{2}  \lambda\rho^{r-1} \frac{\rho+1}{\rho-1}.  \nonumber
\end{align}
where the third equality uses the fact that $p(\rho)=0$ and that $\rho \neq 1$, i.e., that
$$(\rho-1)^{2r-1} = \frac{(-1)^r\lambda\rho^r}{\rho-1}.$$
Consequently,
\be \label{eq:lamb}
\left|\prod\limits_{\rho_t \in \Lambda \setminus \{ \rho \}} (\rho-\rho_t)\right| =  r \lambda \left|\rho^{r-1} \frac{\rho+1}{\rho-1}\right | > 0.
\eeq

For the denominator of \eqref{eq:rho_minus_prime}, recall that for each $\rho_t \in \Lambda, \rho_t\neq \rho$, we have by Lemma \ref{lem:BasicRootsEqns} that
 $$|\rho-\rho_t| \leq  |\rho-1|+|\rho_t-1| \leq 2(1+\sqrt{2}) \lambda^{1/2r}.$$ 
As a result, we can recombine the numerator and denominator of \eqref{eq:rho_minus_prime}  and then invoke \eqref{eq:rhoBound1} of Lemma \ref{lem:BasicRootsEqns} once more to obtain
\begin{align*}
|\rho-\rho' | 
\geq \frac{r \lambda \left|\rho^{r-1} \frac{\rho+1}{\rho-1}\right |} { 2^{2r-2}(1+\sqrt{2})^{2r-2} \lambda^{(r-1)/r}}  
 \geq   \frac {r \left | \rho^{r-1} (\rho+1) \right | }  { 2^{2r-2}(1+\sqrt{2})^{2r-1} }  \lambda^{1/2r}.
\end{align*}
Appealing one additional time to Lemma~\ref{lem:BasicRootsEqns} we now get that
\be
|\rho-\rho' | \geq \frac {r \left | (\rho+1) \right | }  { 2^{2r-2}(1+\sqrt{2})^{4r-3} }  \lambda^{1/2r}.
\label{eq:VenusFlyTrap2}
\eeq

Continuing with our bound, we will now finishing controling the numerator on the right hand side of \eqref{eq:VenusFlyTrap2}. For any $\rho \in \Lambda$ and its associated $c_k$ (defined in \eqref{eq:ck}) we have 
\begin{align} \label{eq:lambdaplusone}
| \rho+ 1| & = |(\sqrt {\rho}+i)(\sqrt {\rho}-i)|=  \left| \frac{c_k+ 2i \pm \sqrt{c_k^2+4}}{2}\cdot\frac{c_k-2i \pm \sqrt{c_k^2+4}}{2} \right| \notag \\ 
& = \left|\frac{-4+2c_ki}{c_k+2i \mp \sqrt{c_k^2+4}} \right|   \left|\frac{-4-2c_ki}{c_k-2i \mp \sqrt{c_k^2+4}} \right| \geq \frac{4 |4+c_k^2|}{64},
\end{align}
where in the inequality we used the fact that $|c_k|<2.$
It remains to obtain a lower bound on the magnitude of $|c_k^2+4|$ above.  We have that
 \begin{align}
 |4+c_k^2| & = \sqrt{ (4 - \lambda^{1/r}\cos (2k\pi/r))^2 + (\lambda^{1/r}\sin (2k\pi/r))^2} \notag \\
 &  = \sqrt{ 16+\lambda^{2/r} -8 \lambda^{1/r}\cos (2k\pi/r)} \notag  \\
 & = \sqrt{ 16- 16 \cos ^2(2k\pi/r)+ (\lambda^{1/r} -4\cos (2k\pi/r))^2} \notag \\ 
 & \geq \sqrt{ 16- 16\cos ^2(2k\pi/r)} \notag \\
  & = 4 \left| \sin (2k\pi/r) \right|. \label{eq:ck1}
 \end{align}
 
 Here we note that $|\rho| \neq 1$ implies that $k = 1, \dots, r-1$ by Lemma~\ref{lem:BasicRootsEqns}.  Furthermore, if $k = r/2$ then $c^2_k = \lambda^{1/r}$ so that
 $|4+c_k^2 | >4$ in that case.  Hence, we may combine \eqref{eq:lambdaplusone} and \eqref{eq:ck1} to see that 
 \begin{align}\label{eq:rho_plus_1}
|\rho+1| \geq \frac{1}{4}\sin(\pi/r) 
\end{align}
holds for all $r \geq 2$.
Combining \eqref{eq:rho_plus_1} with \eqref{eq:VenusFlyTrap2} now yields the desired result where here $c_1 \ge 4(1+\sqrt{2})^4$ and $c_2 \ge \frac{1}{4}r\sin(\pi/r)$ for all $r \ge 2$ (note that since $\pi \ge \frac{1}{4}r\sin(\pi/r)$ for all $r \ge 2$, we can choose $c_2 = \pi/4$, for example).

%%%%%%%%%%%%%%%%%%%%%%%%%%%%%%%%%%%%%%%%%%%%%%%%%%%%%%%%%%%%%%
%%%%%%%%%%%%%%%%%%%%%%%%%%%%%%%%%%%%%%%%%%%%%%%%%%%%%%%%%%%%%%
\section{Proof of Lemma~\ref{lem:coeff}}
\label{sec:Proofoflem:coeff}
\input{appendix_D.tex}

\end{appendix}
\end{document}

%% file: Intro.tex
For a twice differentiable real valued function on $\R^d$, the Laplace operator (or Laplacian) is a second-order differential operator given, in Cartesian coordinates, by $\Delta f  = \sum\limits_{i=1}^d \frac{\partial^2 f} {\partial x_i^2}.$  Perhaps one of the simplest and most well known properties of the Laplace operator is that in the continuous setting of univariate functions on the unit interval, its  eigenfunctions take the explicit form of sinusoidal functions. For example, with the homogenous Dirichlet boundary condition
\begin{align}\label{eq:diri}
\left\{\begin{matrix}
\Delta u_n(x) = -\lambda_n u_n(x),   \\
u_n(0) = u_n(1) = 0,
 \end{matrix}\right.
\end{align}
we have eigenvalues $\lambda_n = n^2\pi^2 $, and eigenfunctions $u_n = \sin(n\pi x)$, $n \in \mathbb{Z}_+$.
Similarly, replacing the above Dirichlet boundary condition with a homogenous mixed boundary condition
%\begin{align}
%\left\{\begin{matrix}
%\Delta u_n(x) = \lambda_n u_n(x),   \\
\begin{equation}u'_n(0) = u_n(1) = 0,\label{eq:mix}\end{equation}
% \end{matrix}\right.
%\end{align}
we have $\lambda_n =\frac{(2n-1)^2\pi^2}{4} $, and $u_n =\sqrt{2}\cos\left(\frac{(2n-1)\pi x}{2}\right) $.
Higher order Laplace operators in the same setting have similar eigen-decompositions. Let $\Delta_r$ be the $r^{th}$ order Laplacian obtained via $r$ successive applications of the Laplace operator, then with the homogenous Dirichlet boundary condition
\begin{align}
\left\{\begin{matrix}
\Delta_r u_n(x) = (-1)^r\lambda_n u_n(x),   \\
u_n(0) = u_n(1) = 0
 \end{matrix}\right.
\end{align}
the eigenfunctions are identical to those for \eqref{eq:diri} while the eigenvalues are simply raised to the $r^{th}$ power, so that $\lambda_n =n^r\pi^r $.

These examples indicate that eigenfunctions of the continuous Laplacian operator, and its higher order counterparts, have \emph{well-spread energy}. In other words, each eigenfunction is \emph{not sparse} and its support is not concentrated in any region of the domain.  Our main interest in this paper is to explore whether, and to what extent, this property is preserved after discretization. Despite its apparent simplicity, this question turns out to be highly non-trivial. We answer this question affirmatively for a specific family of discretizations of $\Delta_r$ of use in signal processing applications. We believe the proof technique used in this paper can be generalized to show similar results for many other high-order finite difference matrices with various boundary conditions. The specific discretizations we shall focus on correspond to a robin-boundary condition, which naturally arises in at least two different scenarios that motivated this work.  Both of these scenarios are related to the quantization and encoding of finite dimensional vectors, and are discussed in detail in Section \ref{sec:application} below. Given that discretizations of the Laplace transform are prevalent in many applications, we hope our basic approach will also be of broader interest. 

Let us now consider the specific discrete problem we are interested in. To that end, define the (bi-diagonal) \textit{difference matrix}, $D \in \mathbbm{R}^{N \times N}$, by its entries
\begin{equation}
D_{i,j} := \left\{ \begin{array}{ll} 1 & \textrm{if}~i=j\\ -1 & \textrm{if}~i = j+1\\ 0 & \textrm{otherwise}\end{array} \right.,
\label{Def:diffMat}
\end{equation}
and note that $D^T D$ can be viewed as a discretized Laplacian, while for integers $r\geq 2$ the matrices $D^r$ are the higher order discretized Laplacians we are interested in (see, e.g., Section \ref{sec:application}). More specifically, for an integer $r\geq 1$, consider the singular value decomposition of $D^r=U\Sigma V^*$ where $U$ and $V$ are orthonormal matrices and $\Sigma$ is a non-negative diagonal matrix. Our goal, ultimately, is to control the $\ell^\infty$ norm of the singular vectors $\v_j$ (resp. $\u_j$), which form the columns of $V$ (resp. $U$). An equivalent version of the question, which we consider herein, is to bound the $\ell^\infty$ norm of the eigenvectors of $(D^r)^T D^r = V\Sigma^2V^T$. 

A few observations are in order to help illustrate the challenge at hand (see also Section~\ref{sec:proofcommentary} below). First, when $r=1$, the problem is relatively easy and the singular values and vectors admit analytic expressions taking the form of simple trigonometric functions with $\|\u_j\|_\infty \approx \|\v_j\|_\infty \approx N^{-1/2}$ (see, e.g., \cite{von1941distribution}). This suggests that a direct approach to the problem when $r\geq 2$ might work out easily, but unfortunately that is not the case. The fundamental issue that arises is that, e.g.,  
\begin{align}\label{eq: D_ineq}(D^TD)^r \neq (D^r)^T D^r,
\end{align}
so that the matrices on either side of the inequality admit different eigen-decompositions. This is in contrast to the continuous case we saw earlier, where the eigenfunctions of higher order operators are preserved and the eigenvalues are simply those of the first order operator raised to the $r$th power. 

While inequality \eqref{eq: D_ineq} holds, it is also true that  the difference $(D^TD)^r - (D^r)^T D^r$ appears relatively well behaved in the sense that it is low-rank and sparse, which gives us hope that we may be able to appeal to matrix perturbation analysis to control the eigen-decomposition of $(D^r)^T D^r$ in terms of that of $(D^TD)^r$. Indeed by appealing to Weyl's inequalities, \cite{Sinan} (see Lemma \ref{lm:sigD}) was thus able to control the \emph{eigenvalues} of $(D^r)^T D^r$. The eigenvectors turn out to be a different matter entirely. The standard approach to eigenvector perturbation is to appeal to some version of the Davis-Kahan $\sin(\Theta)$ theorem \cite{davis1970rotation} (see also \cite{dopico2000note}). In brief, such theorems state that if the (Hermitian) matrix $\widetilde{M}=M+H$ is a perturbation of $M$ by $H$, the subspace spanned by certain eigenvectors of $\widetilde{M}$ is close to the analogous subspace spanned by eigenvectors of ${M}$, provided $\|H\|$ is small compared to the gap between the eigenvalues of $M$ and $\widetilde{M}$. Unfortunately, in our case,  \cite{Sinan} (see Lemma \ref{lm:sigD}) shows that the eigenvalues are quite close to each other, so appealing to Davis-Kahan theorems yields vacuous bounds. In recent years, similar problems in different settings have led to various results yielding eigenvector perturbation bounds, for example when the matrix $H$ is random and the matrix $M$ admits some structure (see, e.g., \cite{fan2018eigenvector, vu2011singular, o2018random, eldridge2018unperturbed}). 

In contrast with most such works, we must work with deterministic perturbations that are very large in norm compared to the minimal spectral gap herein, and thus our approach to obtaining eigenvector bounds is different. Indeed, applying preexisting results to our setting also yields vacuous bounds. Thus, most of the paper is dedicated to proving the following result via a more direct approach.

\begin{theorem}[Main Result]\label{thm:main0}
Suppose that $r \ge 2$, and let $\sigma_j := \sigma_j \left( D^r \right)$ have associated left and right singular vectors $\mathbf{u}_j, \mathbf{v}_j \in \mathbbm{R}^N$ for all $j \in [N]$. There exists absolute constants $C,C_3 > 0$ such that if $N \ge C_3^r$, we have $\max_{j \in [N]} \left\{ \|\mathbf{u}_j\|_{\infty},\|\mathbf{v}_j\|_\infty \right\} \le \frac{(Cr)^{6r}}{\sqrt{N}}$.
\end{theorem}

The proof of Theorem~\ref{thm:main0}, while utilizing relatively elementary techniques, is highly nontrivial.  In addition, we note here that while our analysis is specialized to the case of $D^r$ for the particular finite difference matrix $D$ defined in \eqref{Def:diffMat}, much of the proof technique can also be adapted to other higher order finite difference matrices that implicitly involve different boundary conditions.  We refer the reader to the next section for a more detailed overview of the proof, and to Section~\ref{sec:application} for a discussion regarding why these specific finite difference matrices are so important in some applications.  The actual proof of Theorem~\ref{thm:main0} is then given in Section~\ref{sec:Proof}, with Section~\ref{sec:EqnsForEigenvector}, Section~\ref{sec:proof_BIG_lemma}, and the appendices devoted to the proofs of supporting lemmas.

\subsection{Some Comments on the Proof of Theorem~\ref{thm:main0}}
\label{sec:proofcommentary}

As the reader may have already noticed, the proof of Theorem~\ref{thm:main0} is quite long.  Given this preexisting condition we believe it is appropriate to extend the paper's length just a bit more to try to explain why the proof is so lengthy, and why one can not prove the main result more quickly using, e.g., powerful general purpose perturbation results.  In order to get some intuition for the difficulties involved in bounding the entries of the singular vectors of our difference matrices it is helpful to look at a small example.
For instance, if $r=2$ and $N = 7$ the matrix $(D^2)^TD^2$ whose eigenvectors we must consider becomes  
$$(D^2)^TD^2=
 \left( \begin{array} {ccccccc}
%\begin{ailgn*}
%  6& -4&  1&  0& 0&  0&  0&  0&  0&  0 \\ 
 %-4&  6& -4&  1&  0&  0&  0&  0&  0&  0\\
  %1& -4&  6& -4&  1&  0&  0&  0&  0&  0\\
  %0&  1& -4& 
 6& -4&  1&  0&  0&  0&  0\\
%  0&  0&  1& 
-4&  6& -4&  1&  0&  0&  0\\
% 0&  0&  0&  
1& -4&  6& -4&  1&  0&  0\\
% 0&  0&  0&  
0&  1& -4&  6& -4&  1&  0\\
% 0&  0&  0& 
 0&  0&  1& -4&  6& -4&  1\\
%  0&  0&  0&
  0&  0&  0&  1& -4&  5& -2\\
%  0&  0&  0&
  0&  0&  0&  0&  1& -2&  1
%\end{align*}
\end{array} \right) .$$

Though $(D^2)^TD^2$ above clearly has a lot of nice structure, it is important to notice that it is not quite, e.g., Toeplitz.  Furthermore, by appealing to interlacing results for the eigenvalues of $(D^2)^TD^2$ one can see after some computation that the spectral gaps between neighboring eigenvalues of this matrix are small (on the order of $N^{-4}$ for the smallest eigenvalues).  As a result, the smaller spectral gaps between neighboring eigenvalues tend to go to $0$ relatively rapidly as $N$ grows, making them exceedingly small with respect to the size of the minimal perturbation needed to make $(D^2)^TD^2$, e.g., circulant, or Toeplitz.  Similarly, the smallest eigenvalue gaps of the closest circulant/Toeplitz matrices to $(D^2)^TD^2$ tend to be quite small as well, also going to $0$ polynomially in $1/N$ as $N$ grows.  The upshot is that standard eigenspace perturbation methods such as \cite{davis1970rotation, dopico2000note, fan2018eigenvector, eldridge2018unperturbed} do not appear to yield meaningful bounds on the $\ell^\infty$-norms of the eigenvectors in the setting of Theorem~\ref{thm:main0}. 

We find ourselves in a similar situation if we apply the singular vector perturbation theory to the asymmetric matrix $D^2$. Note that $D^2$ has a singular value gap on the order of $\mathcal{O}(N^{-2})$. If we denote this singular value gap by $\delta$, then state of the art singular vector perturbation results (see e.g. \cite{lyu2020exact}) would bound the distance between the singular vectors of $D^2$ and those of its closest circulant matrices by $O(N^{-1/2}\delta^{-1}) \sim O(N^{1.5})$ if measured in the $\ell^2$-norm, and by  $O(N^{-1}\delta^{-1}) \sim O(N)$ if measured in the $\ell^\infty$ norm.  Crucially, both of these upper bounds blow up as $N\rightarrow \infty$. Additionally, the situation only appears to get worse for $D^r$ if $r$ is chosen to be larger than $2$.

Due to these complications, and inspired by the bravery of, e.g., Strang \cite{strang1999discrete} and B\"{o}ttcher et al. \cite{bottcher2010structure} in more directly assaulting similar eigenvector problems involving related matrices, we follow their example herein.  More specifically, similar to, e.g., \cite{bottcher2010structure} we effectively treat $(D^r)^TD^r$ as a banded Hermitian Toeplitz matrix $H_{\rm Toep}$ with a structured perturbation in its lower right corner. In order to understand both the structure of the eigenvectors of $H_{\rm Toep}$, as well as the general structure of the perturbation $P := H_{\rm Toep} - (D^r)^TD^r$, in Section~\ref{sec:EqnsForEigenvector} we embed each eigenvector of $(D^r)^TD^r$ into the solution of a simple difference equation with prescribed boundary conditions corresponding to $P$.  We then solve this difference equation in order to obtain a formula for the entries of each eigenvector $\v$ of $(D^r)^TD^r$ of the form
\begin{equation} \label{equ:intoFormuladiscuss}
v_i = \sum^{2r}_{k = 1} c'_{k} \cdot \rho^i_{k},
\end{equation}
where the formula parameters $c'_1,\dots, c'_{2r}, \rho_1, \dots, \rho_{2r} \in \mathbbm{C}$ all depend on the (unknown, but bounded) eigenvalue corresponding to $\v$.

The desired $\ell^{\infty}$-bounds having already been obtained for all eigenvectors associated with eigenvalues below a prescribed cutoff in Section~\ref{sec:Proof}, the vast majority of the proof of Theorem~\ref{thm:main0} then involves using \eqref{equ:intoFormuladiscuss} to bound the $\ell^{\infty}$-norms of the eigenvectors associated with eigenvalues above the cutoff.  This portion of the proof is carried out in several phases.  First, the formula parameters $\rho_1, \dots, \rho_{2r}$ in \eqref{equ:intoFormuladiscuss} are bounded (above, below, and away from one another) in Section~\ref{sec:RootsResCharpoly}.  Next, in Section~\ref{sec:proof_BIG_lemma}, the formula parameters $c'_1,\dots, c'_{2r}$ are upper bounded both individually, and in combination with specific powers of their corresponding $\rho_k$ parameters. These upper bounds are quite delicate and involve bounding the solutions to an $\mathcal{O}(r) \times \mathcal{O}(r)$ Vandermonde system coming from the boundary conditions corresponding to $P$.  Finally, the bounds on each $v_i$ are then established using \eqref{equ:intoFormuladiscuss} in combination with the derived bounds on the $c'_k$ and $\rho_k$ parameters.  We refer the reader to Section~\ref{sec:Proof} below for additional details and discussion.

\section{Some Applications of Theorem~\ref{thm:main0} in Signal Processing}\label{sec:application}

Discretized versions of the Laplace operator play an important role in various applications. These applications include numerical analysis, where discrete Laplacians appear as finite-difference approximations of the (continuous) Laplacian operator, and image processing, where they are used, for example, in edge detection. Via discretizations of the Laplace-Beltrami operator, they are also important in various applications involving geometry, including mesh parametrization (see, e.g., \cite{gu2010discrete}). %\cite{Discrete Laplace-Beltrami Operator Determines Discrete Riemannian Metric}. 
Herein, we focus on two applications that are related to quantization schemes in signal processing, and that both specifically benefit from Theorem \ref{thm:main0}.

\subsection{Error Bounds for Sigma-Delta Quantization}
In various signal acquisition systems ranging from classical ones related to audio and image acquisition, to more recent ones like compressed sensing \cite{candes2006stable,donoho2006compressed} (see also \cite{boufounosJKS14}), continuum valued samples of signals need to be converted to digital bits. In order to reduce the quantization error, various quantization schemes have been developed including Memoryless Scalar Quantization (e.g., \cite{goyal1995quantization}), Sigma-Delta Quantization (e.g., \cite{daubechies2003approximating, gunturk2003one, krahmer2012root, Sinan}) and Beta Encoding \cite{chou2016distributed}, among others \cite{plan2013one}. In particular, the Sigma-Delta quantization family has received much attention as it generally enjoys both hardware simplicity and favorable error bounds. Nevertheless, its induced error bounds under certain signal types and measurement conditions are not entirely understood in part due to a lack of bounds on the singular vectors of $D^r$ as considered herein.

Consider, for instance, the problem of measuring and digitizing a signal modeled as a finite dimensional vector $\x$ in $\mathcal{B}^d$, the Euclidean unit ball of $\mathbbm{R}^d$, whose entries $x_i$ are infinite precision real numbers that are, e.g., potentially irrational (see, e.g., \cite{benedetto2006sigma, blum2010sobolev, powell2013quantization}).  In various settings, one can model measurements of $\x$ as inner products with frame vectors $\mathbf{f}_i \in \R^d$, $i=1,...,N$ with $N\geq d$. Thus, one has $\y=F\x \in \R^N$, where the rows of the $N\times d$ matrix $F$ are the vectors $\mathbf{f}_i$.  Having obtained $\y$, one must digitize it, i.e., replace its entries by elements of a finite alphabet $\mathcal{A}$, e.g., $\mathcal{A}=\{\pm 1\}$ in order to store it, transmit it, or process it on digital devices.  To that end, consider a particular family of quantization schemes $\mathcal{Q}: \R^N \to \mathcal{A}^N$, known as Sigma-Delta $\left( \Sigma \Delta \right)$ quantization schemes. The simplest  such scheme is the \textit{first order $\Sigma\Delta$ quantizer} which works as follows. 

Given $\y=F\x$,  one computes a vector ${\bf q} \in \{ -1, 1 \}^N$ via the following recursion with initial condition $u_0=0$:  
%\begin{equation}
%{\bf y} = F {\bf x},~~
%u_0 = 0,
%\label{equ:FirstOrdInit}
%\end{equation}
\begin{equation}
q_i = \textrm{sign} \left( y_i + u_{i-1} \right), %\textrm{and}
\label{equ:FirstOrdq}
\end{equation}
\begin{equation}
u_i =  y_i + u_{i-1} - q_i
\label{equ:FirstOrdu}
\end{equation}
for $i = 1, 2, \dots N$. 
We may succinctly restate the relationships between the vectors ${\bf x}$, ${\bf u}$, and ${\bf q}$  as%$\eqref{equ:FirstOrdInit} - \eqref{equ:FirstOrdu}$ as
\begin{equation*}
D {\bf u} ~=~ F {\bf x} - {\bf q},
\end{equation*}
%Furthermore, a short induction argument shows that $|u_i| \leq 1$ for all $i \in [N]$ provided that $|y_i| \leq 1$ for all $i \in [N]$. 
where the matrix $D$ is exactly the one defined in \eqref{Def:diffMat}.
Higher order versions of the above quantization scheme also exist, and in fact often yield better reconstruction errors.  With stable higher order schemes, equations \eqref{equ:FirstOrdq} and \eqref{equ:FirstOrdu} are modified so that now\begin{equation}D^r \u = F\x -\q, \label{equ:LinSigDelt} \end{equation} where $\|\u\|_\infty$ is bounded independently of $N$ \cite{benedetto2006sigma, blum2010sobolev, krahmer2012root, powell2013quantization}. 
With such an approach, one has represented the underlying vector $\x\in \mathcal{B}^d$ by  $N$ bits (when $\mathcal{A} = \{\pm 1\}$). 

One way to compress (i.e., encode) this $N$-bit representation without compromising reconstruction accuracy is explored in \cite{Iwen2013}. The approach in \cite{Iwen2013} capitalizes on the potential redundancy in the measurements that is due to having $N\geq d$. To encode $\q$, one simply applies the map $\q \mapsto RD^{-r} \q$, where $R\in \{0,1\}^{m\times N}$ with $m<N$ is a \emph{random selector matrix} with exactly one nonzero entry per row, which is selected uniformly
at random.  Observing that $RD^{-r}\q$ can be represented by $\sim rm\log(N)$ bits, we see that this representation can be quite parsimonious compared to storing all of $\q$ when $m\ll N$, provided we can still recover $\x$ faithfully using only $m$ entries from $D^{-r}\q$.  Towards such a reconstruction, rearranging \eqref{equ:LinSigDelt} and applying $R$ yields
\( R\u = RD^{-r}F\x - RD^{-r} q,\) which upon further manipulation yields
\begin{equation}
 \|(RD^{-r}F)^\dagger R\u\|_2 = \| \x - (RD^{-r}F)^\dagger RD^{-r} \q\|_2. \label{eq: rec_error1}\end{equation}
That is, the reconstruction error associated with the above encoding and the decoding $RD^{-r} \q \mapsto (RD^{-r}F)^\dagger RD^{-r}\q$ can be controlled by 
$\|(RD^{-r}F)^\dagger R\u\|_2 \leq \|(RD^{-r}F)^\dagger\| \|R\u\|_2$. This bound is small provided the matrix $RD^{-r}F$ has large singular values.  

When $r=1$, \cite{Iwen2013} shows that when the columns of $F$ are the $d$ singular vectors of $D$ associated with its smallest singular values, the reconstruction error \eqref{eq: rec_error1} decays exponentially  in the number of bits used for the encoding. Extending this  to $r>1$ was left as an open problem, with the stumbling block being the lack of a bound on the singular vectors of $D^{r}$ of the form $\|\v_j\|_\infty \leq \frac{C(r)}{\sqrt{N}}$, which is Theorem \ref{thm:main0} of this paper. Consequently, via minor modifications in the proof of Theorem 3 of \cite{Iwen2013} combined with Theorem \ref{thm:main0} we obtain the following result.
\begin{theorem}
Let $\epsilon,p \in (0,1)$, and $R \in \{ 0,1\}^{m \times N}$ be a random %sampling
selector matrix.  Then, there is a matrix $F \in \mathbbm{R}^{N \times d}$ such that 
$\mathcal{D}:=\left\| {\bf x} - \left(RD^{-r}F \right)^{\dagger} R D^{-r} {\bf q} \right\|_2 \leq C(\epsilon,r,d) N^{-r}$
for all ${\bf x} \in \mathcal{B}^d$ with probability at least $1-p$, when $m \geq C'(r) \epsilon^{-2} d \ln(2d / p)$. Here, $\mathbf q$ is the output of a stable $r$th order $\sd$ quantization scheme applied to $F\x$.  Furthermore, $R D^{-r} {\bf q}$ can  be encoded using $\mathcal{R} \leq m ( r\log_2 N +1 )$ bits.
\label{thm:ResExp}
\end{theorem}

From the above theorem one obtains the relationship $\mathcal{D}(\mathcal{R}) \lesssim \exp \left(-\frac{\mathcal{R}}{C'' d\log(d/p)} \right)$ between the bit-rate $\mathcal{R}$ used to encode $\x$ and the associated error (or distortion) $\mathcal{D}$, where $C''$ depends on $\epsilon, r$, and $d$. {\it Note that in the above application, it is entirely impossible to replace $D^r$, and hence $(D^r)^T D^r$, by easier to analyze matrices (i.e., matrices with nicer boundary conditions) such as $(D^T D)^r$ for two reasons.} The first is that the algorithm naturally works with the matrix $D^r$, and more importantly, the boundary conditions associated with this choice are entirely imposed by the causal nature of the quantization algorithm. {\it Quantizing the current measurement can only depend on past measurements, and not on future ones, so $D^r$ in \eqref{equ:LinSigDelt} (or any realistic substitute) simply must be a lower triangular matrix. }

Theorem~\ref{thm:ResExp} above uses Theorem~\ref{thm:main0} to prove the existence of a matrix $F$ which can be used to compactly and accurately quantize arbitrary vectors.  However, in many applications the matrix $F$ in \eqref{equ:LinSigDelt} is determined by the application, and is not something that the designer of the quantizer can choose however they like.  Thankfully in such settings there are still general formulas that bound the quantization error for a large class of general matrices $F$ (see, e.g., \cite{wang2018sigma}), but their theoretical application again depends on understanding the structure of the singular vectors of $D^r$.   For example, consider the following proposition.
%
%\RSnote{I think we should not look at the relative error in the proposition below (and in the rest of the section), and we should instead stipulate that $\|x\|_2 \leq 1$ as Rongrong does in her paper. The reason is that we're using a finite alphabet quantizer in equations (2.1) and (2.2) so $\|Fx\|_\infty$ needs to be less than 1 which can be guaranteed by the ``normalized rows" conditions and the unit ball condition (and Cauchy-Schwartz). Also, I would edit ``where $\ell$..." to "for any $\ell with$..." to avoid giving the wrong impression that this is one specific fixed $\ell$.} 
\begin{pro}(\cite{wang2018sigma} Proposition 2.3)\label{pro:1} Let $F$ be an $N\times d$ matrix with normalized rows. Then, there exists a decoder such that for any $\x \in \textrm{column-span}(F) \cap \mathcal{B}^d$, the reconstruction $\hat{\x}$ from the $r$-th order Sigma-Delta quantization of $\x$ using this decoder obeys
\[ \| \hat{\x} - \x\|_2 \lesssim \left( \frac{N}{\ell}\right)^{-r}\frac{\sqrt N}{\sigma_{\min}(V_{r,N,\ell}^TF)}
\] 
for any $\ell$ with $d \leq \ell \leq N$.  Here $V_{r,N,\ell} \in \mathbbm{R}^{N \times \ell}$ contains the $\ell$ least significant left singular vectors of the $N \times N$ $r^{\rm th}$ order difference matrix $D^{r}$. Here $\lesssim$ hides a constant independent of $N,l$ and $r$.
\end{pro}

Note that the matrix $V_{r,N,\ell}$ in Proposition~\ref{pro:1} corresponds to the objects of interest in this paper.  As a result, it should not be surprising that Theorem~\ref{thm:main0} can be used in combination with Proposition~\ref{pro:1} in order to make its upper bound on the error $\| \hat{\x} - \x\|_2$ more explicit.  For example, in the critically important case of bandlimited signal quantization via sampling, one can assume that $F$ contains columns of an $N \times N$ DFT matrix.  If we further assume that entries from the bandlimited signal $\x$ are subsampled randomly, then $F$ becomes the composition of an $N \times N$ Discrete Fourier Transform (DFT) matrix with a random sampling matrix $R \in \mathbbm{R}^{m \times N}$ containing exactly one $1$ in each row (in an i.i.d. uniformly random position).  That is, $F = R\tilde{F}$ holds, where $\tilde{F}$ now denotes a full $N \times N$ DFT matrix.  In this setting the following probabilistic lower bound is known for the smallest nonzero singular value $\sigma_{\min}(V_{r,m,\ell}^TF) = \sigma_{\min}(V_{r,m,\ell}^TR\tilde{F})$ appearing in Proposition~\ref{pro:1}.

%\RSnote{There are some differences between the theorem below and how it appears in Rongrong's paper, unless you mean Theorem 2.7 in Rongrong's paper. }
\begin{theorem}(\cite{wang2018sigma} Theorem 2.7 and Theorem 2.8)\label{pro:2} Let $\tilde{F}$ consist of $d$ columns of the $N\times N$ DFT matrix, and $R$ be the operator that randomly samples $m$ rows from a matrix.  Suppose the $\ell$ in Proposition \ref{pro:1} also satisfies  $m/\pi^2 \geq \ell \geq \frac{cm \|V_{r,m,{\ell}}\|_{\max}^2}{{\eta^2}}d\log^3(m/p)$, then with probability at least $1-p$, it holds that
\[
\sigma_{\min}(V_{r,m,\ell}^T R\tilde{F}) \geq (1-\eta)\sqrt \ell,
\]
where $V_{r,m,\ell}$ are the $\ell$ least significant left singular vectors of the $m \times m$ $r^{\rm th}$ order difference matrix $D^{r}$ for all $\ell \leq m$. Here $c$ is an absolute constant.
\end{theorem}

Combining Proposition \ref{pro:1} and Theorem \ref{pro:2} and setting $\ell$ as its allowable lower bound implies the following quantization reconstruction error %estimate depending the elementary-wise energy of singular vectors of the $r$th order difference matrix $D^r$,
\[
\| \hat{\x} - \x\|_2 \lesssim \left( \|V_{r,m,\ell}\|_{\max}^2 d \log^3(m/\epsilon) \right)^{r-\frac{1}{2}}.
\]
This brings the need to control  $\|V_{r,m,\ell}\|_{\max}$ well enough to guarantee that the bound decays for fixed $r$ as $m$ increases. Our Theorem \ref{thm:main0} (stated as a conjecture in \cite{wang2018sigma}) addresses this issue, and gives rise to the polynomially decaying bound in $m$
\[
\| \hat{\x} - \x\|_2 \lesssim \left(\frac{ C(r) d \log^3(m/\epsilon)}{m} \right)^{r-\frac{1}{2}},
\]
where $C(r)$ is a constant that depends on $r$.
This same type of improvement of related error bounds in, e.g., \cite{wang2018sigma} can also be derived for other signal types (i.e., choices of $\tilde{F}$ above) by using our Theorem \ref{thm:main0} together with other existing analogs of Proposition \ref{pro:1} and Theorem \ref{pro:2} proven therein.

%% file: ProofofBIGlemma.tex
Our main goal in this section is to prove Lemma \ref{lm:large} (here restated using notation from Section \ref{sec:EqnsForEigenvector}). 
\newtheorem*{lm:large}{\textbf{\emph{Lemma \ref{lm:large}}}}
\begin{lm:large}
There exist absolute universal constants $C_0, C_2, C_3 \in \mathbbm{R}^+$ (namely, for $C_2,C_3$ the same as in Lemma \ref{lem:coeff} below) such that for all $r \geq 2$ and $N \geq C_3^r$,
$\|\v_j\|_{\infty}\leq \frac{(C_0r)^{4r-5}}{\sqrt{N}}$ holds for all $j \in [N]$ with $\lambda_j^{1/2r}  \geq \frac{C_2 r^6}{N}$.
\end{lm:large}

Proving this lemma will require recalling several results from Section~\ref{sec:EqnsForEigenvector}.  In particular, from \eqref{eq:RRecRel1} and \eqref{equ:RecSoln} we know that the entries of each eigenvector $\v_j$ satisfy
\begin{equation}
\left(\v_{j}\right)_i = \sum^{1}_{\ell=0} \sum^{r-1}_{k = 0} c_{k,\ell} \cdot \rho^i_{k,\ell}
\label{equ:FormulaForEigs}
\end{equation}
for all $i,j \in [N]$, where the $\rho_{k,\ell}$ above are the roots from Lemma~\ref{lem:charPolyProps1} of the polynomial $p$ in \eqref{eq:charPoly} with $\lambda = \lambda_j$.
To prove Lemma~\ref{lm:large} we will use Lemmas~\ref{lem:coeffsmall} and~\ref{lem:coefflarge} below about the coefficients $c_{k,\ell}$ and roots $\rho_{k,\ell}$ in \eqref{equ:FormulaForEigs}.  These lemmas will then allow us to bound the magnitude of each entry of $\v_j$ via \eqref{equ:FormulaForEigs}.

\begin{lm}\label{lem:coeffsmall} 
There exist absolute universal constants $C'_0, C_2, C_3 \in \mathbbm{R}^+$ (for $C_2,C_3$ the same as in Lemma \ref{lem:coeff} below) such that for all $r \geq 2$, $N \geq C_3^r$, and $\lambda_j^{1/2r}  \geq \frac{C_2 r^6}{N}$,
$$\left|c_{k,\ell} \right| \leq  \frac{(C'_0r)^{2r-3}}{\sqrt N}$$
for all $(k,\ell) \in [r-1] \times \{0,1\}$.
\end{lm}

\begin{lm}\label{lem:coefflarge}
There exist absolute universal constants $C''_0, C_2, C_3 \in \mathbbm{R}^+$ (for $C_2,C_3$ the same as in Lemma \ref{lem:coeff} below) such that for all $r \geq 2$, $N \geq C_3^r$, and $\lambda_j^{1/2r}  \geq \frac{C_2 r^6}{N}$,$$\left|\rho^{N+1-r}_{k,\ell}c_{k,\ell} \right| \leq \frac{(C''_0r)^{4r-6}}{\sqrt{N}} \quad \forall |\rho_{k,\ell}|>1.$$
\end{lm}

The proof of Lemma \ref{lm:large} is a fairly simple consequence of \eqref{equ:FormulaForEigs} given Lemmas~\ref{lem:coeffsmall} and~\ref{lem:coefflarge}.

\begin{proof}[\textbf{Proof of Lemma \ref{lm:large}}]
Let $i,j \in [N]$.  Below we will reorder the $2r$ roots $\left\{ \rho_{k,\ell} \right\}_{\ell,k}$ in \eqref{equ:FormulaForEigs} by magnitude so that the resulting reordered sequence
$\rho_l := \rho_{k_l,\ell_l}$ satisfies 
\begin{equation}
|\rho_1| \leq |\rho_2| \leq \dots < |\rho_r| = 1 = |\rho_{r+1}| < |\rho_{r+2}| \leq \dots \leq |\rho_{2r}|.
\label{equ:RootOrder}
\end{equation}
Note that Lemmas~\ref{lem:charPolyProps1} and~\ref{lem:BasicRootsEqns} guarantee that such an ordering of the roots exists.  Similarly, we will also reorder the roots' associated coefficients $c_l := c_{k_l,\ell_l}$ correspondingly so that the resulting reordered sum in \eqref{equ:FormulaForEigs} still satisfies
$$\left(\v_{j}\right)_i = \sum^{2r}_{l = 1} c_{l} \cdot \rho^i_{l}$$
for all $i \in [N]$.
 
We can now see that
$$\left|\left(\v_{j}\right)_i \right| = \left| \sum^{2r}_{l = 1} c_{l} \cdot \rho^i_{l} \right| \leq \sum\limits_{l=1} ^{2r}|c_{l} \rho_{l}^{i}| \leq \sum\limits_{l=1} ^{r+1} |c_l| +\sum\limits_{l=r+2} ^{2r}|c_{l} \rho_{l}^{N}| $$
where our final inequality uses the properties of the ordering of the roots in \eqref{equ:RootOrder}.  Rearranging this last expression and using Lemma~\ref{lem:BasicRootsEqns} now reveals that
$$\left|\left(\v_{j}\right)_i \right| \leq \sum\limits_{l=1} ^{r+1} |c_l| + \left|\rho_{l}\right|^{r-1} \sum\limits_{l=r+2} ^{2r} \left|c_{l} \rho_{l}^{N+1-r} \right| \leq \sum\limits_{j=1} ^{r+1} |c_{l}| +(1+\sqrt{2})^{2r-2}\left(\sum\limits_{l=r+2} ^{2r} \left|c_{l} \rho_{l}^{N+1-r} \right|\right).$$
Employing Lemmas \ref{lem:coeffsmall} and \ref{lem:coefflarge} (since $r \geq 2$, $N \geq C_3^r$, and $\lambda_j^{1/2r}  \geq \frac{C_2 r^6}{N}$ by assumption, their conditions are met), we can now see that
\begin{align*}
\left|\left(\v_{j}\right)_i \right| &\leq (r+1)\frac{(C'_0r)^{2r-3}}{\sqrt N} + (r-1)(1+\sqrt{2})^{2r-2}\frac{(C''_0r)^{4r-6}}{\sqrt{N}} \nonumber \\ 
&\leq (r+1) (1+\sqrt{2})^{2r-2}\frac{(\max \left(C'_0, C''_0,1\right)r)^{4r-6}}{\sqrt{N}} \nonumber \\
&\leq \frac{(C_0r)^{4r-5}}{\sqrt{N}}. 
\end{align*}
where $C_0$ is an absolute constant chosen such that $C_0^{4r-5} \ge (\max \left(C'_0, C''_0,1\right))^{4r-6}(1+\sqrt{2})^{2r-2} \frac{r+1}{r}$ for all $r \ge 2$.
\end{proof}

We will devote the remainder of this section to proving Lemmas \ref{lem:coeffsmall} and \ref{lem:coefflarge}.  In order to do so we will need several supporting results.

\subsection{Supporting Lemmas}

First, we will require the following result about the inverse of a Vandermonde matrix in several places below.

\begin{lm} \label{lm:vandermond}
Suppose $A$ is a Vandermonde matrix
\[
A= \left[ \begin{matrix}1 & x_1& x_1^2 & \cdots & x_1^{n-1} \\ 1 & x_2& x_2^2 & \cdots & x_2^{n-1} \\
\vdots &&&& \vdots\\
1 & x_n& x_n^2 & \cdots & x_n^{n-1}  \end{matrix}\right].
\]
Then $A^{-1}=U^{-1}L^{-1}$ with 
\[
(L^{-1})_{i,j}=\left\{\begin{aligned}& 0, & i<j \\
&1, & i = j = 1\\
& \prod \limits_{k=1,k\neq j}^i \frac{1}{x_j-x_k}, & \textrm{otherwise}
\end{aligned}\right..  
\]
That is,
\[
L^{-1}= \left[\begin{matrix}1&0&0& \cdots 
\\ \frac{1}{x_1-x_2} & \frac{1}{x_2-x_1}&0&\cdots 
\\ \frac{1}{(x_1-x_2)(x_1-x_3)} & \frac{1}{(x_2-x_1)(x_2-x_3)} & \frac{1}{(x_3-x_1)(x_3-x_2)} & \cdots
\\ \vdots & \vdots & \vdots &\vdots \end{matrix}\right].
\]
Moreover,
\[
(U^{-1})_{i,j} = \left\{\begin{aligned}& 1, & i=j \\& 0, & j=1, i \ne j \\& (U^{-1})_{i-1,j-1}-(U^{-1})_{i,j-1}x_{j-1} & \textrm{otherwise}\end{aligned}\right.
\]
where $(U^{-1})_{0,j}$ is considered to be $0$ for the purposes of the recursion. That is,
\[
U^{-1}=\left[\begin{matrix}1 & -x_1 & x_1x_2 & -x_1x_2x_3 & \cdots \\
0&1&-(x_1+x_2)&x_1x_2+x_2x_3+x_3x_1& \cdots \\
0&0&1&-(x_1+x_2+x_3)& \cdots \\
0&0&0&1&\cdots \\
\vdots &\vdots &\vdots&\vdots&\vdots \end{matrix}\right].
\]
This recurrence is equivalent to the following closed form expression for the entries of $U^{-1}$:
\begin{equation}\label{eq:Uclosedform}
(U^{-1})_{i,j} = (-1)^{i+j} \sum_{1 \le a_1 < \cdots < a_{j-i} \le j-1} x_{a_1}\cdots x_{a_{j-i}},
\end{equation}
for all $i \le j$, where the empty sum is defined to be equal to $1$, and 
\[(U^{-1})_{i,j}=0\]
otherwise.
\end{lm}
\begin{proof}
The first two results concerning the entries of $U^{-1}$ and $L^{-1}$ are proven in \cite{Richard}. We prove the third result concerning the closed form expression for the entries of $U^{-1}$ by induction. We first show that for $i>j$, $(U^{-1})_{i,j}=0$. By our earlier result, $(U^{-1})_{i,j} = 0$ for $j=1$, $i \ne j$, hence the result holds for our base case $j=1$. Then, suppose the result holds for all such entries in columns $1,\cdots,j$, and suppose that $i>{j+1}$. Then, we have $(U^{-1})_{i,j+1} = (U^{-1})_{i-1,j}-(U^{-1})_{i,j}x_{j} = 0-0x_j=0$, since $(U^{-1})_{i-1,j} = (U^{-1})_{i,j} = 0$ as $i>j+1$ implies $i-1>j, i>j$. Thus, the desired result holds by induction.

Next, we show the result holds for $i \le j$ by induction. We first see that the result holds for the base case $(U^{-1})_{1,1}$ due to the empty sum being defined as $1$. Now, suppose the result holds for all entries in columns $1,\cdots, j$. First, suppose that $i<j+1$; we then have
\begin{align*}
(U^{-1})_{i,j+1} &= (U^{-1})_{i-1,j}-(U^{-1})_{i,j}x_{j} \\
&= (-1)^{i+j-1} \sum_{1 \le a_1 < \cdots < a_{j-i+1} \le j-1} x_{a_1}\cdots x_{a_{j-i+1}} - x_j\left[(-1)^{i+j} \sum_{1 \le a_1 < \cdots <  a_{j-i} \le j-1} x_{a_1}\cdots x_{a_{j-i}}\right]\\
&= (-1)^{i+j+1} \sum_{1 \le a_1 < \cdots < a_{j-i+1} \le j-1} x_{a_1}\cdots x_{a_{j-i+1}} + (-1)^{i+j+1} \sum_{1 \le a_1 < \cdots <  a_{j-i} \le j-1} x_{a_1}\cdots x_{a_{j-i}}x_j \\
&= (-1)^{i+j+1} \sum_{1 \le a_1 < \cdots < a_{j-i+1} \le j} x_{a_1}\cdots x_{a_{j-i+1}}
\end{align*}
where the second equality holds since $i < j+1$ implies $i-1 \le j, i \le j$, and where the last equality holds since the first term consists of the sum of all products of $j-i+1$ terms consisting of variables indexed in the range $[j-1]$ while the second term consists of the sum of all products of $j-i+1$ terms with variables indexed in the range $[j]$ which contain $x_j$. 

Then, suppose that $i=j+1$. We see that the sum in the right hand side of \eqref{eq:Uclosedform} is the empty sum, hence, it suffices to show that $(U^{-1})_{i,j+1}=1$. We have $(U^{-1})_{i,j+1} = (U^{-1})_{i-1,j}-(U^{-1})_{i,j}x_{j} = 1 - 0x_j = 1$, where $(U^{-1})_{i-1,j}=1$ by our inductive hypothesis and the fact that the sum in the right hand side of \eqref{eq:Uclosedform} is empty, and where $(U^{-1})_{i,j}=0$ by our previous case, since $i=j+1$ implies $i>j$. Thus, the desired result holds by induction. 
\end{proof}

Next, the following result bounds the magnitudes of the coefficients $c_{k,\ell}$ in \eqref{equ:FormulaForEigs} associated with the $r+1$ largest-magnitude roots $\rho_{k,\ell}$ of the polynomial $p$ in \eqref{eq:charPoly} whenever $\lambda = \lambda_j$ is sufficiently large.  It is used to prove Lemma \ref{lem:coeffsmall}.

\begin{lm}\label{lem:coeff} There exist absolute uniform constants $C_1, C_2, C_3 \in \mathbbm{R}^+$ such that for all $r \geq 2$, $N \geq C_3^r$, and $\lambda_j^{1/2r} \geq \frac{C_2 r^6}{N}$,
 \begin{equation} |c_{k,\ell}| \leq  \frac{C_1}{\sqrt{N}} \textrm{  if } |\rho_{k,\ell}| \geq  1,
 \label{equ:bound_large_root}
 \end{equation}
 for all $(k,\ell) \in [r-1] \times \{0,1\}$. 
\end{lm}

In particular, to prove this result, we will prove the following two lemmas, from which Lemma \ref{lem:coeff} immediately follows:
\begin{lm} \label{lem:coeffpart1} There exists absolute uniform constant $C_2 \in \mathbb{R}^+$ such that for all $r \ge 2$ and $1/4 \geq \lambda_j^{1/2r} \geq \frac{C_2r^6}{N}$, 
\[|c_{k,\ell}| \leq \left(\frac{24}{ \min\{|\rho_{k,j}|^{2N},1\} N}\right)^{1/2},\] 
 for all $(k,\ell) \in [r-1] \times \{0,1\}$.
\end{lm}

\begin{lm}\label{lem:coeffpart2}
There exists absolute uniform constant $C_3 \in \mathbb{R}^+$, such that for all $r \ge 2$, $N \geq C_3^{r}$, and $\lambda^{1/2r}>1/4$, 
\[
| c_{k,\ell}| \leq \left(\frac{48}{ \min\{|\rho_{k,j}|^{2N},1\} N}\right)^{1/2},
\]
 for all $(k,\ell) \in [r-1] \times \{0,1\}$.
\end{lm}
\noindent The proofs of Lemmas \ref{lem:coeffpart1} and \ref{lem:coeffpart2} are rather involved, and so have been moved to Appendix~\ref{sec:Proofoflem:coeff}.

Using Lemmas \ref{lm:vandermond} and \ref{lem:coeff}, we now prove Lemmas \ref{lem:coeffsmall} and \ref{lem:coefflarge}, which completes the proof of our main lemma, Lemma \ref{lm:large}. To prove Lemma \ref{lem:coeffsmall} we use the bound on the coefficients $c_{k,\ell}$ corresponding to roots with $ |\rho_{k,\ell}| \geq  1$, and use the boundary conditions \eqref{eq:RRecRel1} to extend this bound to a bound which holds for all $c_{k,\ell}$.

\begin{proof}[\textbf{Proof of Lemma \ref{lem:coeffsmall}}]

As in (\ref{equ:RootOrder}) in the proof of Lemma \ref{lm:large}, we reorder the roots $\rho_{k,\ell}$ in the following way
\[|\rho_1| \leq |\rho_2| \leq \dots < |\rho_r| = 1 = |\rho_{r+1}| < |\rho_{r+2}| \leq \dots \leq |\rho_{2r}|,\]
and similarly rearrange the associated coefficients $c_{k,\ell}$ such that 
\[\left[\tilde{v}(\lambda_j)\right]_i = \sum^{2r}_{l = 1} c_{l} \cdot \rho^i_{l}.\]
 By Lemma \ref{lem:coeff}, since $r \geq 2$, $N \geq C_3^r$, and $\lambda_j^{1/2r}  \geq \frac{C_2 r^6}{N}$ by assumption, the coefficients $|c_{l}| \le \frac{C_1}{\sqrt{N}}$ for $l = r,\cdots,2r$ are bounded by $\frac{C_1}{\sqrt{N}}$, since $|\rho_l| \ge 1$ for these roots.

To establish similar bounds for the remaining $c_{l}$, we need to use the boundary conditions of $\tilde{v}$ from \eqref{eq:RRecRel1}, namely the fact that $\tilde{v}_i =0$ for $1-r\leq i\leq 0$. This is equivalent to 
\begin{equation}
[\tilde{v}(\lambda_j)]_{i} = \sum_{l=1}^{2r} c_{l} \cdot \rho^i_{l} =0, \textrm{ for } i=1-r,...,0.
\label{eq:bnd_first_r}
\end{equation}
We can then rewrite this equation in terms of matrices in the following way:
\begin{equation}\label{eq:matrix_bnd}
\left[\begin{matrix} 1 & \cdots &1 \\ \rho_{1} & \cdots &\rho_{r-1}\\ \vdots &\ddots& \vdots \\
\rho_{1}^{r-2}& \cdots & \rho_{r-1}^{r-2} \end{matrix} \right]  \left( \begin{matrix} \rho_{1}^{1-r} c_{1} \\ \rho_{2}^{1-r}c_{2} \\ \vdots \\ \rho_{r-1}^{1-r} c_{r-1}\end{matrix} \right) = -\left[\begin{matrix} 1 & \cdots &1 \\  \rho_{r} & \cdots & \rho_{2r}\\ \vdots &\ddots& \vdots \\
\rho_{r}^{r-2}& \cdots & \rho_{2r}^{r-2} \end{matrix} \right]  \left( \begin{matrix} \rho_{r}^{1-r} c_{r} \\ \rho_{r+1}^{1-r}c_{r+1} \\ \vdots \\ \rho_{2r}^{1-r} c_{2r}\end{matrix} \right).
\end{equation}
By multiplying both sides of \eqref{eq:matrix_bnd} with $H\in \mathbb{R}^{(r-1)\times (r-1)}$ defined as
\[H_{i,j} = \left\{ 
\begin{array}{ll}
(-1)^{i-j}\binom{i-1}{j-1} &\mbox{for } i\geq j\\
0 & \mbox{for } i<j
\end{array}\right. ,\]
we have 
\begin{equation}\label{eq:bd1}
A_1  \mathbf{c}_1 = -A_2 \mathbf{c}_2
\end{equation}
where $A_1 \in \mathbb{C}^{(r-1) \times (r-1)}$, $A_2 \in \mathbb{C}^{(r-1) \times (r+1)}$. In particular,
\[
\mathbf{c}_1 = \left( \begin{matrix} c_{1} \\ \vdots \\ c_{r-1}\end{matrix} \right) \in \mathbb{C}^{r-1}, \quad \quad \mathbf{c}_2 =\left( \begin{matrix} c_{r} \\ \vdots \\ c_{2r}\end{matrix} \right)\in \mathbb{C}^{r+1}
\]
\[
A_1 = \left[\begin{matrix} 1 & \cdots &1 \\ \rho_{1}-1 & \cdots & \rho_{r-1}-1\\ \vdots &\ddots &\vdots \\
(\rho_{1}-1)^{r-2}& \cdots & (\rho_{r-1}-1)^{r-2} \end{matrix} \right] \left[ \begin{matrix} \rho_{1}^{1-r} & 0 & \cdots & 0 \\  0 & \rho_{2}^{1-r} & \ddots & \vdots \\ \vdots & \ddots & \ddots & 0 \\ 0 & \cdots & 0& \rho_{r-1}^{1-r} \end{matrix} \right] , 
\]
and
\[ A_2 =  \left[\begin{matrix} 1 & \cdots &1 \\  \rho_{r}-1 & \cdots & \rho_{2r}-1\\ \vdots &\ddots &\vdots \\
( \rho_{r}-1)^{r-2}& \cdots & ( \rho_{2r}-1)^{r-2} \end{matrix} \right]\left[ \begin{matrix}  \rho_{r}^{1-r} & 0 & \cdots & 0 \\  0 &  \rho_{r+1}^{1-r} & \ddots & \vdots \\ \vdots & \ddots & \ddots & 0 \\ 0 & \cdots & 0&   \rho_{2r}^{1-r} \end{matrix} \right],
\]
since the $(i,j)$th entry of 
\[H\left[\begin{matrix} 1 & \cdots &1 \\ \rho_{1} & \cdots &\rho_{r-1}\\ \vdots &\ddots& \vdots \\
\rho_{1}^{r-2}& \cdots & \rho_{r-1}^{r-2} \end{matrix} \right]\]
is 
\[\sum_{l=1}^{r-1} H_{i,l}\rho_{j}^{l-1} = \sum_{l=1}^{i} (-1)^{i-l}\binom{i-1}{l-1} \rho_{j}^{l-1} = \sum_{l=0}^{i-1} (-1)^{i-l-1}\binom{i-1}{l} \rho_{j}^{l} = (\rho_{j}-1)^{i-1},\]
and similarly for the right hand side of \eqref{eq:matrix_bnd}.

From \eqref{eq:rhoBound1} in Lemma \ref{lem:BasicRootsEqns}, we know that $| \rho_{l}-1 | \leq (1+\sqrt 2)\lambda^{1/2r}$ for all $l$ and from \eqref{equ:RootOrder} we have $|\rho_{l}| \geq 1$, for any $l \in \{r,\cdots,2r\}$. These therefore imply that for $i \in \{1,\cdots,r-1\}, j \in \{1,\cdots,r+1\}$,
  \begin{align} \label{eq:abound}
|(A_2)_{i,j}| &\leq [(1+\sqrt 2)\lambda^{1/2r}]^{i-1} \nonumber \\ 
&\le (1+\sqrt{2})^{r-2} (\lambda^{1/2r})^{i-1} \nonumber \\
&\le C_{11}^r(\lambda^{1/2r})^{i-1}, 
\end{align}
for $C_{11}>0$ an absolute constant.
Recall that our goal is to bound $\mathbf{c}_1 = - A_1^{-1}A_2\mathbf{c}_2$, so we next seek to bound the operator norm of $A_1^{-1}$ by bounding its entries.

We first note that by \eqref{equ:RootOrder}, $|\rho_{1}|, \cdots, |\rho_{r-1}| < 1$, and therefore $|\rho_{1}^{r-1}|, \cdots |\rho_{r-1}^{r-1}| \le 1$. Hence, to bound the entries of $A_1^{-1}$, it suffices to bound the entries of the inverse of the Vandermonde matrix
\[V := \left[\begin{matrix} 1 & \cdots &1 \\ \rho_{1}-1 & \cdots & \rho_{r-1}-1\\ \vdots &\ddots &\vdots \\
(\rho_{1}-1)^{r-2}& \cdots & (\rho_{r-1}-1)^{r-2} \end{matrix} \right]\]
since the entries of $A_1^{-1}$ will therefore not increase in norm if the inverse diagonal matrix is included. Let $V^T=LU$ be the $LU$ decomposition of $V^T$. By  Lemma \ref{lm:vandermond}, we have
\[
L^{-1}= \left[\begin{matrix}1&0&0& \cdots 
\\ \frac{1}{\rho_{1}- \rho_{2}} & \frac{1}{\rho_{2}-\rho_{1}}&0&\cdots 
\\ \frac{1}{(\rho_{1}-\rho_{2})(\rho_{1}-\rho_{3})} & \frac{1}{(\rho_{2}-\rho_{1})(\rho_{2}-\rho_{3})} & \frac{1}{(\rho_{3}-\rho_{1})(\rho_{3}-\rho_{2})} & \cdots
\\ \vdots & \vdots & \vdots &\vdots \end{matrix}\right].
\]

To bound the entries of $L^{-1}$, recall that from \eqref{eq:upper_and_lower_bound_diff} of Lemma \ref{lm:lambda2}, we know that
$\left| \frac{1}{\rho_{i}-\rho_{j}} \right| \leq t^{-1}(r) \lambda^{-1/2r}$ for any $i\neq j$ and $|\rho_{i}|, |\rho_{j}|\neq 1$, where $t(r) = c_2c_1^{-r}$ in the case of conjugate or inverse roots, and $t(r)=cr^{-2}$ otherwise. As a result, 
\begin{align}\label{eq:Lbound} 
|(L^{-T})_{i,j}| = |(L^{-1})_{j,i}| &= \prod_{k=1, k \ne i}^j \left|\frac{1}{\rho_{i}-\rho_{k}}\right|\nonumber\\
&\leq  \prod_{k=1, k \ne i}^j t^{-1}(r)\lambda^{-1/2r}\nonumber\\
&\leq c_3^{r}(cr^{-2})^{2-j}(\lambda^{-1/2r})^{j-1},
\end{align}
where the last inequality holds because at most one of the pairs $\rho_{i},\rho_{k}$ can be conjugate (since $\rho_{i}$ is held fixed in each term), and none of the pairs can be inverses or conjugate inverses of each other since $|\rho_l| < 1$ for $l \in \{1,\cdots,r-1\}$.  Here $c_3$ is chosen such that $c_3^{r} \ge c_2^{-1}c_1^{r}, c_3^{r} \ge (cr^{-2})^{-1}$ for all $r\ge2$.
Similarly, we can compute the entrywise bound for $U^{-T}$ from its explicit expression derived in Lemma \ref{lm:vandermond} 
 \begin{align}\label{eq:Ubound}  |(U^{-T})_{i,j} | = |(U^{-1})_{j,i}| &= \left|(-1)^{i+j}\sum_{1 \le a_{1} < \cdots < a_{i-j}\le i-1} (\rho_{a_1}-1) \cdots (\rho_{a_{i-j}}-1)\right| \nonumber\\
 &\le \binom{i-1}{i-j} [(1+\sqrt{2}) \lambda^{1/2r}]^{i-j} \nonumber\\
 & \leq 2^r (C_{12}\lambda^{1/2r})^{i-j}.  
\end{align}
for $C_{12}>1$ an absolute constant, $i \ge j$ (where the last inequality holds since $\binom{i-1}{i-j} \le 2^r$ since $i,j < r$).

Since we defined $V^T=LU$, we have that $V^{-1} = L^{-T} U^{-T}$. As a consequence 
\begin{align}\label{eq:first_r_bnd}
|(A_1^{-1}A_2)_{i,j}| &\le \sum_{k=1}^{r-1} \sum_{l=1}^{r-1} |(L^{-T})_{i,k}| |(U^{-T})_{k,l} ||(A_2)_{l,j}| \nonumber \\ 
&\le \sum_{k=1}^{r-1} \sum_{l=1}^{r-1} c_3^{r}(cr^{-2})^{2-k}(\lambda^{-1/2r})^{k-1}2^r (C_{12}\lambda^{1/2r})^{k-l} C_{11}^r(\lambda^{1/2r})^{l-1} \nonumber \\
&\le r^{2r-6}\max(c,c^{3-r}) (2c_3C_{11}C_{12})^r \sum_{k=1}^{r-1} \sum_{l=1}^{r-1} (\lambda^{-1/2r})^{k-1}(\lambda^{1/2r})^{k-l} (\lambda^{1/2r})^{l-1} \nonumber \\ 
&= r^{2r-6}\max(c,c^{3-r}) (2c_3C_{11}C_{12})^r \sum_{k=1}^{r-1} \sum_{l=1}^{r-1} (\lambda^{1/2r})^{(1-k)+(k-l)+(l-1)} \nonumber \\ 
&\le r^{2r-4}\max(c,c^{3-r}) (2c_3C_{11}C_{12})^r \nonumber \\
&\le (C'r)^{2r-4} 
\end{align}
where the first inequality holds because $|(A_1^{-1})_{ij}|\le |(V^{-1})_{ij}|$ for all $i,j$, the second inequality used \eqref{eq:abound}, \eqref{eq:Lbound}  and \eqref{eq:Ubound}, the third inequality follows since $(cr^{-2})^{2-k} \le r^{2r-6}\max(c,c^{3-r})$ since $k \in \{1,\cdots,r-1\}$, and similarly $C_{12}^{k-l} \le C_{12}^{r}$ since $k, l \in \{1,\cdots,r-1\}$ with $k \ge l$ and $C_{12} > 1$. The fourth inequality holds because $\sum_{k=1}^{r-1} \sum_{l=1}^{r-1} 1 = (r-1)^2 \le r^2$, and in the last inequality, $C'>0$ is an absolute constant chosen so that $\max(c,c^{3-r}) (2c_3C_{11}C_{12})^r \le C^{2r-4}$ for all $r \ge 2$. 

Thus, since $\mathbf{c}_1 = - A_1^{-1}A_2\mathbf{c}_2$, \eqref{eq:first_r_bnd} implies that 
\begin{equation*}
\|\mathbf{c}_1\|_{\infty} \leq (C'r)^{2r-4}  (r+1)\|\mathbf{c}_2\|_{\infty} \le (C'r)^{2r-4}  (r+1) \frac{C_1}{\sqrt N} \leq \frac{(C'_0r)^{2r-3}}{\sqrt N}, 
\end{equation*}
where we used $\|\mathbf{c}_2\|_{\infty} \leq \frac{C_1}{\sqrt N}$ (by Lemma \ref{lem:coeff}, and since $|\rho_l| \ge 1$ for $l \in \{r,\cdots,2r\}$ by (\ref{equ:RootOrder})).  Here $C'_0 > 0$ is an absolute constant chosen so that $(C'_0)^{2r-3} \ge (C')^{2r-4}C_1\frac{r+1}{r}$ and $(C'_0)^{2r-3} \ge C_1$ for all $r \ge 2$ (this second condition is so that $\frac{(C'_0r)^{2r-3}}{\sqrt{N}} \ge \frac{C_1}{\sqrt{N}}$; thus the desired bound will hold for all choices of $(k,\ell)$).
\end{proof}

We will now use the bound on $|c_{k,\ell}|$ just proven in Lemma \ref{lem:coeffsmall} to prove Lemma \ref{lem:coefflarge}. Similar to the proof of \ref{lem:coeffsmall}, we extend the bound on a subset of the roots to all of the roots by using the boundary conditions \eqref{eq:RRecRel1BC}.

\begin{proof}[\textbf{Proof of Lemma \ref{lem:coefflarge}}] 

As in \eqref{equ:RootOrder} in Lemma \ref{lm:large}, we reorder the roots $\rho_{k,\ell}$ in the following way
\[|\rho_1| \leq |\rho_2| \leq \dots < |\rho_r| = 1 = |\rho_{r+1}| < |\rho_{r+2}| \leq \dots \leq |\rho_{2r}|,\]
and similarly rearrange the associated coefficients $c_{k,\ell}$ such that 
\[\left[\tilde{v}(\lambda_j)\right]_i = \sum^{2r}_{l = 1} c_{l} \cdot \rho^i_{l}.\]
We first note that
\begin{equation}\label{eq:smallrootcrhobound}
|c_{l}\rho_{l}^i| \leq |c_{l}| \leq \frac{(C'_0r)^{2r-3}}{\sqrt N}
\end{equation}
 for $j=1,...,r+1$, by Lemma \ref{lem:coeffsmall} and the fact that $|\rho_{l}| \le 1$ by \eqref{equ:RootOrder}. We then follow a similar argument as the one used in the proof of Lemma \ref{lem:coeffsmall},  to extend this upper bound of $|c_{l}\rho_{l}^i|$ to $l=r+2,...,2r$ by using the last $r$ boundary conditions  \eqref{eq:RRecRel1BC} and the general solution to the recurrence \eqref{equ:RecSoln}.

Using \eqref{eq:RRecRel1BC} we obtain
\[
0= \sum\limits_{q=0}^r \binom{r}{q}(-1)^q [\tilde{v}(\lambda_j)]_{k-q} = \sum\limits_{q=0}^r \binom{r}{q}(-1)^q  \left( \sum\limits_{l=1}^{2r}  c_{l} \rho_{l}^{k-q}\right), \textrm{ for } k= N+1,...,N+r.
\] 
Factoring $\rho_l^{k-q} = \rho_l^{k-r}\rho_l^{r-q}$ and exchanging the two summations above we obtain
\[
0=  \sum\limits_{l=1}^{2r} c_{l} \rho_{l}^{k-r} \sum\limits_{q=0}^r (-1)^q \binom{r}{q}   \rho_{l}^{r-q} .
\]
This is equivalent to 
\[
\sum\limits_{l=1}^{2r} c_{l} \rho_{l}^{k-r} (\rho_{l}-1)^r = 0,   \textrm{ for } k=N+1,...,N+r,
\]
or to
\[
\sum\limits_{l=r+2}^{2r} c_{l} \rho_{l}^{k-r} (\rho_{l}-1)^r = - \sum\limits_{l=1}^{r+1} c_{l} \rho_{l}^{k-r} (\rho_{l}-1)^r,   \textrm{ for } k=N+1,...,N+r.
\]

As in the proof of Lemma \ref{lem:coeffsmall}, we write the first $r-1$ of these equations in matrix form:
\begin{align*}
\left[\begin{matrix} 1 & \cdots &1 \\ \rho_{r+2} & \cdots &\rho_{2r}\\ \vdots &\ddots &\vdots \\
\rho_{r+2}^{r-2}& \cdots & \rho_{2r}^{r-2} \end{matrix} \right]\left[ \begin{matrix}( \rho_{r+2}-1)^{r} & 0 & \cdots & 0 \\  0 &(\rho_{r+3}-1)^{r} & \ddots & \vdots \\ \vdots & \ddots & \ddots & 0 \\ 0 & \cdots & 0&  (\rho_{2r}-1)^{r} \end{matrix} \right] \left( \begin{matrix}\rho_{r+2}^{N+1-r} c_{r+2} \\ \rho_{r+3}^{N+1-r} c_{r+3}\\ \vdots \\ \rho_{2r}^{N+1-r} c_{2r}\end{matrix} \right) \\ 
= - \left[\begin{matrix} 1 & \cdots &1 \\ \rho_{1} & \cdots &\rho_{r+1}\\ \vdots &\ddots &\vdots \\
\rho_{1}^{r-2}& \cdots & \rho_{r+1}^{r-2} \end{matrix} \right]\left[ \begin{matrix}(\rho_{1}-1)^{r} & 0 & \cdots & 0 \\  0 &(\rho_{2}-1)^{r} & \ddots & \vdots \\ \vdots & \ddots & \ddots & 0 \\ 0 & \cdots & 0&  (\rho_{r+1}-1)^{r} \end{matrix} \right]  \left( \begin{matrix} \rho_{1}^{N+1-r}c_{1} \\  \rho_{2}^{N+1-r}c_{2} \\ \vdots \\ \rho_{r+1}^{N+1-r}c_{r+1}\end{matrix} \right).
\end{align*}
Then, by the same argument as in the proof of Lemma \ref{lem:coeffsmall}, multiplying both sides of the resulting matrix equation by the $(r-1)\times (r-1)$ matrix
\[H_{i,j} = \left\{ 
\begin{array}{ll}
(-1)^{i-j}\binom{i}{j} &\mbox{for } i\geq j\\
0 & \mbox{for } i<j
\end{array}\right. ,\]
results in  
\begin{align*}
\overbrace{\left[\begin{matrix} 1 & \cdots &1 \\ \rho_{r+2}-1 & \cdots &\rho_{2r}-1\\ \vdots &\ddots &\vdots \\
(\rho_{r+2}-1)^{r-2}& \cdots & (\rho_{2r}-1)^{r-2} \end{matrix} \right]}^{A_3} \overbrace{\left[ \begin{matrix}(1- \rho_{r+2})^{r} & 0 & \cdots & 0 \\  0 &(1-\rho_{r+3})^{r} & \ddots & \vdots \\ \vdots & \ddots & \ddots & 0 \\ 0 & \cdots & 0&  (1-\rho_{2r})^{r} \end{matrix} \right] }^{B_3} \left( \begin{matrix}\rho_{r+2}^{N+1-r} c_{r+2} \\ \rho_{r+3}^{N+1-r} c_{r+3}\\ \vdots \\ \rho_{2r}^{N+1-r} c_{2r}\end{matrix} \right) \\ 
= - \underbrace{ \left[\begin{matrix} 1 & \cdots &1 \\ \rho_{1}-1 & \cdots &\rho_{r+1}-1\\ \vdots &\ddots &\vdots \\
(\rho_{1}-1)^{r-2}& \cdots & (\rho_{r+1}-1)^{r-2} \end{matrix} \right]}_{A_4} \underbrace{ \left[ \begin{matrix}(1- \rho_{1})^{r} & 0 & \cdots & 0 \\  0 &(1-\rho_{2})^{r} & \ddots & \vdots \\ \vdots & \ddots & \ddots & 0 \\ 0 & \cdots & 0&  (1-\rho_{r+1})^{r} \end{matrix} \right] }_{B_4} \left( \begin{matrix} \rho_{1}^{N+1-r}c_{1} \\  \rho_{2}^{N+1-r}c_{2} \\ \vdots \\ \rho_{r+1}^{N+1-r}c_{r+1}\end{matrix} \right).
\end{align*}
Since we have 
\[ \left( \begin{matrix}\rho_{r+2}^{N+1-r} a_{r+2} \\ \vdots \\ \rho_{2r}^{N+1-r} a_{2r}\end{matrix} \right)= -B_3^{-1} A_3^{-1} A_4B_4  \left( \begin{matrix} \rho_{1}^{N+1-r}a_{1} \\ \vdots \\ \rho_{r+1}^{N+1-r}a_{r+1}\end{matrix} \right) \]
we begin by bounding the entries of $A_3^{-1}A_4$.

First, we see that essentially identical arguments used to prove (\ref{eq:Lbound}) and (\ref{eq:Ubound}) and bound the entries of $V^{-1}$ in the proof of Lemma \ref{lem:coeffsmall} apply here, and result in the same bounds the on the entries of $A_3^{-1}$, since the indices of the roots are not considered in either argument. The only slight difference is that, in this case, when bounding the entries of $L^{-T}$, all of the roots have norm strictly greater than $1$ rather than strictly less than $1$, but the same argument still holds in this case. Also, by a similar argument to that used in showing (\ref{eq:abound}) we see that 
\begin{equation}\label{eq:a4bound}
(A_4)_{i,j} \le [(1+\sqrt{2})\lambda^{1/2r}]^{i-1} \le (1+\sqrt{2})^{r-2}(\lambda^{1/2r})^{i-1} \le C_{11}^r(\lambda^{1/2r})^{i-1}
\end{equation}
for $C_{11}$ the same constant in (\ref{eq:abound}).
Thus, by essentially the same argument as for (\ref{eq:first_r_bnd}), we have that
\begin{align}\label{eq:a3a4bound}
|(A_3^{-1}A_4)_{i,j}| &= \sum_{k=1}^{r-1} \sum_{l=1}^{r-1} |(L^{-T})_{i,k}| |(U^{-T})_{k,l}| |(A_4)_{l,j}| \nonumber \\ 
&\le \sum_{k=1}^{r-1} \sum_{l=1}^{r-1} c_3^{r}(cr^{-2})^{2-k}(\lambda^{-1/2r})^{k-1}2^r (C_{12}\lambda^{1/2r})^{k-l} C_{11}^r(\lambda^{1/2r})^{l-1} \nonumber \\
&\le (C'r)^{2r-4}, 
\end{align}
for the same constant $C'$.

Then, by Lemma \ref{lem:BasicRootsEqns} we have $|1-\rho_{\ell}|/|1-\rho_{l} | \leq  (1+\sqrt 2)^2$ for all $\ell,l$. Hence, we have 
\[|(-B_3^{-1} A_3^{-1} A_4B_4)_{i,j}| \le (1+\sqrt{2})^{2r}(C'r)^{2r-4}.\]
Also, recall that from \eqref{eq:smallrootcrhobound}, we have 
\[|\rho_{l}^{N+1-r}c_l| \leq \frac{(C'_0r)^{2r-3}}{\sqrt N}\]
for $l \in \{1,\cdots, r+1\}$. 
Therefore,
\begin{align}\label{equ:bound_large_rho}
 \left|\left|\left( \begin{matrix}\rho_{(r+2)}^{N+1-r} c_{(r+2)} \\ \vdots \\ \rho_{(2r)}^{N+1-r} c_{(2r)}\end{matrix} \right)\right|\right|_\infty &=  \left|\left|-B_3^{-1} A_3^{-1} A_4B_4  \left( \begin{matrix} \rho_{(1)}^{N+1-r}c_{(1)} \\ \vdots \\ \rho_{(r+1)}^{N+1-r}c_{(r+1)}\end{matrix} \right)\right|\right|_\infty \nonumber \\
 & \leq (1+\sqrt{2})^{2r}(C'r)^{2r-4} (r+1) \frac{(C'_0r)^{2r-3}}{\sqrt N} \nonumber\\ 
 & \leq \frac{(C''_0r)^{4r-6}}{\sqrt{N}}
\end{align}
for $C''_0$ an absolute constant chosen such that $(C''_0)^{4r-6} \ge (1+\sqrt{2})^{2r} (C')^{2r-4}(C'_0)^{2r-3}\frac{r+1}{r}$ for all $r \ge 2$.

\end{proof}

Having established both Lemma \ref{lem:coeffsmall} and Lemma \ref{lem:coefflarge} now finishes our proof of Lemma \ref{lm:large}.

%% file: appendix_D.tex
%%%%
%%%%   Appendix D from TheoWIPtemp.tex has been reformatted and included below
%%%%

In this appendix, we seek to prove Lemma \ref{lem:coeff}, here restated below for clarity.

\newtheorem*{lem:coeff}{\textbf{\emph{Lemma \ref{lem:coeff}}}}
\begin{lem:coeff}
%\begin{lm}\label{lem:coeff} 
There exist absolute uniform constants $C_1, C_2, C_3 \in \mathbbm{R}^+$ such that for all $r \geq 2$, $N \geq C_3^r$, and $\lambda_j^{1/2r} \geq \frac{C_2 r^6}{N}$,
 \begin{equation} |c_{k,\ell}| \leq  \frac{C_1}{\sqrt{N}} \textrm{  if } |\rho_{k,\ell}| \geq  1,
 \label{equ:bound_large_root}
 \end{equation}
 for all $(k,\ell) \in [r-1] \times \{0,1\}$. 
\end{lem:coeff}

 As mentioned in the discussion after its statement, to prove this result, we will prove separate results for $\lambda^{1/2r} \le 1/4$ and for $\lambda^{1/2r}>1/4$, namely Lemmas \ref{lem:coeffpart1} and \ref{lem:coeffpart2}, again restated below.
 
\newtheorem*{lem:coeffpart1}{\textbf{\emph{Lemma \ref{lem:coeffpart1}}}}
\begin{lem:coeffpart1}
There exists an absolute uniform constant $C_2 \in \mathbb{R}^+$ such that for all $r \ge 2$ and $1/4 \geq \lambda_j^{1/2r} \geq \frac{C_2r^6}{N}$, 
\[|c_{k,\ell}| \leq \left(\frac{24}{ \min\{|\rho_{k,j}|^{2N},1\} N}\right)^{1/2},\] 
 for all $(k,\ell) \in [r-1] \times \{0,1\}$.
\end{lem:coeffpart1}

\newtheorem*{lem:coeffpart2}{\textbf{\emph{Lemma \ref{lem:coeffpart2}}}}
\begin{lem:coeffpart2}
There exists an absolute uniform constant $C_3 \in \mathbb{R}^+$, such that for all $r \ge 2$, $N \geq C_3^{r}$, and $\lambda_j^{1/2r}>1/4$, 
\[
| c_{k,\ell}| \leq \left(\frac{48}{ \min\{|\rho_{k,j}|^{2N},1\} N}\right)^{1/2},
\]
 for all $(k,\ell) \in [r-1] \times \{0,1\}$.
\end{lem:coeffpart2}

The first subsection of this appendix is dedicated to the $\lambda_j^{1/2r}\le1/4$ case, and the second to the $\lambda_j^{1/2r} > 1/4$ case. 

We begin with a toy example that demonstrates our approach for proving Lemma \ref{lem:coeff}, namely the case of two real roots $\rho_1, \rho_2$. Although this case does not occur in practice, it motivates the main ideas used in proving Lemma \ref{lem:coeff}.

\begin{ex}\label{ex:real}[Toy example: two real roots]
Suppose $v \in \mathbb{R}^N $ is a normalized vector $\|v\|_2=1$ with the following element-wise representation $$v_i = \beta_1 \rho_1^i + \beta_2 \rho_2^i,$$ where $\beta_j, \rho_j, j=1,2$ are real numbers.  Assuming $\alpha:= \frac{| \rho_1| } { |\rho_2|} >1$,
then it holds that $$|\beta_j| \leq \frac{2}{\sqrt{N-2\log_2^{-1} \alpha-1}}, \textrm{ if } \rho_j \geq 1.  $$
\end{ex}

To prove this result, we will first prove the following result which bounds the number of indices for which two sequences bounded by geometric progressions will have terms sufficiently far apart.
 
\begin{lm}\label{lm:interval} Let $\rho_1, \rho_2, \beta_1, \beta_2$ be positive numbers with $\alpha = \frac{\rho_1}{\rho_2}>1$. Let $\{B_i\}_{i \in [N]}, \{B'_i\}_{i \in [N]}$, be sequences of positive numbers such that $\frac{B_{i+1}}{B_i} \ge \rho_1$, $\frac{B'_{i+1}}{B'_i} \le \rho_2$ for any $i \in [N]$. Then, for any $q \in \mathbb{N}$, the set of indices $\mathcal{K} \subset [N]$ for which 
\begin{equation}
\mathcal{K} = \left\{ i: \frac{B_i}{B'_i}  \leq 2^{-\frac{q}{2}\log_2 \alpha } \textrm{ or }\frac{B_i}{B'_i}  \geq  2^{\frac{q}{2}\log_2 \alpha } \right\} \notag
 \end{equation}
is of cardinality
\[ 
|\mathcal{K}| \geq N-q-1. \]
%\[
%A_i > B_i 2^{\frac{r}{2}\log \alpha } \textrm{ or } A_i <B_i 2^{-\frac{r}{2}\log \alpha}.
%\]
\end{lm}

\begin{proof}
%We prove the lemma by showing that the number of indices $i \in [N]$  that satisfy
%\[
%\frac{A_i}{B_i} \in [2^{-\frac{r}{2}\log \alpha},2^{\frac{r}{2}\log \alpha}] 
%\] 
%is at most $r$.
Let $\Gamma_i =\frac{B_i}{B'_i} $, and note that since $\alpha>1$ we have $\Gamma_{i+1} = \frac{B_{i+1}}{B'_{i+1}} \ge \frac{p_1B_i}{p_2B'_i} \ge \alpha\Gamma_i > \Gamma_i$, and so $\Gamma_i$ is strictly increasing in $i$.  We have
\begin{equation}\label{eq:gamma}
 \frac{\Gamma_i}{\Gamma_{i+q}} = \prod_{j=i}^{i+q-1} \frac{\Gamma_i}{\Gamma_{i+1}} = \prod_{j=i}^{i+q-1} \left(\frac{B_i}{B_{i+1}}\right) \left(\frac{B'_{i+1}}{B'_i}\right) \le \left(\frac{\rho_2}{\rho_1}\right)^q =  \frac{1}{\alpha^q} =\frac{1}{2^{q \log_2 \alpha}}. 
 \end{equation}
Thus, if $\Gamma_{i} \in I := [2^{-\frac{q}{2}\log_2 \alpha},2^{\frac{q}{2}\log_2 \alpha}] $, then $\Gamma_{i+q+1}$ must be outside the interval $I$. Define $i_0:= \min\{i, \Gamma_i \in I \}$. Then for any  $i \geq i_0+q+1$,  we have $\Gamma_i >\Gamma_{i_0+q} =  2^{q \log_2 \alpha}\Gamma_{i_0} \geq 2^{q/2 \log_2\alpha}$, which implies $[N]\backslash \mathcal{K} \subset \{ i_0,..., i_0+q\}$.
\end{proof}

Using this result, we can now prove the result in the example by using Lemma \ref{lm:interval} to show that for any element in $\mathcal{K}$, $\beta_1 \rho_1^i$ and $\beta_2 \rho_2^i$ differ sufficiently, and using this result to bound the $|\beta_j|$.

\begin{proof}[Proof of Example \ref{ex:real}] Let $\mathcal{K}=\{i: \frac{|\beta_1 \rho_1^i|}{|\beta_2 \rho_2^i|} \geq 2 \textrm{ or }  \frac{|\beta_1 \rho_1^i|}{|\beta_2 \rho_2^i|} \leq 1/2 \}$ be the set of indices at which the two components in the expression of $v$ are sufficiently different (differ by a factor of 2). From Lemma \ref{lm:interval} above, with $B_i = \beta_1\rho_1^i$, $B'_i = \beta_2 \rho_2^i$, and hence $\frac{B_{i+1}}{B_i} = \frac{| \beta_1\rho_1^{i+1}|}{| \beta_1\rho_1^i|} = |\rho_1|$, and similarly $\frac{B'_{i+1}}{B'_i} = \frac{| \beta_2\rho_2^{i+1}|}{ |\beta_2\rho_2^i|} = |\rho_2|$ (and thus $\alpha = \frac{|\rho_1|}{|\rho_2|}$, which we assumed was greater than $1$), the cardinality of $\mathcal{K}$ is of the same order as $N$, in particular, we have $|\mathcal{K}| \geq  N- \frac{2} {\log_2 \alpha} -1$. For each of those indices $i\in \mathcal{K}$, we have by the definition of $\mathcal{K}$ that
\[
|v_i| = | \beta_1 \rho_1^i + \beta_2 \rho_2^i| \geq \frac{1}{2} \max \{ |\beta_1 \rho_1^i|,|\beta_2 \rho_2^i|\} \geq \frac{1}{2}| \beta_j \rho_j^i|, \textrm{ for  } j=1,2.
\]
Summing up the squares of all entries of $v_i$ whose indices are in $\mathcal{K}$, we have
\begin{equation}\label{equ:up}
1\geq \sum\limits_{i\in \mathcal{K}} v_i^2 \geq \sum\limits_{i\in \mathcal{K}} \frac{1}{4} \beta_j^2\rho_j^{2i}.
\end{equation}
If $|\rho_j| \geq 1$, \eqref{equ:up} can be used to show that
\[
\sum\limits_{i\in \mathcal{K}} \frac{1}{4} \beta_j^2\rho_j^{2i} \geq \sum_{i=1}^{|\mathcal{K}|}   \frac{1}{4} \beta_j^2\rho_j^{2i} \geq \sum_{i=1}^{|\mathcal{K}|}   \frac{1}{4} \beta_j^2 =\frac{1}{4} |\mathcal{K}| \beta_j^2.
\]
Together with \eqref{equ:up}, it yields
\[
|\beta_j| \leq \frac{2}{\sqrt {|\mathcal{K}|} } \leq \frac{2}{\sqrt {N-\frac{2}{\log_2 \alpha} -1}}.
\]

\end{proof}

In this example, we were able to directly apply Lemma \ref{lm:interval} since the roots were real. In reality, however, it will not be the case that the roots will be real. In order to resolve this issue, we prove the following lemmas to resolve the issue of complex roots by writing the singular vector as a sum of real-valued terms. First, we prove the following result about the coefficients $c_{k,\ell}$ in the expansion (\ref{equ:RecSoln}).
\begin{lm}\label{lm:conjCoeffs}
For 
\begin{equation} \label{equ:RecSoln2}
[\tilde{v}(\lambda_j)]_{i} = \sum^{1}_{\ell=0} \sum^{r-1}_{k = 0} c_{k,\ell} \cdot \rho^i_{k,\ell}
\end{equation}
as in (\ref{equ:RecSoln}), $c_{0,1}=\overline{c_{0,0}}$, and $c_{k,\ell} = \overline{c_{r-k,\ell}}$ for all $k \in \{1,\cdots,r-1\}$. 
\end{lm}
\begin{proof}
Since the singular vector $\tilde{v}(\lambda_j)$ is chosen to be real (note that such a choice is always possible because $(D^r)^T(D^r)$ is symmetric), we have 

\[\sum^{1}_{\ell=0} \sum^{r-1}_{k = 0} c_{k,\ell} \cdot \rho^i_{k,\ell} = \overline{\sum^{1}_{\ell=0} \sum^{r-1}_{k = 0} c_{k,\ell} \cdot \rho^i_{k,\ell}} =\sum^{1}_{\ell=0} \sum^{r-1}_{k = 0} \overline{c_{k,\ell}} \cdot \overline{\rho_{k,\ell}}^i\]
for all $i \in \mathbb{Z}$.
Then, using the results in Lemmas \ref{lem:charPolyProps1} and \ref{lem:BasicRootsEqns}, we first note that $\rho_{0,0}=\rho_{0,1}^{-1}$ and $|\rho_{0,0}|=|\rho_{0,1}|=1$, hence $\rho_{0,1}=\overline{\rho_{0,0}}$. Then, again by Lemma \ref{lem:charPolyProps1} we also have $\rho_{k,\ell}=\overline{\rho_{r-k,\ell}}$ for all $k \in \{1,\cdots,r-1\}$. Hence, the above equation is equivalent to 

\[c_{0,0}\rho_{0,0}^i + c_{0,1}\rho_{0,1}^i+\sum_{\ell=0}^1\sum_{k=1}^{r-1} c_{k,\ell}\rho^i_{k,\ell} =\overline{c_{0,0}}\rho_{0,1}^i + \overline{c_{0,1}}\rho_{0,0}^i+\sum_{\ell=0}^1\sum_{k=1}^{r-1} \overline{c_{k,\ell}}\rho_{r-k,\ell}^i = \overline{c_{0,0}}\rho_{0,1}^i + \overline{c_{0,1}}\rho_{0,0}^i+\sum_{\ell=0}^1\sum_{k=1}^{r-1} \overline{c_{r-k,\ell}}\rho_{r,\ell}^i\]

where the last equality in the sequence holds by a change of variables. 

Then, since the equality holds for all $i \in \mathbb{Z}$, we thus have 
\[\left[\begin{matrix}1 & 1 & 1 & \cdots & 1\\\rho_{0,0} & \rho_{0,1} & \rho_{1,0} & \cdots & \rho_{r-1,1} \\ \vdots & \vdots & \vdots & \ddots & \vdots \\\rho_{0,0}^{2r-1} & \rho_{0,1}^{2r-1} & \rho_{1,0}^{2r-1} & \cdots & \rho_{r-1,1}^{2r-1}\end{matrix}\right]\left(\begin{matrix}c_{0,0} \\ c_{0,1} \\ c_{1,0} \\ \vdots \\ c_{r-1,1}\end{matrix}\right)=\left[\begin{matrix}1 & 1 & 1 & \cdots & 1\\\rho_{0,0} & \rho_{0,1} & \rho_{1,0} & \cdots & \rho_{r-1,1} \\ \vdots & \vdots & \vdots & \ddots & \vdots \\\rho_{0,0}^{2r-1} & \rho_{0,1}^{2r-1} & \rho_{1,0}^{2r-1} & \cdots & \rho_{r-1,1}^{2r-1}\end{matrix}\right]\left(\begin{matrix}\overline{c_{0,1}} \\ \overline{c_{0,0}} \\ \overline{c_{r-1,0}} \\ \vdots \\\overline{c_{1,1}}\end{matrix}\right).\]
Since 
\[V=\left[\begin{matrix}1 & 1 & 1 & \cdots & 1\\\rho_{0,0} & \rho_{0,1} & \rho_{1,0} & \cdots & \rho_{r-1,1} \\ \vdots & \vdots & \vdots & \ddots & \vdots \\\rho_{0,0}^{2r-1} & \rho_{0,1}^{2r-1} & \rho_{1,0}^{2r-1} & \cdots & \rho_{r-1,1}^{2r-1}\end{matrix}\right] \in \mathbb{C}^{2r \times 2r}\]
is a Vandermonde matrix, we hence have $\det V = \prod_{(j,\ell) \ne (j',\ell')} \left(\frac{1}{\rho_{j,\ell}-\rho_{j',\ell'}}\right)$. By Lemma \ref{lem:Uniquerootsforp}, the roots $\rho$ are distinct, and thus $V$ is invertible. Hence 
\[\left(\begin{matrix}c_{0,0} \\ c_{0,1} \\ c_{1,0} \\ \vdots \\ c_{r-1,1}\end{matrix}\right)=\left(\begin{matrix}\overline{c_{0,1}} \\ \overline{c_{0,0}} \\ \overline{c_{r-1,0}} \\ \vdots \\\overline{c_{1,1}}\end{matrix}\right)\]
and therefore $c_{0,1}=\overline{c_{0,0}}$ and $c_{k,\ell} = \overline{c_{r-k,\ell}}$ for all $k \in \{1,\cdots,r-1\}$, as desired.
\end{proof}
%Since the singular vector $\tilde{v}(\lambda)$ is real and that complex roots $\rho_{k,j}$ appear in conjugate paris, the corresponding coefficients $a_{k,j}$, $a_{r-k,j}$ (with $k=1,...,r-1$, $j=0,1$) are also conjugate to each other. 

%Using this property and the notation $\rho_{k,j} = |\rho_{k,j}| \mathbbm{e}^{\theta_{k,j}\mathbbm{i}}$, $a_{k,j} = |a_{k,j}|\mathbbm{e}^{\gamma_{k,j}\mathbbm{i}}$, we can rewrite the expansion \eqref{equ:RecSoln}  in real basis as follows.
%
%\begin{align}\label{equ:realrep}
%[\tilde{v}(\lambda)]_{i} &= Re\left[\sum^{1}_{j=0} \sum^{r-1}_{k = 0} a_{k,j} \cdot \rho^i_{k,j} \right]= \sum^{1} _{j=0} \sum^{\lfloor (r-1)/2 \rfloor}_{k = 1} 2 |a_{k,j}| \cdot |\rho_{k,j}|^i \cos (i \theta_{k,j}+\gamma_{k,j}) \notag \\ & + 2 |a_{0,0}| \cos (i \theta_{0,0}+ \gamma_{0,0}) +
%\sum^{1}_{j=0}\sum\limits_{k=\lceil r/2 \rceil}^{\lfloor r/2 \rfloor} a_{k,j} \rho_{k,j}^i
%\end{align}
%where the first and second sum on the right hand side accounts for the contribution from  complex $\rho_{k,j}$s, and the third term accounts for that from real $\rho_{k,j}$s, which only exist when $r$ is divisible by 2, according to Lemma \ref{lem:charPolyProps1}. The second term represent the contribution from the conjugate pair of the unit roots. 

Using this result, we can now prove the following lemma, which allows us to write $[\tilde{v}(\lambda_j)]_i$ (the $i$th entry of $\tilde{v}$ corresponding to $\lambda=\lambda_j$) in terms of real-valued terms. 

\begin{lm}\label{lm:realrep}
Let $\rho_{k,\ell} = |\rho_{k,\ell}| \mathbbm{e}^{\theta_{k,\ell}\mathbbm{i}}$, $c_{k,\ell} = |c_{k,\ell}|\mathbbm{e}^{\gamma_{k,\ell}\mathbbm{i}}$, for $\rho_{k,\ell}$ and $c_{k,\ell}$ as in \eqref{equ:RecSoln2}, and let 
\[\tilde{c}_{k,\ell} := \begin{cases} |c_{k,\ell}| & \text{ if } r \in 2\mathbb{Z} \text{ and } \ell=r/2 \\ 2|c_{k,\ell}| & \text{ otherwise.}\end{cases}\] Then 
\begin{equation}\label{eq:compact_expre}
[\tilde{v}(\lambda_j)]_{i}=\tilde{c}_{0,0}|\rho_{0,0}|^i\cos(i\theta_{0,0}\gamma_{0,0})+\sum^{1} _{\ell=0} \sum^{\lfloor r/2 \rfloor}_{k = 1}\tilde{c}_{k,\ell} |\rho_{k,\ell}|^i \cos (i \theta_{k,\ell}+\gamma_{k,\ell}). 
\end{equation}
\end{lm}
\begin{proof}
From the results in Lemmas \ref{lem:charPolyProps1}, \ref{lem:BasicRootsEqns}, and \ref{lm:conjCoeffs} we have $\rho_{0,1}=\overline{\rho_{0,0}}$ and $c_{0,1}=\overline{c_{0,0}}$; we also have $\rho_{k,\ell}=\overline{\rho_{r-k,\ell}}$ and $c_{k,\ell} = \overline{c_{r-k,\ell}}$ for all $k \in \{1,\cdots,r-1\}$. Hence, we see that $c_{0,1}\rho_{0,1}^i=\overline{c_{0,0}\rho_{0,0}^i}$, and $c_{k,\ell}\rho_{k,\ell}^i = \overline{c_{r-k,\ell}\rho_{r-k,\ell}^i}$ for all $k \in \{1,\cdots,r-1\}$. We also note that $c_{k,\ell} = \overline{c_{r-k,\ell}}$ implies that if $r$ is even, $c_{r/2,\ell}$ is real, and since by Corollary \ref{cor:onlyRealRoots}, $\rho_{r/2,\ell}$ is also real, we can see that $c_{r/2,\ell}\rho_{r/2,\ell}^i$ is real as well. 

Then, (for $k \ne r/2$) since $c_{k,\ell}\rho_{k,\ell}^i = \overline{c_{r-k,\ell}\rho_{r-k,\ell}^i}$, we have \[c_{k,\ell}\rho_{k,\ell}^i +c_{r-k,\ell}\rho_{r-k,\ell}^i = 2 Re(c_{k,\ell}\rho_{k,\ell}^i) = 2Re(|c_{k,\ell}||\rho_{k,\ell}|^i\mathbbm{e}^{(i\theta_{k,\ell}+\gamma_{k,\ell})\mathbbm{i}})=2|c_{k,\ell}||\rho_{k,\ell}|^i\cos(i\theta_{k,\ell}+\gamma_{k,\ell}),\]
where the second equality holds since we have assumed that $\rho_{k,\ell} = |\rho_{k,\ell}| \mathbbm{e}^{\theta_{k,\ell}\mathbbm{i}}$, $c_{k,\ell} = |c_{k,\ell}|\mathbbm{e}^{\gamma_{k,\ell}\mathbbm{i}}$. Similarly, since $c_{0,1}\rho_{0,1}^i=\overline{c_{0,0}\rho_{0,0}^i}$, we have 
\[c_{0,0}\rho_{0,0}^i+c_{0,1}\rho_{0,1}^i = Re(c_{0,0}\rho_{0,0}^i)+Re(c_{0,1}\rho_{0,1}^i) = 2|c_{0,0}||\rho_{0,0}|\cos(i \theta_{0,0}+\gamma_{0,0})\] Finally, if $r$ is even, and hence $c_{r/2,\ell}\rho_{r/2,\ell}^i$ is real, we thus have 

\[\sum_{\ell=0}^1 c_{r/2,\ell}\rho_{r/2,\ell}^i = \sum_{\ell=0}^1 Re(c_{r/2,\ell}\rho_{r/2,\ell}^i) = \sum_{\ell=0}^1 |c_{r/2,\ell}||\rho_{r/2,\ell}|^i \cos(i\theta_{r/2,\ell}+\gamma_{r/2,\ell}).\]

We then see that we can break the expansion (\ref{equ:RecSoln}) into parts as follows:
\[[\tilde{v}(\lambda_j)]_{i} = \sum^{1}_{\ell=0} \sum^{r-1}_{k = 0} c_{k,\ell} \cdot \rho^i_{k,\ell} = \sum^{1} _{\ell=0} \sum^{\lfloor (r-1)/2 \rfloor}_{k = 1} c_{k,\ell}\rho_{k,\ell}^i +c_{r-k,\ell}\rho_{r-k,\ell}^i  + c_{0,0}\rho_{0,0}^i+c_{0,1}\rho_{0,1}^i +
\sum^{1}_{\ell=0}\sum\limits_{k=\lceil r/2 \rceil}^{\lfloor r/2 \rfloor} c_{k,\ell} \rho_{k,\ell}^i.\]
Here, the last term in the sum will be empty if $r$ is odd, and $c_{r/2,\ell} \rho_{r/2,\ell}^i$ if $r$ is even. Thus, combining the three equations above, we have 
\begin{align*}[\tilde{v}(\lambda_j)]_{i} =& \sum^{1} _{\ell=0} \sum^{\lfloor (r-1)/2 \rfloor}_{k = 1} 2 |c_{k,\ell}| |\rho_{k,\ell}|^i \cos (i \theta_{k,\ell}+\gamma_{k,\ell}) +  2|c_{0,0}||\rho_{0,0}|\cos(i \theta_{0,0}+\gamma_{0,0})+ \\&
\sum^{1}_{\ell=0}\sum\limits_{k=\lceil r/2 \rceil}^{\lfloor r/2 \rfloor} |c_{k,\ell}||\rho_{k,\ell}|^i \cos(i\theta_{k,\ell}+\gamma_{k,\ell}).
\end{align*}

Thus, by our definition of $\tilde{c}_{k,\ell}$, we have 
\[[\tilde{v}(\lambda_j)]_{i} =\tilde{c}_{0,0}|\rho_{0,0}|^i\cos(i\theta_{0,0}\gamma_{0,0})+ \sum^{1} _{\ell=0} \sum^{\lfloor r/2 \rfloor}_{k = 1}\tilde{c}_{k,\ell} |\rho_{k,\ell}|^i \cos (i \theta_{k,\ell}+\gamma_{k,\ell}),\]
as desired.

\end{proof}

\subsection{The $\lambda_j^{1/2r} \leq 1/4$ case: Proof of Lemma \ref{lem:coeffpart1}}\label{subsec:leq14}

The proof of Lemma \ref{lem:coeffpart1} follows the same idea as in our toy example, Example 1. Our goal is still to find the set of indices for which the components in the expansion of $\tilde{v}$ are well separated, and to show that this set has a large enough cardinality. These were easy to prove in the toy example because each time the index increases by 1, the increments of the two components in the expression of $v$ has a ratio that is lower bounded by some positive number independent of $N$, which fulfills the assumptions of Lemma \ref{lm:interval}, leading to the desired coefficient bound. 

In the general scenario, however, the presence of complex roots and the cosine functions in the expansion (\ref{eq:compact_expre}) of $\tilde{v}$ prevent such lower bounds from existing. Indeed, we observe that it is the large oscillation of the components from index to index that destroys the lower bound, but at the same time, we observe that the average oscillation over an interval of indices is much smaller. This motivates us to look at the increments of components of $\tilde{v}$ from interval to interval instead of from index to index, with the hope of finding a positive lower bound that is sufficiently large. We first establish the intervals in the following definition by dividing $[N]$ into a number of subsets of equal length $\Delta N$.
%changed from i\DelN+1
\begin{definition}\label{def:intervals}
Given $\Delta N > 0$ and $m:=\lfloor N/\Delta N \rfloor $, let $I_0,...,I_{m-1}$ be intervals with $I_i=\{A_i,...,A_{i+1}-1\}$, where $A_i = i\Delta N+1$ for all $i \in \{0,\cdots,m\}$. 
\end{definition}

As described before, in a similar fashion to Example \ref{ex:real}, we seek to show the summands in the expansion (\ref{eq:compact_expre}) differ sufficiently from each other. However, because of the presence of $\cos (i \theta_{k,\ell}+\gamma_{k,\ell})$ in the summands, we must consider sets of summands rather than individual summands. We begin by defining the following terminology to discuss sets of summands.

\begin{definition} \label{def:appDnotation} In \eqref{eq:compact_expre}, we call each summand a \emph{component} of $\tilde{v}(\lambda_j)$, we call the squared summand at index $i$, i.e., $\tilde{c}_{k,\ell}^2 |\rho_{k,\ell}|^{2i} \cos^2(i\theta_{k,\ell}+\gamma_{k,\ell}) $ the \emph{energy of the $(k,\ell)$'th component at $i$}, and the sum $\sum_{i\in I_n}\tilde{c}_{k,\ell}^2 |\rho_{k,\ell}|^{2i} \cos^2(i\theta_{k,\ell}+\gamma_{k,\ell}) $ the \emph{energy of the $(k,\ell)$'th component over interval $I_{n}$}.
\end{definition}

We then define the following sets which make precise the comparison of the increments of components of $\tilde{v}$. 

\begin{definition} \label{def:appDC} For some fixed $(k,\ell)\neq (k',\ell')$ with $(k,\ell),(k',\ell') \in \{(0,0)\} \cup (\{1,\cdots, \lfloor r/2 \rfloor\} \times \{0,1\})$, we define the index set $C_{k,\ell; k',\ell'}$ where for each $i\in C_{k,\ell; k',\ell'}$ the energy of the $(k,\ell)$th component on $I_i$ is comparable to the $(k',\ell')$th component on $I_i$, i.e., $$C_{k,\ell; k',\ell'} =  \left\{i: \frac{1}{4r^2}\leq \frac{\tilde{c}_{k,\ell} ^2 \sum\limits_{n \in I_i} |\rho_{k,\ell}|^{2n} \cos^2 (n \theta_{k,\ell}+\gamma_{k',\ell'}) }{\tilde{c}_{k',\ell'} ^2 \sum\limits_{n \in I_{i}} |\rho_{k',\ell'}|^{2n} \cos^2 (n \theta_{k',\ell'}+\gamma_{k',\ell'}) } \leq 4r^2 \right\}.  $$ Also, we define 
\begin{equation}\label{eq:def_C}
C = \bigcap_{(k,\ell) \ne (k',\ell')} C_{k,\ell; k',\ell'}^c,
\end{equation}
where $(k,\ell), (k',\ell')$ again run over $\{(0,0)\} \cup (\{1,\cdots, \lfloor r/2 \rfloor\} \times \{0,1\})$.
\end{definition}

In the previous definition, the bound $4r^2$ is chosen out of consideration for a later argument. Also, on the complement set of $C_{k,\ell; k',\ell'}$, the $(k,\ell)$th and $(k',\ell')$th components are sufficiently different, more specifically, one component will be $4r^2$ times larger than the other. The idea is that we compute $C_{k,\ell; k',\ell'}$ for all $(k,\ell)$ and $(k',\ell')$, so that on $C$, all components will be sufficiently different from each other, and therefore there must exist a dominant component that behaves similarly to their sum over all of the components, which is the square of the $\ell^2$ norm of the singular vector of $D^r$ we are considering. The coefficient of this dominant component can in turn be bounded by the energy of the singular vector on the corresponding interval. The claim is, if we have enough intervals in $C$, then we can expect a small bound for the coefficient (Lemma \ref{lem:coeffpart1}). To compute $|C|$ we first find a lower bound for $|C_{k,\ell;k',\ell'}^c|$, for each $k,\ell;k',\ell'$ quadruple. 

We begin by supposing that $|\rho_{k,\ell}|, |\rho_{k',\ell'}|>1$.
In this case, define 
\begin{equation}\label{equ:est_C}
D_{k,\ell;k',\ell'} =\left\{i: \frac{\widetilde{c}_{k,\ell}^2 |\rho_{k,\ell}|^{2A_{i+1}}  \sum\limits_{n\in I_i} \cos^2(n\theta_{k,\ell} +\gamma_{k,\ell})} {\widetilde{c}_{k',\ell'}^2 |\rho_{k',\ell'}|^{2A_i} \sum\limits_{n\in I_i}  \cos^2(n\theta_{k',\ell'} +\gamma_{k',\ell'})} \leq \frac{1} {4r^2}, \textrm{ or} \frac{\widetilde{c}_{k,\ell}^2 |\rho_{k,\ell}|^{2A_{i}}  \sum\limits_{n\in I_i} \cos^2(n\theta_{k,\ell} +\gamma_{k,\ell})} {\widetilde{c}_{k',\ell'}^2 |\rho_{k',\ell'}|^{2A_{i+1}} \sum\limits_{n\in I_i}  \cos^2(n\theta_{k',\ell'} +\gamma_{k',\ell'})}  \geq 4r^2  \right\}
\end{equation}
Then, since we assumed $|\rho_{k,\ell}|, |\rho_{k',\ell'}|>1$, for $l \in A_i$ we have $|\rho_{k,\ell}|^{2A_i} \le |\rho_{k,\ell}|^{2l} \le |\rho_{k,\ell}|^{2A_{i+1}}$ (and likewise for $\rho_{k',\ell'}$), by the definitions of $I_i$ and $A_i$ in Definition \ref{def:intervals}. Therefore any $i\in D_{k,\ell;k',\ell'}$ must be in $C_{k,\ell;k',\ell'}^c$, and hence $|C_{k,\ell;k',\ell'}^c| \geq |D_{k,\ell;k',\ell'}|$, and $|D_{k,\ell;k',\ell'}|$ is easier to compute as $|\rho_{k,\ell}|$ is pulled out of the sum.

In the other cases, for example for $|\rho_{k,\ell}| \le 1, |\rho_{k',\ell'}| > 1$, we will have analogous expressions for $D_{k,\ell;k',\ell'}$ such as

\[D_{k,\ell;k',\ell'} =\left\{i: \frac{\widetilde{c}_{k,\ell}^2 |\rho_{k,\ell}|^{2A_{i}}  \sum\limits_{n\in I_i} \cos^2(n\theta_{k,\ell} +\gamma_{k,\ell})} {\widetilde{c}_{k',\ell'}^2 |\rho_{k',\ell'}|^{2A_i} \sum\limits_{n\in I_i}  \cos^2(n\theta_{k',\ell'} +\gamma_{k',\ell'})} \leq \frac{1} {4r^2}, \textrm{ or} \frac{\widetilde{c}_{k,\ell}^2 |\rho_{k,\ell}|^{2A_{i+1}}  \sum\limits_{n\in I_i} \cos^2(n\theta_{k,\ell} +\gamma_{k,\ell})} {\widetilde{c}_{k',\ell'}^2 |\rho_{k',\ell'}|^{2A_{i+1}} \sum\limits_{n\in I_i}  \cos^2(n\theta_{k',\ell'} +\gamma_{k',\ell'})}  \geq 4r^2  \right\}\]
As a result, the proof of Lemma \ref{lm:card_D} will differ slightly depending on the case we are considering, and these differences will be briefly mentioned in its proof. 

To simplify the notation, we denote the accumulated sums of the angles in the $i$'th interval $I_i$ by $\phi(i,\theta,\gamma)$, i.e.,
 \begin{equation}\label{equ:def_phi}
 \phi(i, \theta, \gamma): = \sum\limits_{n \in I_i} \cos^2(n\theta+\gamma).
 \end{equation}
and define a close approximation $\Gamma_{k,\ell; k',\ell'}^i$ to the complex ratio in the definition \eqref{equ:est_C} of $D_{k,\ell;k',\ell'}$
\begin{equation}\label{equ:def_gamma}
\Gamma_{k,\ell; k',\ell'}^i =  \frac{\widetilde{c}_{k,\ell}^2 |\rho_{k,\ell}|^{2A_{i}}  \sum\limits_{n\in I_i} \cos^2(n\theta_{k,\ell} +\gamma_{k,\ell})} {\widetilde{c}_{k',\ell'}^2 |\rho_{k',\ell'}|^{2A_i} \sum\limits_{n\in I_i}  \cos^2(n\theta_{k',\ell'} +\gamma_{k',\ell'})} = \frac{\widetilde{c}_{k,\ell}^2 |\rho_{k,\ell}|^{2A_{i}}   \phi(i,\theta_{k,\ell},\gamma_{k,\ell})} {\widetilde{c}_{k',\ell'}^2 |\rho_{k',\ell'}|^{2A_i} \phi(i,\theta_{k',\ell'},\gamma_{k',\ell'})}
\end{equation}
Note that in $\Gamma_{k,\ell; k',\ell'}^i$, $|\rho_{k,\ell}|$ and $|\rho_{k',\ell'}|$ have the same power while they do not in the definition \eqref{equ:est_C} of $D_{k,\ell;k',\ell'}$. We will need to take this difference into account when using $\Gamma_{k,\ell; k',\ell'}^i$ to bound $|D_{k,\ell;k',\ell'}|$. 

Using these definitions, we will now prove Lemma \ref{lem:coeff} for $\lambda_j^{1/2r} \le 1/4$. Explicitly, we shall prove the following main lemma. 

\begin{lem:coeffpart1}
There exists absolute uniform constant $C_2 \in \mathbb{R}^+$ such that for all $r \ge 2$ and $1/4 \geq \lambda_j^{1/2r} \geq \frac{C_2r^6}{N}$, 
\[|c_{k,\ell}| \leq \left(\frac{24}{ \min\{|\rho_{k,j}|^{2N},1\} N}\right)^{1/2},\] 
 for all $(k,\ell) \in [r-1] \times \{0,1\}$.
\end{lem:coeffpart1}

In order to prove this result, we will show the following two intermediate results. The first result is used to bound $\phi(i,\theta,\gamma)$. The second result is the key result that allows us to prove the main lemma; it allows us to ensure that $|C|$ is sufficiently large, i.e. that there is a sufficiently large set of intervals such that the energy of every component is sufficiently different from the energy of every other component. 

%\newtheorem*{lm:card_D}{\textbf{\emph{Lemma \ref{lm:card_D}}}}
%%\begin{lm}\label{lm:card_D}
%\begin{lm:card_D}
%For any $(k,j)\neq (k',j')$, if
% \[  \frac{\Gamma_{k,j;k',j'}^{i} }{\Gamma_{k,j;k',j'}^{i+1}}  \geq 4/3 \textrm{ for all $i$,}  \textrm{ or  } \frac{\Gamma_{k,j;k',j'}^{i} }{\Gamma_{k,j;k',j'}^{i+1}} \leq 3/4 \textrm{ for all } i  \]
% then $|D_{k,j;k',j'}| \geq N/\Delta N-16(1+\sqrt 2)\lambda^{1/2r} \Delta N -\log (4^9 r^{16}) $.
%%\end{lm}
%\end{lm:card_D}

%\newtheorem*{lm:small_lambda}{\textbf{\emph{Lemma \ref{lm:small_lambda}}}}
%%\begin{lm} \label{lm:small_lambda}
%\begin{lm:small_lambda}
%Assume $\lambda^{1/2r} \leq 1/4$, then either of the following two cases hold for any $k=0,..,r-1$ and $j=0,1$
%\begin{itemize}
%\item $\sin(\theta_{k,j} )= 0$;
%\item $\sin(\theta_{k,j}) \geq c \lambda^{1/2r}$.
%\end{itemize}
% Here $c=\sin(\pi /2r) /25$. 
%%\end{lm}
%\end{lm:small_lambda}

%\begin{lm}\label{lm:tau} 
\newtheorem*{lm:tau}{\textbf{\emph{Lemma \ref{lm:tau}}}}
\begin{lm:tau}
Suppose $\lambda_j^{1/2r} \leq 1/4$. For all $(k,\ell) \ne (k',\ell')$ with $(k,\ell), (k',\ell') \in \{(0,0)\} \cup (\{1,\cdots, \lfloor r/2 \rfloor\} \times \{0,1\})$, there exists absolute universal constant $C_4 \in \mathbb{R}^+$ such that if $\Delta N \geq \max \left\{ \frac{75}{\sin(\pi/2r) \lambda^{1/2r}},   \frac{r^2 + C_4\lambda^{1/2r}}{C_4\lambda^{1/2r}} \log 12\right\} =: E(r) \lambda^{-1/2r} $, either \[  \frac{\Gamma_{k,\ell;k',\ell'}^{i} }{\Gamma_{k,\ell;k',\ell'}^{i+1}}  \geq 4/3 \textrm{ for all $i \in [m-1]$,}  \textrm{ or  } \frac{\Gamma_{k,\ell;k',\ell'}^{i} }{\Gamma_{k,\ell;k',\ell'}^{i+1}} \leq 3/4 \textrm{ for all } i\in[m-1], \]
(for $m=\lfloor N/\Delta N \rfloor$ as defined in Definition \ref{def:intervals}).
Additionally,
\[
\phi(i,\theta_{k,\ell},\gamma_{k,\ell}) \geq \Delta N/3.
\]
for all $i \in [m], (k,\ell) \in \{(0,0)\} \cup (\{1,\cdots, \lfloor r/2 \rfloor\} \times \{0,1\})$.
\end{lm:tau}
%\end{lm}

%%\begin{lm}\label{lm:C_lbound} 
%\newtheorem*{lm:C_lbound}{\textbf{\emph{Lemma \ref{lm:C_lbound}}}}
%\begin{lm:C_lbound}
%Suppose $N\geq \nu_1(r) \lambda^{-1/2r} $, then $|C| \geq m/2 $, $m$ is the total number of intervals. Here $\nu_1(r)=4r(2r-1)(\log (4^9r^{16})(E(r)+1/2)+ 16(1+\sqrt 2)((E(r))(E(r)+1/2)+\frac{1}{2}(E(r)+1/2)))+2(E(r)+1/2)$ where $E(r)\equiv \max\{ \frac{75}{\sin(\pi/2r)},   \frac{r^2 + C_8\lambda^{1/2r}}{C_8}\log 12\}$, for $C_8$ the same absolute universal constant defined in Lemma \ref{lm:tau}. 
%\end{lm:C_lbound}

%\begin{corollary}\label{lm:C_lboundCor}
\newtheorem*{lm:C_lboundCor}{\textbf{\emph{Corollary \ref{lm:C_lboundCor}}}}
\begin{lm:C_lboundCor}
Suppose that the assumptions of both Lemmas \ref{lm:C_lbound} and \ref{lm:card_D} hold. Then, there exists an absolute universal constant $C_2$ such that if $N\geq C_2r^6\lambda_j^{-1/2r} $, then $|C| \geq m/2 $, where $m=\lfloor N/\Delta N \rfloor \ge 2$ is the total number of intervals.
\end{lm:C_lboundCor}
%\end{corollary}

Using these results, we have the following proof of Lemma \ref{lem:coeffpart1}.
 
\begin{proof}[Proof of Lemma \ref{lem:coeffpart1}]
For $(k,\ell) \in \{(0,0)\} \cup (\{1,\cdots, \lfloor r/2 \rfloor\} \times \{0,1\})$, we have $|c_{k,\ell}| \le |\tilde{c}_{k,\ell}|$ by its definition in Lemma \ref{lm:realrep}. Additionally, for $(k,\ell) \not \in \{(0,0)\} \cup (\{1,\cdots, \lfloor r/2 \rfloor\} \times \{0,1\})$, then we have $c_{0,1} = \overline{c_{0,0}}$, and $c_{k,\ell} = \overline{c_{r-k,\ell}}$, by Lemma \ref{lm:conjCoeffs}. Therefore it suffices to show  that $|\wtl c_{k,\ell}| \leq \left(\frac{48}{ \min\{|\rho_{k,\ell}|^{2N},1\} N}\right)^{1/2}$ for all $(k,\ell) \in \{(0,0)\} \cup (\{1,\cdots, \lfloor r/2 \rfloor\} \times \{0,1\})$. On each $I_t$ with $t\in C$, by the definition of $C$ (Definition \ref{def:appDC}), all of the components are sufficiently different from each other. Hence we can reorder them according to their energy on $I_t$. For simplicity of notation, we denote the component (corresponding to index $(k,\ell)$) which has the $q$th largest energy over $I_t$ as $\tilde{a}_{k,\ell} \cdot |\rho_{k,\ell}|^i \cos (i \theta_{k,\ell}+\gamma_{k,\ell}) =: \tilde{a}_{t}^{(q)}|\rho_{t}^{(q)}|^i\cos(i\theta_{t}^{(q)}+\gamma_{t}^{(q)}) =: P_t^{(q)}(i)$. We also denote $[\tilde{v}(\lambda_j)]_{i}$ as $\tilde{v}(i)$. Using this notation, \eqref{eq:compact_expre} becomes
\[\tilde{v}(i)=\sum^{2\lfloor r/2 \rfloor + 1} _{q=1} \tilde{a}_{t}^{(q)} \cdot |\rho_{t}^{(q)}|^i \cos (i \theta_{t}^{(q)}+\gamma_{t}^{(q)}) = \sum^{2\lfloor r/2 \rfloor + 1} _{q=1} P_t^{(q)}(i) ,\ \ \  i\in I_t. \]
%Suppose $t\in C$, then by the definition of $C$ in Definition \ref{def:appDC}, the energy of the $l$th component over $I_t$ is at least $(4r^2)$ times larger than the energy of the $(l+1)$th component over $I_t$. We can then use this property to see that
We can now see that
\begin{align}
\sum \limits_{i \in I_t} \tilde{v}^2(i) & =  \sum \limits_{i\in I_t} \left(\sum \limits_{q=1}^{2\lfloor r/2 \rfloor + 1} P_t^{(q)}(i) \right)^2  \notag \\
& = \sum \limits_{i\in I_t}  \left[  (P_t^{(1)}(i))^2+ 2 P_t^{(1)}(i) \sum\limits_{q=2}^{2\lfloor r/2 \rfloor + 1} P_t^{(q)}(i)  +  \left(\sum\limits_{q=2}^{2\lfloor r/2 \rfloor + 1} P_t^{(q)}(i)\right)^2\right] \notag \\
& \geq \sum \limits_{i\in I_t}  \left[\frac{1}{2}(P_t^{(1)}(i))^2  -  \left(\sum\limits_{q=2}^{2\lfloor r/2 \rfloor + 1} P_t^{(q)}(i)\right)^2 \right] \notag \\
& \geq \frac{1}{2}\sum \limits_{i\in I_t}   (P_t^{(1)}(i))^2  -  (2\lfloor r/2 \rfloor)   \sum\limits_{q=2}^ {2\lfloor r/2 \rfloor+1} \sum \limits_{i\in I_t} (P_t^{(q)}(i))^2 \notag 
\end{align}
where the first inequality is obtained by applying $2|ab| \leq \frac{1}{2} a^2 + 2b^2$ (derived from the Arithmetic Mean/Geometric Mean Inequality) to the cross term  followed by a direct simplification, and the second inequality uses the Cauchy-Schwarz inequality on the second term.

Now, suppose that $t \in C$. By the definition of $C$ in Definition \ref{def:appDC}, the energy of the $q$th component over $I_t$ is at least $(4r^2)$ times larger than the energy of the $(q+1)$th component over $I_t$, and hence $\sum_{i \in I_t} (P_t^{(q)}(i))^2 \le \frac{1}{4r^2} \sum_{i \in I_t} (P_t^{(1)}(i))^2$. Therefore, we have 

\[\frac{1}{2}\sum \limits_{i\in I_t}   (P_t^{(1)}(i))^2  -  (2\lfloor r/2 \rfloor)   \sum\limits_{q=2}^ {2\lfloor r/2 \rfloor+1} \sum \limits_{i\in I_t} (P_t^{(q)}(i))^2 \ge \frac{1}{4} \sum \limits_{i \in I_t} (P_t^{(1)}(i))^2.\]
Furthermore, since $P_t^{(q)}(i)$ is the component with $l$th largest energy over $I_t$, $\sum_{i \in I_t} (P_t^{(1)}(i))^2 \ge \sum_{i \in I_t} (P_t^{(q)}(i))^2$ for all $q \in \{1,...,2\lfloor r/2 \rfloor+1\}$. Hence, we have 

\begin{align}\label{eq:sum}
\sum \limits_{i \in I_t} \tilde{v}^2(i) 
& \geq \frac{1}{2}\sum \limits_{i\in I_t}   (P_t^{(1)}(i))^2  -  (2\lfloor r/2 \rfloor)   \sum\limits_{q=2}^ {2\lfloor r/2 \rfloor+1} \sum \limits_{i\in I_t} ((P_t^{(q)}(i))^2 \notag \\& \geq  \frac{1}{4} \sum \limits_{i \in I_t} (P_t^{(1)}(i))^2 \notag \\ &\geq  \frac{1}{4} \sum \limits_{i\in I_t}  (P_t^{(q)}(i))^2, \textrm{ for  any }  q \in \{1,...,2\lfloor r/2 \rfloor+1\} \notag\\
& =  \frac{1}{4} \sum \limits_{i\in I_t}  \wtl c_{k,\ell}^2 |\rho_{k,\ell}|^{2i}\cos^2(i\theta_{k,\ell}+\gamma_{k,\ell}), \textrm{ for  all } (k,\ell) \in \{(0,0)\} \cup (\{1,\cdots, \lfloor r/2 \rfloor\} \times \{0,1\}).
\end{align}

Now we can see that, since the above calculation holds for any $I_t$ with $t\in C$, we can sum over all such $I_t$ and use the fact that $|\mathbf{v}_j|=1$ and the above result \eqref{eq:sum} to get 
\[ 1 \geq \sum\limits_{t\in C} \sum\limits_{i\in I_t}  \tilde{v}^2(i) \geq \frac{1}{4}  \sum\limits_{t \in C} \wtl c_{k,\ell} ^2  \sum\limits_{i\in I_t}  (|\rho_{k,\ell}|^{2i} \cos^2(i\theta_{k,\ell} +\gamma_{k,\ell})  ).\]
Then, since $i \le N$, we must have $|\rho_{k,\ell}|^{2i} \geq \min\{|\rho_{k,\ell}|^{2N},1\}$, and by Lemma \ref{lm:tau}, we have $\phi(t,\theta_{k,\ell},\gamma_{k,\ell}) =  \sum\limits_{i \in I_t} \cos^2(i\theta_{k,\ell}+\gamma_{k,\ell}) \geq \Delta N/3$, hence
\[\frac{1}{4}  \sum\limits_{t \in C} \wtl c_{k,\ell} ^2  \sum\limits_{i\in I_t}  (|\rho_{k,\ell}|^{2i} \cos^2(i\theta_{k,\ell} +\gamma_{k,\ell})  ) \ge  \frac{1}{4} \frac{\Delta  N}{3} \sum\limits_{t \in C} \wtl c_{k,\ell}^2 \min\{|\rho_{k,\ell}|^{2N},1\}.\] 
Finally, Corollary \ref{lm:C_lboundCor} implies that $|C| \ge m/2$ for $m = \lfloor N/\Delta N\rfloor$, and it also implies that $m \ge 2$. Since $m \ge 2$, we also have $\frac{\lfloor N/\Delta N\rfloor}{2} \ge \frac{N/\Delta N}{4}$. Hence, 
\begin{align*}
\frac{1}{4} \frac{\Delta  N}{3} \sum\limits_{t \in C} \wtl c_{k,\ell}^2 \min\{|\rho_{k,\ell}|^{2N},1\} & \ge  \frac{1}{4} \frac{\Delta  N}{3} \frac{\lfloor N/\Delta N\rfloor}{2} \wtl c_{k,\ell}^2 \min\{|\rho_{k,\ell}|^{2N},1\} \\
& \ge \frac{1}{4} \frac{\Delta  N}{3} \frac{N/\Delta N}{4} \wtl c_{k,\ell}^2 \min\{|\rho_{k,\ell}|^{2N},1\}\\
& \ge \frac{N}{48} \wtl c_{k,\ell}^2 \min\{|\rho_{k,\ell}|^{2N},1\}.
\end{align*}
Combining these results, we therefore see that 
\begin{align*}
1 & \ge \frac{1}{4}  \sum\limits_{t \in C} \wtl c_{k,\ell} ^2  \sum\limits_{i\in I_t}  (|\rho_{k,\ell}|^{2i} \cos^2(i\theta_{k,\ell} +\gamma_{k,\ell})  ) \\
& \ge \frac{1}{4} \frac{\Delta  N}{3} \sum\limits_{t \in C} \wtl c_{k,\ell}^2 \min\{|\rho_{k,\ell}|^{2N},1\} \\
& \ge \frac{N}{48} \wtl c_{k,\ell}^2 \min\{|\rho_{k,\ell}|^{2N},1\}
\end{align*}
and rearranging this inequality produces the desired result.
\end{proof}

Therefore, to complete the proof of Lemma \ref{lem:coeffpart1} it suffices to prove our two intermediate results. We begin by proving the first result, Lemma \ref{lm:tau}. In order to prove this result, we first prove the following lemma, which will in turn allow us to bound $\sin(\theta_{k,\ell})$, which will allow us to bound $\phi(i,\theta_{k,\ell},\gamma_{k,\ell})$ using trigonometric identities. 

\begin{lm} \label{lm:small_lambda}
Assume $\lambda_j^{1/2r} \leq 1/4$, then one of the following two cases holds for any $k=0,..,r-1$ and $\ell=0,1$:
\begin{itemize}
\item if $r$ is even and $k=r/2$, $\sin(\theta_{k,\ell} )= 0$,
\item otherwise, $\sin(\theta_{k,\ell}) \geq \frac{\lambda_j^{1/2r}\sin(\pi /2r)}{25}$.
\end{itemize}
\end{lm}

\begin{proof}
 Recall that $\theta_{k,\ell}$ is defined in Lemma \ref{lm:realrep} such that $\rho_{k,\ell} = |\rho_{k,\ell}| \mathbbm{e}^{\theta_{k,\ell}\mathbbm{i}}$. Thus, if $\rho_{k,\ell}$ is a real root (which by Corollary \ref{cor:onlyRealRoots} happens iff $r$ is even and $k=r/2$) then $\sin(\theta_{k,\ell} )= 0$. Otherwise, as discussed in the proof of Lemma \ref{lem:BasicRootsEqns}, $\sqrt {\rho_{k,\ell}} = \frac{c_k\pm \sqrt{c_k^2+4}}{2}$ with $c_k$ complex, where we recall that 
 
 \[c_k=\pm i \lambda_j^{1/2r}e^{k\pi i/r}.\] 
Hence
\begin{equation} \label{eq:sqrtLbound} |\textrm{Im}(\sqrt {\rho_{k,\ell}})|= \frac{1}{2}\left|\frac{c_k - \overline{c_k} }{2} \mp \frac{\sqrt{c_k^2+4}-\sqrt{\overline{c_k^2+4}}}{2}\right|= \left|\frac{c_k-\overline{c_k}}{4}\left(1\mp \frac{c_k+\overline{c_k}}{\sqrt{c_k^2+4}+\sqrt{\overline{c_k}^2+4}}\right)\right|\geq \frac{|c_k-\overline{c_k}|}{8},
\end{equation} 
where the last inequality follows since $|c_k|=\lambda_j^{1/2r}\le1/4$, by our assumption on $\lambda_j$, and therefore $|c_k+\overline{c_k}| \le \frac{1}{2}$ and $|\sqrt{c_k^2+4}+\sqrt{\overline{c_k}^2+4}| = 2 |\textrm{Re}(\sqrt{c_k^2+4})| \ge \frac{3}{2} \ge 1$ .

Using a similar argument to bound $1\mp \frac{c_k+\overline{c_k}}{\sqrt{c_k^2+4}+\sqrt{\overline{c_k}^2+4}}$ from above, we see that 
\begin{equation} \label{eq:sqrtUbound}|\textrm{Im}(\sqrt {\rho_{k,\ell}})|= \left|\frac{c_k-\overline{c_k}}{4}\left(1\mp \frac{c_k+\overline{c_k}}{\sqrt{c_k^2+4}+\sqrt{\overline{c_k}^2+4}}\right)\right|\leq \frac{|c_k-\overline{c_k}|}{3}.
\end{equation}
This implies
 \[
\left| \sin\left(\frac{\theta_{k,\ell}}{2}\right) \right| = \frac{|\textrm{Im}(\sqrt{ \rho_{k,\ell}})|}{|\sqrt{\rho_{k,\ell}} |}  \geq \frac{|c_k-\overline{c_k}|}{24} = \frac{|\textrm{Im}(\pm i \lambda_j^{1/2r}e^{k\pi i/r})|}{12} =  \frac{|\textrm{Re}(\pm \lambda_j^{1/2r}e^{k\pi i/r})|}{12},
 \]
 where the second inequality follows from (\ref{eq:sqrtLbound}) plus the fact that $|\sqrt{\rho_{k,\ell}}| \le 1+\sqrt{2} \le 3$ from (\ref{eq:rhoBound2}). First, since we are assuming that $\rho_{k,\ell}$ is not real, by Corollary \ref{cor:onlyRealRoots} we must have $k \ne r/2$, and hence $e^{k\pi i/r}$ must have a nonzero real part. We then see that the smallest nonzero value of $|\textrm{Re}(\pm \lambda_j^{1/2r}e^{k\pi i/r})|$ is greater than or equal to the smallest nonzero value of $|\textrm{Re}(\pm \lambda_j^{1/2r}e^{k\pi i/2r})|$, which is clearly $\lambda_j^{1/2r}\cos(\frac{\pi}{2}-\frac{\pi}{2r}) = \lambda_j^{1/2r}\sin(\pi/2r)$ since $e^{k\pi i/2r}$ lies on the imaginary axis for $k=r$ (hence, setting $k=r-1$ gives the desired result). 
 
This bound therefore implies
 \[
\left| \sin\left(\frac{\theta_{k,\ell}}{2}\right) \right|\ge  \frac{|\textrm{Re}(\pm \lambda_j^{1/2r}e^{k\pi i/r})|}{12} \ge \lambda_j^{1/2r} \frac{\sin(\pi/2r)}{12}.
 \]
We also have 
  \[
\left| \sin\left(\frac{\theta_{k,\ell}}{2}\right) \right| = \frac{|\textrm{Im}(\sqrt{\rho_{k,\ell}}|)}{|\sqrt{\rho_{k,\ell}}| }  \leq |c_k-\overline{c_k}| \leq 1/2,
 \]
 where the first inequality follows from (\ref{eq:sqrtUbound}) and since $|\sqrt{\rho_{k,\ell}}| \ge (1+\sqrt{2})^{-1} \ge \frac{1}{3}$ from (\ref{eq:rhoBound2}). The second inequality follows since $|c_k| \le 1/4$.
Thus, 
\[ |\sin(\theta_{k,\ell})| =  2\left| \sin\left(\frac{\theta_{k,\ell}}{2}\right)\right| \left| \sqrt{ 1- \sin^2 \left(\frac{\theta_{k,\ell}}{2}\right)} \right| \geq 2\left(\lambda_j^{1/2r} \frac{\sin(\pi/2r)}{12}\right) \left(\sqrt{3/4}\right) \ge \lambda_j ^{1/2r} \frac{ \sin(\pi/2r)}{25}.\]
\end{proof}

Using this bound, we can now prove Lemma \ref{lm:tau}, one of the two lemmas needed to prove the main result. 

\begin{lm}\label{lm:tau} Suppose $\lambda_j^{1/2r} \leq 1/4$. There exists absolute universal constant $C_4 \in \mathbb{R}^+$ such that if $\Delta N \geq \max \left\{ \frac{75}{\sin(\pi/2r) \lambda_j^{1/2r}},   \frac{r^2 + C_4\lambda_j^{1/2r}}{C_4\lambda_j^{1/2r}} \log 12\right\} =: E(r) \lambda_j^{-1/2r} $, for all $(k,\ell) \ne (k',\ell')$ with $(k,\ell), (k',\ell') \in \{(0,0)\} \cup (\{1,\cdots, \lfloor r/2 \rfloor\} \times \{0,1\})$ either \[  \frac{\Gamma_{k,\ell;k',\ell'}^{i} }{\Gamma_{k,\ell;k',\ell'}^{i+1}}  \geq 4/3 \textrm{ for all $i \in [m-1]$,}  \textrm{ or  } \frac{\Gamma_{k,\ell;k',\ell'}^{i} }{\Gamma_{k,\ell;k',\ell'}^{i+1}} \leq 3/4 \textrm{ for all } i\in[m-1], \]
(for $m=\lfloor N/\Delta N \rfloor$ as defined in Definition \ref{def:intervals}).
Additionally,
\[
\phi(i,\theta_{k,\ell},\gamma_{k,\ell}) \geq \Delta N/3.
\]
for all $i \in [m], (k,\ell) \in \{(0,0)\} \cup (\{1,\cdots, \lfloor r/2 \rfloor\} \times \{0,1\})$.\end{lm}

\begin{proof}
Recall that in \eqref{equ:def_phi}, $\phi(i,\theta,\gamma)$ is defined as accumulated sums of the angles in the $i$'th interval $I_i$, i.e.,
 $$
 \phi(i, \theta, \gamma) = \sum\limits_{n \in I_i} \cos^2(n\theta+\gamma).
 $$
We first derive bounds on $\phi(i,\theta_{k,\ell},\gamma_{k,\ell})$. As $\lambda_j^{1/2r} \leq 1/4$ by assumption, from Lemma \ref{lm:small_lambda} we know that either $|\sin(\theta_{k,\ell})|=0$ (if $r$ is even with $k=r/2$) or $|\sin(\theta_{k,\ell})| \geq \frac{\sin(\pi/2r)}{25} \lambda_j^{1/2r}$. First, we consider the case $|\sin(\theta_{k,\ell})|=0$. Since $r$ is even, with $k=r/2$, by the proof of Lemma \ref{lm:realrep}, we also have $a_{k,\ell}$ real, and hence, by the definitions of $\theta_{k,\ell}$, $\gamma_{k,\ell}$ in Lemma \ref{lm:realrep}, $\theta_{k,\ell},\gamma_{k,\ell} \in \pi \mathbb{Z}$. Hence, $\cos^2(t\theta_{k,\ell}+\gamma_{k,\ell})=1$ for all $t \in \mathbb{Z}$. As a result,
\[\phi(i,\theta_{k,\ell},\gamma_{k,\ell}) = \sum \limits_{t=A_{i}}^{A_{i+1}-1} \cos^2(t\theta_{k,\ell}+\gamma_{k,\ell}) = \sum\limits_{t=A_{i+1}}^{A_{i+2}-1} 1 = \Delta N.\]

Next, we consider the case where $|\sin(\theta_{k,\ell})| \geq \frac{\lambda_j^{1/2r}\sin(\pi /2r)}{25}$. Since $\Delta N \ge \frac{75}{\sin(\pi/2r) \lambda_j^{1/2r}}$, we have that $\frac{1}{\Delta N} \le \frac{\sin(\pi/2r) \lambda_j^{1/2r}}{75} \le \frac{|\sin(\theta_{k,\ell})|}{3}$, and hence $|\sin(\theta_{k,\ell})| \geq \frac{3}{\Delta N}$. 
As a result,
\begin{align}
\phi(i,\theta_{k,\ell},\gamma_{k,\ell}) &= \sum \limits_{t=A_{i}}^{A_{i+1}-1} \cos^2(t\theta_{k,\ell}+\gamma_{k,\ell}) \notag \\
& = \frac{\Delta N}{2}+\sum\limits_{t=A_{i}}^{A_{i+1}-1} \frac{1}{2}\cos (2t\theta_{k,\ell}+2\gamma_{k,\ell}) \notag \\
& = \frac{\Delta N}{2}+ \frac{\sin ((2A_{i+1}-1)\theta_{k,\ell}+2\gamma_{k,\ell})-\sin((2A_{i}-1)\theta_{k,\ell} +2\gamma_{k,\ell})}{4\sin \theta_{k,\ell}} \notag \\
& \geq \frac{\Delta N}{2} - \frac{1}{ 2|\sin(\theta_{k,\ell})|}\notag \\
& \geq \frac{\Delta N}{2} - \frac{\Delta N}{6} \notag\\
& \geq \Delta N/3. \label{eq:DeltaN}
\end{align}
where the second-to-last inequality follows since $|\sin(\theta_{k,\ell})| \geq \frac{3}{\Delta N}$.

Thus, since it holds in all cases, $\phi(i,\theta_{k,\ell},\gamma_{k,\ell}) \ge \Delta N/3$ for all $i \in [m]$, $(k,\ell) \ne (k',\ell')$ with $(k,\ell),(k',\ell') \in \{(0,0)\} \cup (\{1,\cdots, \lfloor r/2 \rfloor\} \times \{0,1\})$. Also noting that $\phi(i,\theta_{k,\ell},\gamma_{k,\ell})\leq \Delta N$ holds for all $i \in [m]$, we have 
 \begin{equation}\label{eq:phi}
1/3\leq \frac{ \phi(i,\theta_{k,\ell},\gamma_{k,\ell})}{ \phi(i+1,\theta_{k,\ell},\gamma_{k,\ell})} \leq 3
  \end{equation}
for all $i \in [m-1]$.

Suppose that $|\rho_{k,\ell}| > |\rho_{k',\ell'}|$, with $(k,\ell),(k',\ell') \in \{(0,0)\} \cup (\{1,\cdots, \lfloor r/2 \rfloor\} \times \{0,1\})$. We first note that the our choices of $(k,\ell),(k',\ell')$ do not allow $\rho_{k,\ell} , \rho_{k',\ell'}$ to be conjugates, since they are paired in \eqref{eq:compact_expre}.  Lemma \ref{lm:lambda}, which states that there exists $c,C$ absolute positive constants such that for any two roots $\rho$, $\tilde{\rho}$, if $\rho$ and $\tilde{\rho}$ are not conjugates, inverses or conjugate inverses, then $c r^{-2}\lambda_j^{1/2r}\leq \left||\tilde{\rho}|-|{\rho}|\right|\leq C \lambda_j^{1/2r}$ and (\ref{eq:rhoBound2}), which states that any root $\rho$ has $(1+\sqrt{2})^{-2} \le |\rho| \le (1+\sqrt{2})^2$, imply that there exists an absolute constant $C_4$ such that $\frac{|\rho_{k,\ell}|-|\rho_{k',\ell'}|}{|\rho_{k',\ell'}|} \ge \frac{C_4\lambda_j^{1/2r}}{r^2}$. Note that we can apply Lemma \ref{lm:lambda} because the two roots cannot be conjugates, and if the roots are inverses or conjugate inverses, their magnitudes must be separated by $1$, hence $||\rho_{k,\ell}|-|\rho_{k,\ell}^{-1}||\le ||\rho_{k,\ell}|-1|$, $||\rho_{k,\ell}|-|\overline{\rho_{k,\ell}}^{-1}|| \le ||\rho_{k,\ell}|-1|$ and since there must be a root (namely $\rho_{0,0}$) with magnitude $1$, we can apply Lemma \ref{lm:lambda} to $||\rho_{k,\ell}|-|\rho_{0,0}||$. (Note that we cannot have $|\rho_{k,\ell}|=|\rho_{k,\ell}^{-1}|=1$ or $|\rho_{k,\ell}|=|\overline{\rho_{k,\ell}}^{-1}|=1$, since such roots would be conjugate or equal, respectively.) 

We then see that
\beq \label{eq:abs_lambda}
\left(\frac{|\rho_{k,\ell}|}{|\rho_{k',\ell'}|}\right)^{\Delta N} =\left(1+\frac{|\rho_{k,\ell}|-|\rho_{k',\ell'}|}{|\rho_{k',\ell'}|} \right)^{\Delta N} \geq \left(1+  \frac{C_4\lambda_j^{1/2r}}{r^2} \right)^{\Delta N} \geq e^{\frac{C_4\lambda_j^{1/2r}}{r^2+C_4 \lambda_j^{1/2r}}\cdot \Delta N }\geq 12
\eeq
where we also used the fact that $\Delta N \ge \frac{r^2 + C_4\lambda_j^{1/2r}}{C_4\lambda_j^{1/2r}} \log 12$ and the fact that $\log (1+x) \geq \frac{x}{1+x}$ for $x\geq -1$. If $|\rho_{k,\ell}| < |\rho_{k',\ell'}|$, we apply the same argument to $\left(\frac{|\rho_{k',\ell'}|}{|\rho_{k,\ell}|}\right)^{\Delta N}$ to find that $\left(\frac{|\rho_{k',\ell'}|}{|\rho_{k,\ell}|}\right)^{\Delta N} \ge 12$ and hence that $\left(\frac{|\rho_{k,\ell}|}{|\rho_{k',\ell'}|}\right)^{\Delta N} \le \frac{1}{12}$. 

To complete the proof, we again assume first that $|\rho_{k,\ell}| > |\rho_{k',\ell'}|$. Then, by the definition \eqref{equ:def_gamma} of $\Gamma^i_{k,\ell;k',\ell'}$ we have
\[
\frac{\Gamma_{k,\ell;k',\ell'}^{i} }{\Gamma_{k,\ell;k',\ell'}^{i+1}} = \left(\frac{|\rho_{k,\ell}|}{|\rho_{k',\ell'}|}\right)^{-\Delta N}\frac{ \phi(i,\theta_{k,\ell},\gamma_{k,\ell})}{ \phi(i+1,\theta_{k,\ell},\gamma_{k,\ell})} \frac{ \phi(i+1,\theta_{k',\ell'},\gamma_{k',\ell'})}{ \phi(i,\theta_{k',\ell'},\gamma_{k',\ell'})} \leq \left(\frac{1}{12}\right)(3)(3)= \frac{3}{4}
\]
for all $i \in [m-1]$.
In the case where $|\rho_{k,\ell}| < |\rho_{k',\ell'}|$, we will instead get 
\[
\frac{\Gamma_{k,\ell;k',\ell'}^{i} }{\Gamma_{k,\ell;k',\ell'}^{i+1}} \geq (12)\left(\frac{1}{3}\right)\left(\frac{1}{3}\right) = \frac{4}{3}
\]
for all $i \in [m-1]$.
\end{proof}

Next, we will prove the second main result, Lemma \ref{lm:C_lboundCor}. In order to prove this bound, we recall from (\ref{equ:est_C}) that we can bound $|D_{k,\ell;k',\ell'}|$ and use the fact that $|C^c_{k,\ell;k',\ell'}| \ge |D_{k,\ell;k',\ell'}|$ to bound $|C^c_{k,\ell;k',\ell'}|$ and ultimately bound $|C|$ in Lemma \ref{lm:C_lbound}. We therefore prove the following lemma which uses the result in Lemma \ref{lm:tau} to bound $|D_{k,\ell;k',\ell'}|$. 

%Proofread cases argument...
\begin{lm}\label{lm:card_D}
For any $(k,\ell)\neq (k',\ell')$ with $(k,\ell), (k',\ell') \in \{(0,0)\} \cup (\{1,\cdots, \lfloor r/2 \rfloor\} \times \{0,1\})$, if either
 \[  \frac{\Gamma_{k,\ell;k',\ell'}^{i} }{\Gamma_{k,\ell;k',\ell'}^{i+1}}  \geq 4/3 \textrm{ for all $i \in [m-1]$,}  \textrm{ or  } \frac{\Gamma_{k,\ell;k',\ell'}^{i} }{\Gamma_{k,\ell;k',\ell'}^{i+1}} \leq 3/4 \textrm{ for all } i \in [m-1] \]
 then $|D_{k,\ell;k',\ell'}| \geq \lfloor N/\Delta N \rfloor -16(1+\sqrt 2)\lambda_j^{1/2r} \Delta N -\log (4^{10} r^{16}) $.
\end{lm}
\begin{proof} As mentioned in the discussion following Definition \ref{def:appDC}, this proof will have slightly different arguments depending whether $|\rho_{k,\ell}|, |\rho_{k',\ell'}|$ are greater than or less than $1$.  

Without loss of generality, assume $\frac{\Gamma_{k,\ell;k',\ell'}^{i} }{\Gamma_{k,\ell;k',\ell'}^{i+1}}  \geq 4/3$ (As mentioned in a later note, the same argument holds if $\frac{\Gamma_{k,\ell;k',\ell'}^{i} }{\Gamma_{k,\ell;k',\ell'}^{i+1}}  \leq 3/4$).
Let $S= \frac{2(1+\sqrt 2)\lambda_j^{1/2r} \Delta N + \log 4 r^2}{\log 2} $, then $2^{S}= 4r^2 e^{2(1+\sqrt 2)\lambda_j^{1/2r} \Delta N}$.
By Lemma \ref{lm:interval} with $B_i = 1$, $B'_i=\Gamma_{k,\ell;k',\ell'}^i$ (for all $i$) and hence $\rho_1=1$, $\rho_2=3/4$, and $\alpha =4/3$ (note that the hypotheses of the lemma are fulfilled since $\frac{\Gamma_{k,\ell;k',\ell'}^{i} }{\Gamma_{k,\ell;k',\ell'}^{i+1}}  \geq 4/3$ for all $i$), setting $S=\frac{q}{2}\log_2(4/3)$ and solving for $q$, we see that except for $\lceil 2S/\log_2 (4/3)\rceil+1 \ge 2S/\log_2 (4/3) + 2$ consecutive subsets, it holds 
\begin{equation}\label{eq:interval}
\Gamma_{k,\ell;k',\ell'}^{i} \le 2^{-S } \textrm{ or }  \Gamma_{k,\ell;k',\ell'}^{i} \ge 2^{S }.
\end{equation}
Note that Lemma \ref{lm:interval} can also be used (with $B_i = \Gamma_{k,\ell;k',\ell'}^i$, $B'_i=1$, $\rho_1=4/3$, $\rho_2=1$, and $\alpha = 4/3$) to derive the same result in the case of $\frac{\Gamma_{k,\ell;k',\ell'}^{i} }{\Gamma_{k,\ell;k',\ell'}^{i+1}} \leq 3/4$.

If $|\rho_{k,\ell}| \le 1, |\rho_{k',\ell'}| > 1$ or $|\rho_{k,\ell}| > 1, |\rho_{k',\ell'}| \le 1$, then the expressions in $D_{k,\ell;k',\ell'}$ will be either $\Gamma^i_{k,\ell;k',\ell'}$ or $\Gamma^{i+1}_{k,\ell;k',\ell'}$. As a result, since $2^S = 4r^2 e^{2(1+\sqrt 2)\lambda_j^{1/2r} \Delta N} \ge 4r^2$, we have the desired conclusion in all but $2S/\log_2 (4/3)+2$ consecutive subsets.
%paragraph break because different cases?

Next we note that in the case that $|\rho_{k,\ell}|, |\rho_{k',\ell'}| > 1$, we then note that 
\begin{equation}\label{eq:exprhobound}
|\rho_{k,\ell}|^{2\Delta N} = (1+ |\rho_{k,\ell}|-1)^{2\Delta N} \leq (1+ (1+\sqrt 2)\lambda_j^{1/2r})^{2\Delta N} \leq e^{2(1+\sqrt 2) \lambda_j^{1/2r}\Delta N }, 
\end{equation}
where the first inequality follows by Lemma \ref{lem:BasicRootsEqns} and the second inequality follows since $1+x \le e^x$ for all $x$. 
Now suppose that $\Gamma_{k,\ell;k',\ell'}^{i} \le 2^{-S}$. If this is the case, then we have
\begin{align}\label{eq:larger}
\wtl c_{k',\ell'}^2 |\rho_{k',\ell'}|^{2A_i}  \phi (i, \theta_{k',\ell'}, \gamma_{k',\ell'}) & \ge 2 ^{S} \wtl c_{k,\ell}^2 |\rho_{k,\ell}|^{2A_i} \phi (i, \theta_{k,\ell}, \gamma_{k,\ell}) \notag\notag\\
&=2 ^{S} \wtl c_{k,\ell}^2 |\rho_{k,\ell}|^{2A_{i+1}} |\rho_{k,\ell}|^{-2\Delta N}  \phi (i, \theta_{k,\ell}, \gamma_{k,\ell}) \notag \\ &\geq 2^{S}e^{-2(1+\sqrt 2) \lambda_j^{1/2r}\Delta N } \wtl c_{k,\ell}^2 |\rho_{k,\ell}|^{2A_{i+1}} \phi (i, \theta_{k,\ell}, \gamma_{k,\ell}) \notag \\
& = 4r^2 \wtl c_{k,\ell}^2 |\rho_{k,\ell}|^{2A_{i+1}} \phi (i, \theta_{k,\ell}, \gamma_{k,\ell}).
\end{align}
where the first inequality is from the definition (\ref{equ:def_gamma}) of $\Gamma_{k,\ell;k',\ell'}^{i}$, and where the last equality is due to the definition of $S$. A similar argument using \eqref{eq:exprhobound} on $|\rho_{k',\ell'}|$ shows that 
\[ \wtl c_{k,\ell}^2 |\rho_{k,\ell}|^{2A_{i}} \phi (i, \theta_{k,\ell}, \gamma_{k,\ell}) \ge 4r^2 \wtl c_{k',\ell'}^2 |\rho_{k',\ell'}|^{2A_{i+1}}  \phi (i, \theta_{k',\ell'}, \gamma_{k',\ell'})\]
if $ \Gamma_{k,\ell;k',\ell'}^{i} \ge 2^{S }$. 

Also, an analogous argument implies the same result if $|\rho_{k,\ell}|, |\rho_{k',\ell'}| \le 1$, since using Lemma \ref{lem:charPolyProps1} to note that $\rho_{k,0}=\rho_{k,1}^{-1}$, we can apply the previous result \eqref{eq:exprhobound} to $\rho_{k,\ell'}=\rho_{k,\ell}^{-1}$ and therefore we see that
\[
|\rho_{k,\ell}|^{-2\Delta N}  \leq e^{2(1+\sqrt 2) \lambda_j^{1/2r}\Delta N }.
\]
 
Finally, the bound on the number of the excluded intervals 
\begin{align*}
2S/\log_2 (4/3)+2 &= \frac{4(1+\sqrt 2)\lambda_j^{1/2r} \Delta N + 2\log 4 r^2}{\log 2\log_2(4/3)}+2 \\
&\leq \frac{4(1+\sqrt 2)\lambda_j^{1/2r} \Delta N + 2\log 4 r^2}{\log(4/3)}+2\log(4) \\
&\leq 16(1+\sqrt 2)\lambda_j^{1/2r} \Delta N + 8\log 4 r^2 + \log(4^2)\\
&\leq 16(1+\sqrt 2)\lambda_j^{1/2r} \Delta N +\log (4^{10} r^{16}) 
\end{align*}
now implies the conclusion of this lemma.
\end{proof}

Given this bound on $|D_{k,\ell;k',\ell'}|$, and hence on $|C^c_{k,\ell;k',\ell'}|$, we can now prove a lemma which produces the desired bound $|C|$ for sufficiently large $N$. We note that we need the conditions of Lemma \ref{lm:tau} to be satisfied (i.e., we need to choose $\Delta N$ to be large enough), since the conclusions of Lemma \ref{lm:tau} are the assumptions of Lemma \ref{lm:card_D}.

\begin{lm}\label{lm:C_lbound} Suppose $\Delta N = \lceil E(r) \lambda_j^{-1/2r} \rceil+1$, where $E(r)$ is defined as in Lemma \ref{lm:tau}. Then for $N\geq \max(\nu_1(r) \lambda_j^{-1/2r},2\Delta N)$, $|C| \geq m/2 $, where $m=\lfloor N/\Delta N\rfloor \ge 2$ is the total number of intervals, and $\nu_1(r):=4r(2r-1)(\log (4^{10}r^{16})(E(r)+1/2)+ 16(1+\sqrt 2)((E(r))(E(r)+1/2)+\frac{1}{2}(E(r)+1/2)))+2(E(r)+1/2)$. \end{lm}
\begin{proof}
First, we clearly note that if $N \ge 2 \Delta N$, then $m =\lfloor N/\Delta N\rfloor \ge \lfloor 2 \Delta N/\Delta N\rfloor \ge 2$. Hence, it suffices to show that if $N \ge \nu_1(r) \lambda_j^{-1/2r}$, then $|C| \ge m/2$.

Set $\Delta N = \lceil E(r) \lambda_j^{-1/2r} \rceil+1$, where $E(r) = \max\{ \frac{75}{\sin(\pi/2r)},  \frac{r^2 + C_4\lambda_j^{1/2r}}{C_4}\log 12\}$ as per Lemma \ref{lm:tau}. We then note that 
\[(E(r)+1/2)\lambda_j^{-1/2r}=E(r)\lambda_j^{-1/2r}+(1/2)\lambda_j^{-1/2r} \ge (\Delta N-2) + 2 = \Delta N,\] 
since $\lambda_j^{1/2r}\le 1/4$ by assumption.
By Lemmas \ref{lm:card_D} and \ref{lm:tau}, we now have that
$$|D_{k,\ell;k',\ell'}| \geq \lfloor N/\Delta N \rfloor- 16(1+\sqrt 2)\lambda_j^{1/2r} \Delta N -\log (4^{10}r^{16}) . $$

Recall from \eqref{eq:def_C} in Definition \ref{def:appDC} that $C = \bigcap_{(k,\ell) \ne (k',\ell')} C_{k,\ell; k',\ell'}^c$ (where $(k,\ell), (k',\ell')$ run over $\{(0,0)\} \cup (\{1,\cdots, \lfloor r/2 \rfloor\} \times \{0,1\})$), and that $|C_{k,\ell;k',\ell'}^c| \geq |D_{k,\ell;k',\ell'}|$. Therefore, using the above bound on $|D_{k,\ell;k',\ell'}|$, the desired result follows since
\begin{align*}
|C| & =|\cap C_{k,\ell; k',\ell'}^c| \geq \lfloor N/\Delta N \rfloor -|\cup C_{k,\ell; k',\ell'}|  \\ & \geq \lfloor N/\Delta N \rfloor - \sum_{k,\ell,k',\ell', (k,\ell)\neq (k',\ell')}|C_{k,\ell;k',\ell'}|  \\ &\geq \lfloor N/\Delta N \rfloor - \sum_{k,\ell,k',\ell', (k,\ell)\neq (k',\ell')}(\lfloor N/\Delta N \rfloor- |D_{k,\ell;k',\ell'}|) \\ & \geq \lfloor N/\Delta N \rfloor - 2r(2r-1)(16(1+\sqrt 2)\lambda_j^{1/2r} \Delta N +\log (4^{10} r^{16}) ) \\& \ge \frac{\lfloor N/\Delta N \rfloor}{2} = m/2.
\end{align*}

In particular, the last inequality holds by substituting $\nu_1(r)=4r(2r-1)(\log (4^{10}r^{16})(E(r)+1/2)+ 16(1+\sqrt 2)((E(r))(E(r)+1/2)+\frac{1}{2}(E(r)+1/2)))+2(E(r)+1/2)$ and noting that 
\begin{align*}
\lfloor N/\Delta N \rfloor &\ge \frac{1}{\Delta N}[\nu_1(r) \lambda_j^{-1/2r}] - 1\\ 
&=\frac{1}{\Delta N}\Bigg[4r(2r-1)\left(\log (4^{10}r^{16})\lambda_j^{-1/2r}(E(r)+1/2)+ 16(1+\sqrt 2) \lambda_j^{-1/2r}((E(r))(E(r)+1/2)+\frac{1}{2}(E(r)+1/2))\right)\\
&\hspace{0.5in} +2\lambda_j^{-1/2r}(E(r)+1/2)\Bigg] - 1\\
&\ge \frac{1}{\Delta N}\left[4r(2r-1)\left(\log (4^{10}r^{16})\Delta N+ 16(1+\sqrt 2) \left(E(r)\Delta N+\frac{1}{2}\Delta N\right) \right)+2\Delta N\right] - 1\\
&= 4r(2r-1)\left(\log (4^{10}r^{16})+ 16(1+\sqrt 2) \left(E(r)+\frac{1}{2}\right)\right)+1\\
&\ge 4r(2r-1)\left(\log (4^{10}r^{16})+ 16(1+\sqrt 2)(\lambda_j^{1/2r}\Delta N)\right)+1
\end{align*} 
where the second and last inequalities follow from the result shown earlier: $(E(r)+1/2)\lambda_j^{-1/2r} \ge \Delta N$ (and hence $\lambda_j^{1/2r}\Delta N \le E(r)+\frac{1}{2}$).
\end{proof}

Finally, using this lemma, we can now show the desired result, Corollary \ref{lm:C_lboundCor}, by simplifying the bound for $\Delta N$ in Lemma \ref{lm:C_lbound}.

\begin{corollary}\label{lm:C_lboundCor} Suppose that the assumptions of both Lemmas \ref{lm:C_lbound} and \ref{lm:card_D} hold. Then, there exists an absolute universal constant $C_2$ such that if $N\geq C_2r^6\lambda_j^{-1/2r} $, then $|C| \geq m/2 $, where $m=\lfloor N/\Delta N \rfloor \ge 2$ is the total number of intervals.
\end{corollary}

\begin{proof}
By Lemma \ref{lm:C_lbound}, it suffices to show that there exists some $C_2$ such that $C_2r^6 \ge \nu_1(r)=4r(2r-1)\left[\log (4^{10}r^{16})(E(r)+1/2)+ 16(1+\sqrt 2)((E(r))(E(r)+1/2)+\frac{1}{2}(E(r)+1/2))\right]+2(E(r)+1/2)$, and $C_2r^6\lambda_j^{-1/2r} \ge 2(\lceil E(r) \lambda_j^{-1/2r} \rceil+1)$ for all positive integers $r$. We then see that there exists some absolute universal constant $C_5$ such that $C_5r^2 \ge E(r) = \max\{ \frac{75}{\sin(\pi/2r)},  \frac{r^2 + C_4\lambda_j^{1/2r}}{C_4}\log 12\}$ for all positive integers $r$, since $x - \sin(x) \le \frac{1}{2} x^2 \le \frac{\pi}{4}x$ (by Taylor's Remainder Theorem) for any $0 \le x \le \pi/2$, so $\sin(\pi/2r) \ge \frac{\pi(1-\pi/4)}{2r}$ for any $r \ge 1$ (and hence there exists constant $C_6$ such that $C_6 r \ge \frac{75}{\sin(\pi/2r)}$ for all positive integers $r$), and since $\lambda_j^{1/2r} \le \frac{1}{4}$ (so there exists constant $C_7$ such that $C_7r^2 \ge\frac{r^2 + C_4\lambda_j^{1/2r}}{C_4}\log 12$ for all positive integers $r$). Also, by a similar argument we see that there exists constant $C_8$ such that $C_8r^2 \ge 2E(r)+1 \ge 2E(r)+4\lambda_j^{1/2r} \ge 2(E(r) \lambda_j^{-1/2r}+2)\lambda_j^{1/2r} \ge 2(\lceil E(r) \lambda_j^{-1/2r} \rceil+1)\lambda_j^{1/2r}$. Hence, substituting into the expression for $\nu_1(r)$, and combining the two bounds, the desired $C_2$ must indeed exist.

\end{proof}

\subsection{The $\lambda_j^{1/2r} > 1/4$ case: Proof of Lemma \ref{lem:coeffpart2}}
In a similar fashion to the previous section, we will seek to prove Lemma \ref{lem:coeff} in the case that $\lambda_j^{1/2r}>1/4$. In particular, our main result is as follows:

\begin{lem:coeffpart2} There exists absolute uniform constant $C_3 \in \mathbb{R}^+$, such that for all $r \ge 2$, $N \geq C_3^{r}$, and $\lambda_j^{1/2r}>1/4$, 
\[
| c_{k,\ell}| \leq \left(\frac{48}{ \min\{|\rho_{k,j}|^{2N},1\} N}\right)^{1/2},
\]
 for all $(k,\ell) \in [r-1] \times \{0,1\}$.
\end{lem:coeffpart2}

To prove this result, we will use similar reasoning to the argument in the previous section. As before, we will prove two intermediate results, which we can use to prove the main result. In particular, we use the following two results, one which bounds $|\sin(\theta_{k,\ell})|$ and one which bounds $|C|$ (as defined in (\ref{eq:def_C}) in Definition \ref{def:appDC}):

\newtheorem*{lm:angle_bound}{\textbf{\emph{Lemma \ref{lm:angle_bound}}}}
\begin{lm:angle_bound}
Assume $\lambda_j^{1/2r}>1/4$. There exists absolute constant $C_9>0$ such that one of the following two cases holds for $k=1,...,r-1$, $\ell=0,1$: 
\begin{itemize}
\item if $r$ is even and $k=r/2$, $\sin(\theta_{k,\ell} )= 0$;
\item otherwise, $\sin(\theta_{k,\ell}) \geq  C_9^{-r} \lambda_j^{1/2r}$. 
\end{itemize}
Moreover, $|\sin(\theta_{0,0}) |=|\sin(\theta_{0,1}) | \ge \min\left\{\left|\sin\left(\frac{2\pi}{2N+1}\right)\right|, \frac{3\sqrt{7}}{32}\right\}$.\end{lm:angle_bound}
\newtheorem*{lm:C_lbound_large}{\textbf{\emph{Lemma \ref{lm:C_lbound_large}}}}
\begin{lm:C_lbound_large}
Let $\Delta N = \lceil C_{10}^2r^4 \rceil$  where $C_{10}>0$ is an absolute constant that is the same as in Lemma \ref{lm:tau_large} below, then $$|C| \geq  \lfloor N/\Delta N \rfloor - 2r(2r-1)(32 \Delta N + 8\log 16 r^2 + 2)  . $$ 
\end{lm:C_lbound_large} 

Note the similarity between Lemma \ref{lm:angle_bound} and Lemma \ref{lm:small_lambda} from the previous section, as well as between Lemma \ref{lm:C_lbound_large} and Corollary \ref{lm:C_lboundCor}. In particular, the reasoning used in this section is very similar to the reasoning used in the previous case. The main difference is that when $\lambda_j^{1/2r} > 1/4$, Lemma \ref{lm:small_lambda} no longer holds, so we have to use a weaker result, namely Lemma \ref{lm:angle_bound}. This weakness is later compensated for with a larger value of $|C|$, and ultimately, by a larger bound on $N$. 

Using these two main results, the proof of Lemma \ref{lem:coeffpart2} is fairly straightforward.

\begin{proof}[Proof of Lemma \ref{lem:coeffpart2}]
First, we note that by the same reasoning as in Lemma \ref{lem:coeffpart1}, it suffices to show that $| \tilde{c}_{k,\ell}| \leq \left(\frac{48}{ \min\{|\rho_{k,\ell}|^{2N},1\} N}\right)^{1/2}$ for all $(k,\ell) \in \{(0,0)\} \cup (\{1,\cdots, \lfloor r/2 \rfloor\} \times \{0,1\})$. 

We take $\Delta N = \lceil C_{10}^2r^4 \rceil$, for $C_{10}$ as in Lemma \ref{lm:tau_large} below. Therefore there exists an absolute universal constant $C_3$ such that 
\begin{equation}\label{eq:constassump}
C_3^r \ge \max\left\{6C_9^r,12(\Delta N + 2r(2r-1)(\Delta N)(32 \Delta N + 8\log 16 r^2 + 2)),\frac{16}{\sqrt{7}}\right\},
\end{equation}
for all $r>1$, and let $N \ge C_3^r$.

The proof of this theorem follows the same line of reasoning as the proof of Lemma \ref{lem:coeffpart1}, until the end of \eqref{eq:sum}, which reads 
\[
\sum\limits_{i\in I_t } \tilde{v}^2(i)\geq \frac{1}{4} \sum\limits_{i\in I_t} \wtl c_{k,\ell}^2 |\rho_{k,\ell}|^{2i} \cos^2(i \theta_{k,\ell}+\gamma_{k,\ell}), \textrm{ for  all } (k,\ell) \in \{(0,0)\} \cup (\{1,\cdots, \lfloor r/2 \rfloor\} \times \{0,1\})
\]
and for all $t\in C$ where $$|C| \geq  \lfloor N/\Delta N \rfloor - 2r(2r-1)(32 \Delta N + 8\log 16 r^2 + 2), $$  since  the assumption of Lemma \ref{lm:C_lbound_large} is fulfilled for the given choice of $\Delta N$. Therefore, analogously to Lemma \ref{lem:coeffpart1}, we see that
\begin{align}\label{eq:bound1} 1 &\geq \sum\limits_{t\in C} \sum\limits_{i\in I_t}  \tilde{v}^2(i) \notag \\ &\geq \frac{1}{4}  \sum\limits_{t \in C} \wtl c_{k,\ell} ^2  \sum\limits_{i\in I_t}  (|\rho_{k,\ell}|^{2i} \cos^2(i\theta_{k,\ell} +\gamma_{k,\ell})  ) \notag \\ & \geq \frac{1}{4} \sum\limits_{t \in C} \wtl c_{k,\ell}^2 \min\{|\rho_{k,\ell}|^{2N},1\}  \sum\limits_{i\in I_t}  \cos^2(i\theta_{k,\ell}  +\gamma_{k,\ell} ) \notag \\  &\geq \frac{1}{4}\wtl c_{k,\ell}^2 \min\{|\rho_{k,\ell}|^{2N},1\} \left( \sum\limits_{i \in [N]} \cos^2(i\theta_{k,\ell} +\gamma_{k,\ell}) - |\{i: i\in I_t, t\in C^c\}|\right).
\end{align}
where the last inequality holds since $\cos^2(i\theta_{k,\ell}  +\gamma_{k,\ell} ) \le 1$ and by rearranging the sums.

If $|\sin(\theta_{k,\ell})|=0$, $\rho_{k,\ell}$ is therefore real, and hence as discussed in Lemma \ref{lm:tau}, we thus have $\cos^2(t\theta_{k,\ell}+\gamma_{k,\ell})=1$ for all $t$, and hence $ \sum\limits_{i \in [N]} \cos^2(i\theta_{k,\ell} +\gamma_{k,\ell}) = N \ge \frac{N}{6}$, and if $|\sin(\theta_{k,\ell})| \ne 0$, we then have 

\begin{align}\label{eq:coeff_one}
 \sum\limits_{i \in [N]} \cos^2(i\theta_{k,\ell} +\gamma_{k,\ell}) &= \frac{N}{2} + \frac{1}{2} \sum_{i \in [N]} \cos(2i\theta_{k,\ell} +2\gamma_{k,\ell}) \notag\\
 &= \frac{ N}{2}+\frac{\sin(N \theta_{k,\ell} )\cos (( N+1)\theta_{k,\ell}+2\gamma_{k,\ell})}{2\sin \theta_{k,\ell}} \notag\\ 
& \geq \frac{N}{2} - \frac{1}{2|\sin \theta_{k,\ell}|} \notag\\
& \ge \frac{N}{2} - \max\left\{2C_9^{r},\frac{1}{2}\left|\sin\left(\frac{2\pi}{2N+1}\right)\right|^{-1},\frac{16}{3\sqrt{7}}\right\} \notag\\
&\geq \frac{N}{2} -  \max\left\{2C_9^{r},\frac{5N}{24},\frac{16}{3\sqrt{7}}\right\}\notag\\
&\geq \frac{N}{6},
\end{align}
where the third to last inequality in \eqref{eq:coeff_one} used Lemma \ref{lm:angle_bound}, namely that $|\sin(\theta_{k,\ell})| \ge \min\left\{\frac{1}{4}C_9^{-r},\left|\sin\left(\frac{2\pi}{2N+1}\right)\right|, \frac{3\sqrt{7}}{32}\right\}$, while the second to last inequality used the fact that $|\sin(x)| \ge \frac{3}{\pi}{|x|}$ for $|x| \le \frac{\pi}{6}$ and the fact that $N \ge 6$ implies that $\frac{2\pi}{2N+1} \le \frac{\pi}{6}$, so therefore $\left|\sin\left(\frac{2\pi}{2N+1}\right)\right|^{-1} \le \frac{\pi}{3}\frac{2N+1}{2\pi} = \frac{2N+1}{6}\le\frac{5N}{12}$ (where the last inequality holds since $N \ge 2$) and the last inequality used the assumption \eqref{eq:constassump} on $N$. 
Hence, we see that $ \sum\limits_{i \in [N]} \cos^2(i\theta_{k,\ell} +\gamma_{k,\ell}) \ge \frac{N}{6}$ in all cases. 

In addition, we have 
\begin{align*}
|\{i: i\in I_t, t\in C^c\}| & \leq  N - \Delta N |C| \\ &\leq (N-\lfloor N/\Delta N \rfloor \Delta N) + 2r(2r-1)(32 \Delta N + 8\log 16 r^2 + 2)  \\ &\leq \Delta N + 2r(2r-1)(\Delta N)(32 \Delta N + 8\log 16 r^2 + 2),
\end{align*}
since, as discussed earlier, we have
$$|C| \geq  \lfloor N/\Delta N \rfloor - 2r(2r-1)(32 \Delta N + 8\log 16 r^2 + 2). $$

Thus, combining these two results, and the earlier result \eqref{eq:bound1}, we have 
\begin{align}\label{eq:coeff_unit}
1 &\ge \frac{1}{4}\wtl c_{k,\ell}^2 \min\{|\rho_{k,\ell}|^{2N},1\} ( \sum\limits_{i \in [N]} \cos^2(i\theta_{k,\ell} +\gamma_{k,\ell}) - |\{i: i\in I_t, t\in C^c\}|)\notag \\ &\geq  \frac{1}{4}\wtl c_{k,\ell}^2 \min\{|\rho_{k,\ell}|^{2N},1\}\left(\frac{N}{6} -( \Delta N + 2r(2r-1)(\Delta N)(32 \Delta N + 8\log 16 r^2 + 2))\right)  \notag \\
& \geq \frac{N}{48}\wtl c_{k,\ell}^2 \min\{|\rho_{k,\ell}|^{2N},1\},
\end{align}
where the last inequality follows from the assumption \eqref{eq:constassump} on $N$.
Rearranging \eqref{eq:coeff_unit} thus completes the proof.
\end{proof}

Thus, to complete this proof, it suffices to prove the intermediate results Lemmas \ref{lm:angle_bound} and \ref{lm:C_lbound_large}. We first prove Lemma \ref{lm:angle_bound}, which we use to bound $|\sin(\theta_{k,\ell})|$.

\begin{lm} \label{lm:angle_bound}
Assume $\lambda_j^{1/2r}>1/4$. There exists absolute constant $c>0$ such one of the following two cases holds for $k=1,...,r-1$, $\ell=0,1$: 
\begin{itemize}
\item if $r$ is even and $k=r/2$, $\sin(\theta_{k,\ell} )= 0$;
\item otherwise, $\sin(\theta_{k,\ell}) \geq  C_9^{-r} \lambda_j^{1/2r}$.\end{itemize}
Moreover, $|\sin(\theta_{0,0}) |=|\sin(\theta_{0,1}) | \ge \min\left\{\left|\sin\left(\frac{2\pi}{2N+1}\right)\right|, \frac{3\sqrt{7}}{32}\right\}$\end{lm}
\begin{proof}
From Lemma \ref{lm:sigD}, we know that $\lambda_j^{1/2r} \leq 2\cos\left(\frac{\pi}{2N+1}\right)$. Also, by our assumption, we have $\lambda_j^{1/2r} > 1/4$. Hence, we have $\frac{1}{4^{2r}} < \lambda_j \le 4^r\cos^{2r}\left(\frac{\pi}{2N+1}\right)$. From Lemmas  \ref{lem:charPolyProps1} and \ref{lem:BasicRootsEqns}, we know that $\rho_{0,0} = \frac{2-\lambda_j^{1/r} +\sqrt{\lambda_j^{2/r}-4\lambda_j^{1/r}}}{2}$ and  $\rho_{0,1} = \frac{2-\lambda_j^{1/r} -\sqrt{\lambda_j^{2/r}-4\lambda_j^{1/r}}}{2}$, $\rho_{0,0} = \overline{\rho_{0,1}}$, and $|\rho_{0,0}|=|\rho_{0,1}|=1$. Hence 
\[
|\sin(\theta_{0,0})| = |\sin(\theta_{0,1})| = \left|\frac{\rho_{0,0}-\rho_{0,1}}{2}\right| = \frac{\sqrt{4\lambda_j^{1/r}-\lambda_j^{2/r}}}{2},
\]
where the expression inside the square root was changed to $4\lambda_j^{1/r}-\lambda_j^{2/r}$ since $\lambda_j \in (0,4^r)$ implies $\lambda_j^{2/r}-4\lambda_j^{1/r}<0$.
Now, note that $\sqrt{4\lambda_j^{1/r}-\lambda_j^{2/r}}$ has one critical point in $(0,4^r)$ at $\lambda_j=2^r$ (a maximum), so $\sqrt{4\lambda_j^{1/r}-\lambda_j^{2/r}}$ takes on its minimum value in $\left[\frac{1}{4^{2r}}, 4^r\cos^{2r}\left(\frac{\pi}{2N+1}\right) \right] \subset (0,4^r)$ at one of the endpoints. 
Hence, since 
\begin{align*}
\sqrt{4\lambda_j^{1/r}-\lambda_j^{2/r}}\Big|_{\lambda_j=4^r\cos^{2r}\left(\frac{\pi}{2N+1}\right)} &= \sqrt{16\cos^2\left(\frac{\pi}{2N+1}\right)-16\cos^4\left(\frac{\pi}{2N+1}\right)}\\ &=4\sqrt{\cos^2\left(\frac{\pi}{2N+1}\right)\sin^2\left(\frac{\pi}{2N+1}\right)} = 4\left|\cos\left(\frac{\pi}{2N+1}\right)\sin\left(\frac{\pi}{2N+1}\right)\right|\\ &= 2\left|\sin\left(\frac{2\pi}{2N+1}\right)\right|,
\end{align*}
we thus have 
\[|\sin(\theta_{0,0})| = |\sin(\theta_{0,1})| \ge \min\left\{\left|\sin\left(\frac{2\pi}{2N+1}\right)\right|, \frac{3\sqrt{7}}{32}\right\}.\]
Then, for $k=1,...,r-1,\ell=0,1$, if $\rho_{k,\ell} \in \mathbb{R}$ (which by Corollary \ref{cor:onlyRealRoots} happens iff $r$ is even and $k=r/2$) then $\sin(\theta_{k,\ell})=0$. Otherwise, if $\rho_{k,\ell} \not \in \mathbb{R}$, Lemmas \ref{lm:lambda2} and \ref{lem:BasicRootsEqns} (which guarantee, respectively, that there exist constants $c_1$, $c_2$ such that for any conjugate roots $\rho_{k,\ell}, \overline{\rho_{k,\ell}}$ whose norms are not $1$, $|\rho_{k,\ell}-\overline{\rho_{k,\ell}}|\ge c_2c_1^{-r}\lambda_j^{1/2r}$, and that $|\rho_{k,\ell}| \le (1+\sqrt{2})^2$) imply that there exists an absolute constant $C_9>0$ such that 
\[|\sin(\theta_{k,\ell})| = \frac{|\rho_{k,\ell}- \overline{\rho_{k,\ell}}|}{2|\rho_{k,\ell}|} \geq C_9^{-r} \lambda_j^{1/2r}.\]
\end{proof}

Next, we seek to prove Lemma \ref{lm:C_lbound_large}. In the same vein as Lemma \ref{lm:tau} in the previous section, we will first bound the quotient $\frac{\Gamma_{k,\ell;k',\ell'}^{i} }{\Gamma_{k,\ell;k',\ell'}^{i+1}}$. To this end, we will prove the following lemma which bounds $\frac{\phi(i,\theta,\gamma)}{\phi(i+1,\theta,\gamma)}$, which will, in turn, allow us to prove bounds on $\frac{\Gamma_{k,\ell;k',\ell'}^{i} }{\Gamma_{k,\ell;k',\ell'}^{i+1}}$.

\begin{lm}\label{lm:phi_bd_large} If $\Delta N \ge 10$, then for all $i \in [m-1]$, where $m=\lfloor N/\Delta N \rfloor$ as defined in Definition \ref{def:intervals}), we have 
\[
\frac{1}{3 \Delta N^3}  \leq  \frac{ \phi(i,\theta,\gamma)}{ \phi(i+1,\theta,\gamma)}  \leq 3\Delta N^3,\]
where $\phi$ is as defined in (\ref{equ:def_phi}). 
\end{lm}

\begin{proof}

Let $\Delta N \ge 10$. 
We first recall from (\ref{equ:def_phi}) that
 $$
 \phi(i, \theta, \gamma): = \sum\limits_{n \in I_i} \cos^2(n\theta+\gamma).
 $$
We then break the argument into cases, depending on the distance from $\theta$ to an integer multiple of $\pi$, since the process for bounding $\phi(i,\theta,\gamma)$ will be different in each case.
\\[2\baselineskip]
\noindent \textbf{Case 1}: $\min\limits_{j\in \mathbb{Z}} | \theta- j \pi |  <\pi/\Delta N $. 
Suppose first that $\min\limits_{j\in \mathbb{Z}} | \theta- j \pi |  <\pi/\Delta N $. We will begin by proving by contradiction that
 \begin{align}\label{eq:ineq} 
\phi(i,\theta,\gamma)= \sum\limits_{n \in I_i} \cos^2(n\theta+\gamma) \geq \sin^2(\theta). 
 \end{align}
First, since $\Delta N \ge 10$, note that it suffices to show that the first three terms of \eqref{eq:ineq} satisfy

\begin{equation}\label{eq:intineq}\cos^2(A_i\theta+\gamma) + \cos^2((A_i+1)\theta+\gamma) + \cos^2((A_i+2)\theta+\gamma) \ge \sin^2(\theta).
\end{equation}

Let $\theta = \theta_0+k_\theta \pi$ where $k_\theta \in \mathbb{Z},\theta_0 \in [-\pi/\Delta N,\pi/\Delta N]$ (since $\min\limits_{j\in \mathbb{Z}} | \theta- j \pi |  <\pi/\Delta N $). Thus, since each of our functions of interest ($\cos^2, \sin^2, \cdots$) is $\pi$-periodic, it suffices to prove the result with $\theta$ replaced by $\theta_0$.
Now, assume the inequality \eqref{eq:intineq} (where $\theta_0$ is substituted for $\theta$) does not hold. 
Then it must be the case that $\cos^2(A_i\theta_0+\gamma),\cos^2((A_i+2)\theta_0+\gamma) < \sin^2(\theta_0)$, and thus we have
 \begin{equation}\label{eq:ineq1}
 |\sin(A_i \theta_0+\gamma+\pi/2)| =   |\cos(A_i \theta_0+\gamma)| < |\sin(\theta_0)|,
 \end{equation}
  and
   \begin{equation}\label{eq:ineq2}
 |\sin((A_i+2) \theta_0+\gamma+\pi/2)| =   |\cos((A_i+2) \theta_0+\gamma)| < |\sin(\theta_0)|.
 \end{equation} 
  From \eqref{eq:ineq1} we deduce that
   \begin{equation}\label{eq:intv1}A_i\theta_0+\gamma+\pi/2 \in \bigcup_{j\in \mathbb{Z}} (j\pi-\theta_0,j\pi+\theta_0). \end{equation}
Similarly, \eqref{eq:ineq2} implies that
 $$ (A_i+2)\theta_0+\gamma+\pi/2  \in \bigcup_{j\in \mathbb{Z}} (j\pi-\theta_0,j\pi+\theta_0), $$ or equivalently 
 \begin{equation}\label{eq:intv2}A_i\theta_0+\gamma+\pi/2 \in \bigcup_{j\in \mathbb{Z}} (j\pi-3\theta_0,j\pi-\theta_0). \end{equation}
which contradicts \eqref{eq:intv1}, since $\theta_0 \in [-\pi/\Delta N,\pi/\Delta N]$ and $\Delta N \ge10$ ensures that $|4\theta_0| < \pi$, and hence none of the intervals in \eqref{eq:intv1} and \eqref{eq:intv2} overlap.

Now that we have proved \eqref{eq:ineq} which shows $\phi$ has a lower bound,  we now seek to find an upper bound for $\phi$.
We first observe that 
\begin{align} \label{eq:trig_1}
 \phi(i,\theta,\gamma)&= \sum\limits_{n \in I_{i+1}}\cos^2(n\theta+\gamma-\Delta N \theta) \notag \\
 & =  \sum\limits_{n \in I_{i+1}}(\cos(n\theta+\gamma)\cos(\Delta N \theta)+\sin(\Delta N \theta)\sin(n\theta+\gamma))^2.
 \end{align}
 
Then, we note that 
\[(\cos(n\theta+\gamma)\cos(\Delta N \theta)+\sin(\Delta N \theta)\sin(n\theta+\gamma))^2
\le \cos^2(n\theta+\gamma)+|2\cos(n\theta+\gamma)\cos(\Delta N \theta)\sin(\Delta N \theta)\sin(n\theta+\gamma)|+\sin^2(\Delta N \theta),\]
and since 
\[|2\cos(n\theta+\gamma)\cos(\Delta N \theta)\sin(\Delta N \theta)\sin(n\theta+\gamma)| \le 2|\cos(n\theta+\gamma)\sin(\Delta N \theta)| \le \cos(n\theta+\gamma)^2+\sin^2(\Delta N \theta)\]
by the AM-GM Inequality, we thus have 
\[(\cos(n\theta+\gamma)\cos(\Delta N \theta)+\sin(\Delta N \theta)\sin(n\theta+\gamma))^2 \le 2\cos^2(n\theta+\gamma)+2\sin^2(\Delta N \theta),\]
and so
\begin{equation}\label{eq:trigamgm}
\sum\limits_{n \in I_{i+1}}(\cos(n\theta+\gamma)\cos(\Delta N \theta)+\sin(\Delta N \theta)\sin(n\theta+\gamma))^2 \le \sum\limits_{n \in I_{i+1}}2\cos^2(n\theta+\gamma)+2\sin^2(\Delta N \theta).
\end{equation}

Then, by the definition of $\phi$, the fact that $|I_{i+1}| = \Delta N$, and the fact that $\sin^2(\Delta N \theta) \le (\Delta N)^2 \sin^2(\theta)$ (which follows from repeated use of the identity $|\sin(a+b)| \le |\sin(a)|+|\sin(b)|$, which in turn holds since $|\sin(a+b)| \le |\sin(a)\cos(b)|+|\sin(b)\cos(a)| \le |\sin(a)|+|\sin(b)|$),
\begin{equation}\label{eq:trigjensen}
\sum\limits_{n \in I_{i+1}}2\cos^2(n\theta+\gamma)+2\sin^2(\Delta N \theta) \le 2(\Delta N)\sin^2(\Delta N \theta) + \sum\limits_{n \in I_{i+1}}2\cos^2(n\theta+\gamma) \le 2\phi(i+1,\theta,\gamma) + 2(\Delta N)^3 \sin^2(\theta)
\end{equation} 
Finally, combining \eqref{eq:trig_1}, \eqref{eq:trigamgm}, and \eqref{eq:trigjensen}, we see that 

\begin{equation}\label{eq:phi_ub}
\phi(i,\theta,\gamma) \le 2\phi(i+1,\theta,\gamma)+2(\Delta N)^3 \sin^2(\theta)
\end{equation}
By combining the lower and upper bounds [\eqref{eq:ineq} and \eqref{eq:phi_ub} respectively], we may now derive a bound on $\frac{ \phi(i,\theta,\gamma)}{ \phi(i+1,\theta,\gamma)} $, in particular,
 \begin{equation*}
\frac{ \phi(i,\theta,\gamma)}{ \phi(i+1,\theta,\gamma)} 
 \leq 2+\frac{2(\Delta N)^3\sin^2(\theta)}{\phi(i+1,\theta,\gamma)} \le 2+2(\Delta N)^3 \le 3 \Delta N^3.
  \end{equation*}
  where the last inequality used the fact that $(\Delta N)^3 \ge 2$.
Then, noting that 
\begin{align*} 
 \phi(i+1,\theta,\gamma)&= \sum\limits_{n \in I_{i}}\cos^2(n\theta+\gamma+\Delta N \theta) \notag \\
 & =  \sum\limits_{n \in I_{i+1}}(\cos(n\theta+\gamma)\cos(\Delta N \theta)-\sin(\Delta N \theta)\sin(n\theta+\gamma))^2
 \end{align*}
 and following the same arguments used to prove \eqref{eq:trigamgm}, and \eqref{eq:trigjensen}, we also have
 \begin{equation*}
\phi(i+1,\theta,\gamma) \le 2\phi(i,\theta,\gamma)+2(\Delta N)^3 \sin^2(\theta)
\end{equation*}
and hence (using the same lower bound \eqref{eq:ineq}) we see that
 \begin{equation}\label{eq:case1}
\frac{1}{3 \Delta N^3}  \leq  \frac{ \phi(i,\theta,\gamma)}{ \phi(i+1,\theta,\gamma)}  \leq 3\Delta N^3.
 \end{equation}
\noindent \textbf{Case 2}: $\min\limits_{j \in \mathbb{Z}} |\theta- j\pi| \geq \pi /\Delta N$
 
Recall from \eqref{eq:DeltaN} in the proof of Lemma \ref{lm:tau} that
\begin{equation} \label{eq:DeltaN2}
\phi(i,\theta,\gamma) = \frac{\Delta N}{2}+ \frac{\sin ((2A_{i+1}-1)\theta+2\gamma)-\sin((2A_{i}-1)\theta+2\gamma)}{4\sin \theta},
\end{equation}
and so, using sum-to-product identities (namely $\sin(\theta_1)-\sin(\theta_2) = 2 \sin(\frac{\theta_1-\theta_2}{2})\cos(\frac{\theta_1+\theta_2}{2})$), we see that 
\begin{align}\label{eq:sumtoprod} \frac{\Delta N}{2}+ \frac{\sin ((2A_{i+1}-1)\theta+2\gamma)-\sin((2A_{i}-1)\theta+2\gamma)}{4\sin \theta} &= \frac{\Delta N}{2} + \frac{\sin(\Delta N \theta)\cos((A_{i+1}+A_{i}-1)\theta+2\gamma)}{2\sin \theta} \notag \\
&= \frac{\Delta N}{2}\left(1 + \frac{\sin(\Delta N \theta)\cos((A_{i+1}+A_{i}-1)\theta+2\gamma)}{\Delta N \sin \theta}\right). 
\end{align}
We assume in this case that $\min\limits_{j \in \mathbb{Z}} |\theta- j\pi| \geq \pi /\Delta N$, which implies $| \sin(\theta) | \geq |\sin (\pi/\Delta N)| = \sin(\pi/\Delta N)$ (where the last equality holds because $\Delta N \ge 10$). Therefore, we have 
\[\left|\frac{\sin(\Delta N \theta)\cos((A_{i+1}+A_{i}-1)\theta+2\gamma)}{\Delta N \sin \theta}\right| \le \frac{1}{\Delta N \sin(\pi/\Delta N)}.\]

Continuing, we see that since $\frac{\sin x}{x} \ge \frac{1}{2}$ for $|x| \le 1$, and since $\Delta N \ge 10$, we thus have 
\[\frac{1}{\Delta N \sin(\pi/\Delta N)} \le \frac{2}{\pi}.\]
Combining this result with our earlier results \eqref{eq:DeltaN2} and \eqref{eq:sumtoprod}, we see that 
\begin{equation}\label{eq:lowerbd2}
\phi(i,\theta,\gamma) \ge \frac{\Delta N}{2}\left(1-\frac{2}{\pi}\right)
\end{equation}
We also note that clearly 
\begin{equation}\label{eq:upperbd2}
\phi(i,\theta,\gamma) = \sum\limits_{n \in I_i} \cos^2(n\theta+\gamma) \le \Delta N.
\end{equation}
Hence, by a similar argument to that used in the proof of Lemma \ref{lm:tau}, since our bounds \eqref{eq:lowerbd2} and \eqref{eq:upperbd2} on $\phi(i,\theta,\gamma)$ are independent of $i$, they also hold for $\phi(i+1,\theta,\gamma)$, and hence the bound for the ratio in this case is
 \begin{equation}\label{eq:case2}
\frac{1-2/\pi}{2} \leq \frac{ \phi(i,\theta,\gamma)}{ \phi(i+1,\theta,\gamma)} \leq \frac{2}{1-2/\pi}.
 \end{equation}
Note that $\frac{2}{1-2/\pi} < 6 < 3\Delta N^3$ since $\Delta N \ge 10$, and hence, we have \eqref{eq:case1} in all cases, which is our desired result.

\end{proof}
Using Lemma \ref{lm:phi_bd_large}, we now prove the following lemma (which is similar to Lemma \ref{lm:tau}) that bounds $\frac{\Gamma_{k,\ell;k',\ell'}^{i} }{\Gamma_{k,\ell;k',\ell'}^{i+1}}$. This lemma will then be used to bound $|D_{k,\ell;k',\ell'}|$, which ultimately allow us to prove our second main result, the bound on $|C|$ (Lemma \ref{lm:C_lbound_large}).

\begin{lm}\label{lm:tau_large} Suppose $\lambda_j^{1/2r} > 1/4$. There is a absolute constant $C_{10}$ such that if $\Delta N \geq C_{10} r^2 \log \Delta  N $, for all $(k,\ell) \ne (k',\ell')$ with $(k,\ell), (k',\ell') \in \{(0,0)\} \cup (\{1,\cdots, \lfloor r/2 \rfloor\} \times \{0,1\})$ either \[  \frac{\Gamma_{k,\ell;k',\ell'}^{i} }{\Gamma_{k,\ell;k',\ell'}^{i+1}}  \geq 4/3 \textrm{ for all $i \in [m-1]$,}  \textrm{ or  } \frac{\Gamma_{k,\ell;k',\ell'}^{i} }{\Gamma_{k,\ell;k',\ell'}^{i+1}} \leq 3/4 \textrm{ for all } i\in[m-1], \]
(for $m=\lfloor N/\Delta N \rfloor$ as defined in Definition \ref{def:intervals}).
\end{lm}
\begin{proof} 

As in Lemma \ref{lm:tau}, by the definition \eqref{equ:def_gamma} of $\Gamma_{k,\ell;k',\ell'}^{i}$, we have
  \begin{equation*}
\frac{\Gamma_{k,\ell;k',\ell'}^{i} }{\Gamma_{k,\ell;k',\ell'}^{i+1}} = \left(\frac{|\rho_{k,\ell}|}{|\rho_{k',\ell'}|}\right)^{-\Delta N}\frac{ \phi(i,\theta_{k,\ell},\gamma_{k,\ell})}{ \phi(i+1,\theta_{k,\ell},\gamma_{k,\ell})} \frac{ \phi(i+1,\theta_{k',\ell'},\gamma_{k',\ell'})}{ \phi(i,\theta_{k',\ell'},\gamma_{k',\ell'})} 
  \end{equation*}
  so, to bound this quantity, we need to consider bounds on $\left(\frac{|\rho_{k,\ell}|}{|\rho_{k',\ell'}|}\right)^{\Delta N}$.
 If $|\rho_{k,\ell}| > |\rho_{k',\ell'}|$, recall from \eqref{eq:abs_lambda} in Lemma \ref{lm:tau} that
 \[\left(\frac{|\rho_{k,\ell}|}{|\rho_{k',\ell'}|}\right)^{\Delta N} \ge \left(1 + \frac{C_4 \lambda_j^{1/2r}}{r^2}\right)^{\Delta N},\]
 for $C_4$ an absolute constant.
 Using the assumption $\lambda_j^{1/2r}\geq 1/4$, and the fact that $\log (1+x) \geq \frac{x}{1+x}$ for $x\geq -1$, we have
\begin{equation} \label{eq:rjrl}\left(\frac{|\rho_{k,\ell}|}{|\rho_{k',\ell'}|}\right)^{\Delta N} \ge \left(1 + \frac{C_4/4}{r^2}\right)^{\Delta N} \ge e^{\frac{C_4}{C_4+4r^2}\Delta N}.
\end{equation}
We now choose the absolute constant $C_{10}$ to be such that 
 \[C_{10}r^2 \log \Delta N \ge \max\left\{\frac{C_4+4r^2}{C_4}(6 \log \Delta N + \log 9 - \log(3/4)),10\right\},\]
 and suppose that $\Delta N \ge C_{10}r^2 \log \Delta N$.
We therefore have $\Delta N \ge 10$, so, by Lemma \ref{lm:phi_bd_large}, we have
 \[
\frac{1}{3 \Delta N^3}  \leq  \frac{ \phi(i,\theta,\gamma)}{ \phi(i+1,\theta,\gamma)}  \leq 3\Delta N^3.\] 
 Combining the above bound with \eqref{eq:rjrl},  we have
  \begin{equation*}
\frac{\Gamma_{k,\ell;k',\ell'}^{i} }{\Gamma_{k,\ell;k',\ell'}^{i+1}} = \left(\frac{|\rho_{k,\ell}|}{|\rho_{k',\ell'}|}\right)^{-\Delta N}\frac{ \phi(i,\theta_{k,\ell},\gamma_{k,\ell})}{ \phi(i+1,\theta_{k,\ell},\gamma_{k,\ell})} \frac{ \phi(i+1,\theta_{k',\ell'},\gamma_{k',\ell'})}{ \phi(i,\theta_{k',\ell'},\gamma_{k',\ell'})} \leq 9(\Delta N)^6 e^{-\frac{C_4}{C_4+4r^2}\Delta N} \leq 3/4,
  \end{equation*}
where the last inequality used the fact that $\Delta N \ge C_{10}r^2 \log \Delta N \ge \frac{C_4+4r^2}{C_4}(6 \log \Delta N + \log 9 - \log(3/4))$, which is equivalent to $e^{\Delta N} \ge [\frac{4}{3}(9(\Delta N)^6)]^{\frac{C_4+4r^2}{C_4}}$, i.e. $e^{\frac{C_4}{C_4+4r^2}\Delta N} \ge \frac{4}{3}[(9(\Delta N)^6)]$.
As in the proof of Lemma \ref{lm:tau}, in the case that $|\rho_{k,\ell}| < |\rho_{k',\ell'}|$, we can apply the same arguments to $\left(\frac{|\rho_{k',\ell'}|}{|\rho_{k,\ell}|}\right)^{\Delta N}$ to find a lower bound for $\left(\frac{|\rho_{k,\ell}|}{|\rho_{k',\ell'}|}\right)^{-\Delta N}$; combining this with the lower bounds on $\frac{ \phi(i,\theta_{k,\ell},\gamma_{k,\ell})}{ \phi(i+1,\theta_{k,\ell},\gamma_{k,\ell})} ,\frac{ \phi(i+1,\theta_{k',\ell'},\gamma_{k',\ell'})}{ \phi(i,\theta_{k',\ell'},\gamma_{k',\ell'})}$, we find that $\frac{\Gamma_{k,\ell;k',\ell'}^{i} }{\Gamma_{k,\ell;k',\ell'}^{i+1}} \geq 4/3$ in this case, and the desired result therefore holds.

\end{proof}
Then, as in the previous section, it is sufficient to bound $|D_{k,\ell;k',\ell'}|$ and use the fact that $|C^c_{k,\ell;k',\ell'}| \ge |D_{k,\ell;k',\ell'}|$ to bound $|C^c_{k,\ell;k',\ell'}|$ and then bound $|C|$. Thus, just as we used Lemma \ref{lm:tau} to prove Lemma \ref{lm:card_D}, we use Lemma \ref{lm:tau_large} to prove the following bound on $|D_{k,\ell;k',\ell'}|$. 
\begin{lm}\label{lm:card_D_large}
For any $(k,\ell)\neq (k',\ell')$, if
 \[  \frac{\Gamma_{k,\ell;k',\ell'}^{i} }{\Gamma_{k,\ell;k',\ell'}^{i+1}}  \geq 4/3 \textrm{ for all $i$,}  \textrm{ or  } \frac{\Gamma_{k,\ell;k',\ell'}^{i} }{\Gamma_{k,\ell;k',\ell'}^{i+1}} \leq 3/4 \textrm{ for all } i  \]
 then $|D_{k,\ell;k',\ell'}| \geq \lfloor N/\Delta N \rfloor - 32 \Delta N - 8\log 16 r^2 - 2.$
\end{lm}
\begin{proof} Without loss of generality, assume $\frac{\Gamma_{k,\ell;k',\ell'}^{i} }{\Gamma_{k,\ell;k',\ell'}^{i+1}}  \geq 4/3$ (just as in Lemma \ref{lm:card_D}, Lemma \ref{lm:interval} will imply the same result if  $\frac{\Gamma_{k,\ell;k',\ell'}^{i} }{\Gamma_{k,\ell;k',\ell'}^{i+1}}  \leq 3/4$). Also, as in the proof of Lemma \ref{lm:card_D}, there are slightly different arguments depending on whether or not $|\rho_{k,\ell}|, |\rho_{k',\ell'}|$ are greater than $1$ or not. 

Let $S= \frac{4 \Delta N + \log 16 r^2}{\log 2}$, so that $2^{S}\geq 4r^2 e^{4\Delta N}$.
As in Lemma \ref{lm:card_D}, by Lemma \ref{lm:interval} (with the same choices of $B_i, B'_i, \rho_1, \rho_2$, and $\alpha$, since the assumptions are the same), we know that except for $\lceil 2S/\log_2 (4/3) \rceil+1\ge 2S/\log_2 (4/3) + 2$ consecutive subsets, it holds that
\begin{equation}\label{eq:interval}
\Gamma_{k,\ell;k',\ell'}^{i} <2^{-S } \textrm{ or }  \Gamma_{k,\ell;k',\ell'}^{i} >2^{S }.
\end{equation}
Analogously to Lemma \ref{lm:card_D}, if $|\rho_{k,\ell}| \le 1, |\rho_{k',\ell'}|>1$ or $|\rho_{k,\ell}| > 1, |\rho_{k',\ell'}| \le 1$, the expressions in $D_{k,\ell;k',\ell'}$ are $\Gamma_{k,\ell;k',\ell'}^{i}$ or $\Gamma_{k,\ell;k',\ell'}^{i+1}$. As a result, since $2^S \geq 4r^2 e^{4\Delta N} \ge 4r^2$, we have the desired conclusion for all but $2S/\log_2(4/3)+2$ consecutive subsets.

Next, we note that in the case that $|\rho_{k,\ell}|, |\rho_{k',\ell'}| > 1$, we see that 
\begin{equation}\label{eq:rhoboundlarge}
|\rho_{k,\ell}|^{2\Delta N} \le (1+\sqrt{2})^{4\Delta N} \leq e^{4\Delta N},
\end{equation}
where the first inequality follows by Lemma \ref{lem:BasicRootsEqns} and the second inequality follows since $1+\sqrt{2}<e$. Now suppose that $\Gamma_{k,\ell;k',\ell'}^{i} \le 2^{-S}$. If this is the case, then we have
\begin{align}\label{eq:larger2}
\wtl c_{k',\ell'}^2 |\rho_{k',\ell'}|^{2A_i}  \phi (i, \theta_{k',\ell'}, \gamma_{k',\ell'}) & \ge 2 ^{S} \wtl c_{k,\ell}^2 |\rho_{k,\ell}|^{2A_i} \phi (i, \theta_{k,\ell}, \gamma_{k,\ell}) \notag\notag\\
&=2 ^{S} \wtl c_{k,\ell}^2 |\rho_{k,\ell}|^{2A_{i+1}} |\rho_{k,\ell}|^{-2\Delta N}  \phi (i, \theta_{k,\ell}, \gamma_{k,\ell}) \notag \\ &\geq 2^{S}e^{-4\Delta N } \wtl c_{k,\ell}^2 |\rho_{k,\ell}|^{2A_{i+1}} \phi (i, \theta_{k,\ell}, \gamma_{k,\ell}) \notag \\
& \ge 4r^2 \wtl c_{k,\ell}^2 |\rho_{k,\ell}|^{2A_{i+1}} \phi (i, \theta_{k,\ell}, \gamma_{k,\ell}).
\end{align}
where the first inequality is from the definition (\ref{equ:def_gamma}) of $\Gamma_{k,\ell;k',\ell'}^{i}$, and where the last equality is due to the definition of $S$. A similar argument using \eqref{eq:rhoboundlarge} on $|\rho_{k',\ell'}|$ shows that 
\[\wtl c_{k,\ell}^2 |\rho_{k,\ell}|^{2A_{i}} \phi (i, \theta_{k,\ell}, \gamma_{k,\ell}) \ge 4r^2 \wtl c_{k',\ell'}^2 |\rho_{k',\ell'}|^{2A_{i+1}}  \phi (i, \theta_{k',\ell'}, \gamma_{k',\ell'})  \]
if $ \Gamma_{k,\ell;k',\ell'}^{i} \ge 2^{S }$. 

Also, an analogous argument implies the same result if $|\rho_{k,\ell}|,|\rho_{k',\ell'}| \le 1$, since
\[|\rho_{k,\ell}|^{-2\Delta N} \le (1+\sqrt{2})^{4\Delta N} \leq e^{4\Delta N},\] 
where the first inequality follows by Lemma \ref{lem:BasicRootsEqns} and the second inequality follows since $1+\sqrt{2}<e$.
 
Finally, the number of the excluded intervals 
\begin{align*}
2S/\log_2 (4/3)+2 &= \frac{8 \Delta N + 2\log 16 r^2}{\log 2\log_2(4/3)}+2 \\
&= \frac{8 \Delta N + 2\log 16 r^2}{\log(4/3)}+2 \\
&\leq 32 \Delta N + 8\log 16 r^2 + 2
\end{align*}
which implies the conclusion of this lemma.
\end{proof}

Finally, we now use this bound on $|D_{k,\ell;k',\ell'}|$ to prove our second main result, which provides a bound on $|C|$ (note the similarities with Lemma \ref{lm:C_lbound}).
\begin{lm}\label{lm:C_lbound_large} Let $\Delta N = \lceil C_{10}^2r^4 \rceil$  where $C_{10}$ is the same as in Lemma \ref{lm:tau_large}, then $$|C| \geq  \lfloor N/\Delta N \rfloor - 2r(2r-1)(32 \Delta N + 8\log 16 r^2 + 2)  . $$ 
\end{lm}
\begin{proof}
Since $\Delta N = \lceil C_{10}^2r^4 \rceil$, it can be verified that this $\Delta N$ satisfies the assumption of Lemma \ref{lm:tau_large}, since $\Delta N \ge C_{10}^2r^4 \ge C_{10}r^2 \log(C_{10}^2r^4+1) \ge C_{10}r^2 \log\Delta N$ where the second inequality follows since $x \ge \log(x^2+1)$ for all positive $x$. 
Thus Lemma \ref{lm:tau_large} implies the bounds on $\frac{\Gamma_{k,\ell;k',\ell'}^{i}}{\Gamma_{k,\ell;k',\ell'}^{i+1}}$ used in the assumption of Lemma \ref{lm:card_D_large}. Then, by Lemma \ref{lm:card_D_large} we have a lower bound on the cardinality of $D_{k,\ell;k',\ell'}$, $|D_{k,\ell;k',\ell'}| \geq \lfloor N/\Delta N \rfloor - 32 \Delta N - 8\log 16 r^2 - 2$. Thus, using this bound we can proceed to compute,
\begin{align*}
|C| & =|\cap C_{k,\ell; k',\ell'}^c| \geq \lfloor N/\Delta N \rfloor -|\cup C_{k,\ell; k',\ell'}|  \\ & \geq \lfloor N/\Delta N \rfloor - \sum_{k,\ell,k',\ell', (k,\ell)\neq (k',\ell')}|C_{k,\ell;k',\ell'}|  \\ &\geq \lfloor N/\Delta N \rfloor - \sum_{k,\ell,k',\ell', (k,\ell)\neq (k',\ell')}(\lfloor N/\Delta N \rfloor- |D_{k,\ell;k',\ell'}|) \\ & \geq \lfloor N/\Delta N \rfloor - 2r(2r-1)(32 \Delta N + 8\log 16 r^2 + 2).  
\end{align*}
\end{proof}